\documentclass[10pt,leqno]{article}

\usepackage[francais]{babel}
\usepackage{amsmath,amssymb,amscd}
\usepackage{a4wide}
\usepackage[latin1]{inputenc}
\usepackage[all]{xy} 
\usepackage{hyperref}
\usepackage{chngcntr}

\title{Comptage de fibrés de Hitchin pour le groupe $\mathrm{SL}(n)$}
\author{Pierre-Henri Chaudouard}
\date{}

\newenvironment{paragr}[1][]{\refstepcounter{subsubsection} \noindent \textbf{\thesubsubsection . \ #1}}{\medskip}

\newenvironment{theoreme}{ \medskip\refstepcounter{theo}  \noindent\textbf{Th\'eor\`eme \thetheo}. ---\em}{\em \medskip}
\newenvironment{proposition}{\medskip\refstepcounter{theo}   \noindent\textbf{Proposition \thetheo}. ---\em}{\em\medskip}
\newenvironment{corollaire}{\medskip\refstepcounter{theo}  \noindent\textbf{Corollaire \thetheo}. ---\em}{\em\medskip}

\newenvironment{lemme}{\medskip\refstepcounter{theo}   \noindent\textbf{Lemme \thetheo}. ---\em}{\em\medskip}

\newenvironment{preuve}[1][]{\noindent \textbf{Démonstration.} #1 --- }{\hfill
  \ensuremath{\square} \medskip}

\newenvironment{remarque}{\medskip\refstepcounter{theo}  \noindent\textbf{Remarque \thetheo}. ---}{\medskip}

\DeclareMathOperator{\Fr}{\mathfrak{F}}
\DeclareMathOperator{\D}{D}
\DeclareMathOperator{\IC}{IC}
\DeclareMathOperator{\Perv}{Perv}
\DeclareMathOperator{\Amp}{Amp}
\DeclareMathOperator{\Dbc}{D_{\mathrm{c}}^{\mathrm{b}}}

\DeclareMathOperator{\GL}{GL}

\DeclareMathOperator{\reg}{reg}

\DeclareMathOperator{\vol}{vol}
\DeclareMathOperator{\rang}{rang}

\DeclareMathOperator{\ad}{ad}

\DeclareMathOperator{\el}{ell}

\DeclareMathOperator{\Gal}{Gal}
\DeclareMathOperator{\Hom}{Hom}

\DeclareMathOperator{\Id}{Id}

\DeclareMathOperator{\aff}{aff}

\DeclareMathOperator{\Res}{Res}

\DeclareMathOperator{\Ker}{Ker}

\DeclareMathOperator{\SL}{SL}

\DeclareMathOperator{\Ext}{Ext}
\DeclareMathOperator{\Exp}{Exp}
\DeclareMathOperator{\Log}{Log}
\DeclareMathOperator{\trace}{trace}

\newcommand{\ZZ}{\mathbb{Z}}

\newcommand{\SG}{\mathfrak{S}}
\newcommand{\NN}{\mathbb{N}}

\newcommand{\AAA}{\mathbb{A}}
\newcommand{\CC}{\mathbb{C}}

\newcommand{\QQ}{\mathbb{Q}}

\newcommand{\Fq}{\mathbb{F}_q}

\newcommand{\Fqd}{\mathbb{F}_{q^d}}
\newcommand{\Fqm}{\mathbb{F}_{q^m}}

\newcommand{\Ql}{\mathbb{Q}_{\ell}}
\newcommand{\Qlb}{\overline{\mathbb{Q}}_{\ell}}
\newcommand{\Qb}{\overline{\mathbb{Q}}}

\newcommand{\ga}{\gamma}

\newcommand{\zc}{\mathcal{Z}}
\newcommand{\oc}{\mathcal{O}}

\newcommand{\Sc}{\mathcal{S}}

\newcommand{\ec}{\mathcal{E}}
\newcommand{\gc}{\mathcal{G}}

\newcommand{\lc}{\mathcal{L}}

\newcommand{\fc}{\mathcal{F}}

\newcommand{\hc}{\mathcal{H}}
\newcommand{\nc}{\mathcal{N}}

\newcommand{\kc}{\mathcal{K}}
\newcommand{\Ic}{\mathcal{I}}

\newcommand{\ggo}{\mathfrak{g}}

\newcommand{\Lgo}{\mathfrak{L}}

\newcommand{\mgo}{\mathfrak{m}}

\newcommand{\hgo}{\mathfrak{h}}

\newcommand{\Xgo}{\mathfrak{X}}

\newcommand{\al}{\alpha}
\newcommand{\be}{\beta}

\newcommand{\om}{\omega}
\newcommand{\Om}{\Omega}

\newcommand{\La}{\Lambda}
\newcommand{\la}{\lambda}

\newcommand{\Ga}{\Gamma}
\newcommand{\back}{\backslash}

\newcommand{\eps}{\varepsilon}

\renewcommand{\leq}{\leqslant}
\renewcommand{\geq}{\geqslant}

\begin{document}
\counterwithin{equation}{subsubsection}
\maketitle

\begin{abstract}
  Let $C$ be a smooth projective curve of genus $g$ over a finite field $\mathbb{F}_q$ and let $D$ be a divisor on $C$ of degree $>2g-2$. We assume that the characteristic of $\mathbb{F}_q$ is sufficiently large. Let $n$ be an integer and let $\beta$ be a line bundle on $C$ of degree $e$, coprime to $n$.  We give a formula for  the  number of stable ($D$-twisted) Hitchin bundles over $C$ of rank $n$ and determinant $\beta$ in terms of the number of  stable  Hitchin bundles over $C'$ of rank $n/d$ and degree $e$  where $C'$  ranges over cyclic covers $C'$ of $C$ of degree $d$ dividing $n$. Using a work by Mozgovoy-O'Gorman, we derive a closed formula for the following invariants of the moduli space of  ($D$-twisted) Hitchin bundles over $C$ of rank $n$, trace $0$ and  determinant $\beta$: its number of points over finite extensions of $\mathbb{F}_q$ and its $\ell$-adic Poincaré polynomial and Euler-Poincaré characteristic. Our main tools are the fundamental lemma of automorphic induction and a support theorem for the relative cohomology of a local system on the Hitchin fibration for the group $\mathrm{GL}(n)$.
  \end{abstract}

\tableofcontents
\section{Introduction}

\subsection{Comptage de fibrés de Hitchin}

\begin{paragr}
  Soit $C$  une courbe projective,  lisse et  géométriquement connexe sur un corps fini $\Fq$ de cardinal $q$. Soit $D$ un diviseur sur $C$. On supposera dans tout l'article qu'on a $\deg(D)>2g-2$ où $g$ est le genre de la courbe $C$.   Un fibré de Hitchin est un couple $(\ec,\theta)$  formé d'un $\oc_C$-module $\ec$ localement libre de rang fini   et d'un endomorphisme tordu $\theta :\ec \rightarrow \ec \otimes_{\mathcal{O}_{C}}\mathcal{O}_{C}(D)$. Par définition, le rang, le déterminant et la trace  de $(\ec,\theta)$ sont respectivement le rang de $\ec$, le fibré en droites donné par la puissance extérieure maximale de $\ec$  et la trace de $\theta$ qui est un élément de $H^0(C,\oc_C(D))$.   Soit $n\geq 1$ un entier et $e\in \ZZ$. On suppose $e$ et $n$ premiers entre eux. On fixe un fibré en droites $\be$ de degré $e$ sur $C$. Soit $N_n^\be(C,D)$ l'espace de modules des fibrés de Hitchin $(\ec,\theta)$ de trace  nulle, de rang $n$ et de déterminant $\be$. Il est connu, depuis les travaux de Nitsure, \cite{Nitsure}, que  $N_n^\be(C,D)$ est une variété quasi-projective et lisse sur $\Fq$.
\end{paragr}

\begin{paragr} L'ensemble des $\Fq$-points $N_n^\be(C,D)(\Fq)$ s'identifie à l'ensemble des classes d'isomorphisme de fibrés de Hitchin stables (au sens du § \ref{S:stabilite}), de rang $n$, de trace nulle et de déterminant $\be$. Dans cet article, on donne une formule explicite pour le cardinal $|N_n^\be(C,D)(\Fq)|$ de l'ensemble $N_n^\be(C,D)(\Fq)$ lorsque, du moins, la caractéristique $\mathrm{car}(\Fq)$ de $\Fq$ est grande; précisément,  elle doit vérifier la condition:
\begin{align}
  \label{eq:intro-hyp-p}\mathrm{car}(\Fq) >\max(n, 3+n(2g-2)+n(n-1)\deg(D)).
  \end{align}
   Comme la cohomologie de $N_n^\be(C,D)$ est pure, la formule qu'on obtient (appliquée aux extensions finies de $\Fq$) nous permet de calculer des invariants topologiques de $N_n^\be(C,D)$, à savoir le polynôme de Poincaré et  la caractéristique d'Euler-Poincaré.
\end{paragr}

\begin{paragr} Soit $M_n^e(C,D)$, resp. $M_n^\be(C,D)$,  l'espace de modules des fibrés de Hitchin $(\ec,\theta)$  dont le rang  est $n$ et le degré du  déterminant est $e$, resp. et le déterminant est $\be$. Ce sont encore des variétés quasi-projectives et lisses sur $\Fq$. Puisqu'on suppose $\deg(D)>2g-2$ et $\mathrm{car}(\Fq)>n$, on a
    \begin{align*}
    |M_n^\be(C,D)(\Fq)|&=|H^0(C,\oc_C(D))|\cdot | N_n^\be(C,D)(\Fq)|\\
    &=q^{\deg(D)+1-g} | N_n^\be(C,D)(\Fq)|.
  \end{align*}
 Il s'agit donc de trouver une formule pour   $|M_n^\be(C,D)(\Fq)|$. En genre $g=0$, les fibrés en droites sont simplement classifiés par leur degré et on  a
  \begin{align*}
    |M_n^\be(C,D)(\Fq)|&= |M_n^e(C,D)(\Fq)| \ \ \ \ (g=0).
  \end{align*}
  Une formule explicite quoique assez compliquée pour $|M_n^e(C,D)(\Fq)|$ (sans hypothèse sur $\mathrm{car}(\Fq)$) a été obtenue par Mozgovoy et Schiffmann dans  \cite{Moz-Schiff}. Suivant une méthode due à Mellit \cite{Mellit-invent},  Mozgovoy et  O'Gorman en ont déduit une formule sensiblement plus simple, cf. \cite[théorème 1.1 et théorème 4.6]{MoOgor},  inspirée de travaux précurseurs de  Hausel-Rodriguez-Villegas \cite{HRV}. Plus précisément, Mozgovoy et  O'Gorman définissent un polynôme \og universel \fg{}
  \begin{align}\label{eq:intro-MO}
    \mathbf{H}_{g',n',p'}\in  \ZZ[x_1^{\pm1},\dots,x_{g'}^{\pm1},z]
  \end{align}
  qui ne dépend que d'entiers $g',n',p'$ et dont on rappelle la définition au § \ref{S:Log}. Pour tout $m\geq 1$, ils établissent  la formule
  \begin{align}
     |M_n^e(C,D)(\Fqm)| = \mathbf{H}_{g,n,p}(\la_1^m,\ldots,\la_g^m, q^m)
  \end{align}
  où  $p=\deg(D)-(2g-2)$ et le $2g$-uplet $(\la_1,\dots,\la_{g},q \la_{1}^{-1},\dots,q\la_{g}^{-1})\in (\Qlb^\times)^{2g}$ est constitué des valeurs propres du Frobenius $\Fr$ agissant sur le premier espace de cohomologie $\ell$-adique $H^1(C,\Qlb )$ de la courbe $C$. On a fixé ici un nombre premier $\ell\not=\mathrm{car}(\Fq)$. Le polynôme $\mathbf{H}_{g,n,p}$ est assez symétrique pour que l'évaluation ci-dessus ne dépende pas de la façon d'indexer ces valeurs propres.  Cela termine notre calcul en genre $g=0$. Dans la suite de l'introduction, on suppose donc  $g\geq 1$.
\end{paragr}

\begin{paragr}\label{S:intro-Picard} Soit $P$ le schéma de Picard de $C$ et, pour tout $i\in \ZZ$, soit $P^i$ la composante de $P$ formée des fibrés en droites de degré $i$. On identifie $\be$ à un élément de $P^e(\Fq)$.  Soit $m\geq 1$ et  $Q_m^1$ l'ensemble des caractères $\chi$ de $P(\Fqm)$ qui sont d'ordre divisant $n$ et tels que le noyau $\Ker(\chi)$ vérifie $\Ker(\chi)\cap P^1(\Fqm)\not=\emptyset$. Soit $\chi\in Q_m^1$ et $\al\in P^1(\Fqm)$ tel que $\chi(\al)=1$. Soit $AJ_{\al}:C\times_{\Fq} \Fqm \to P^0\times_{\Fq} \Fqm$ l'application obtenue par composition de l'application d'Abel-Jacobi $C\to P^1$ avec la tensorisation par $\al^{-1}$. Soit $\tilde C_\chi \to C\times_{\Fq}\Fqm$ le revêtement   galoisien de groupe $P^0(\Fqm)$ qu'on obtient en tirant en arrière par $AJ_{\al}$ le torseur de Lang sur $P^0\times_{\Fq}\Fqm$. En divisant par le sous-groupe $\Ker(\chi)\cap P^0(\Fqm) \subset P^0(\Fqm)$, on obtient un revêtement cyclique  $\rho_\chi: C_\chi \to C\times_{\Fq}\Fqm$ de groupe de Galois noté $\Ga_\chi\simeq   P(\Fqm)/ \Ker(\chi)$. Soit $D_\chi=\rho_\chi^*(D\times_{\Fq} \Fqm)$. Ces constructions ne dépendent pas du  choix de $\al$, à $\Fqm$-isomorphisme près. Soit $d_\chi$ l'ordre de $\chi$ et  $n_\chi=n/d_\chi$. On peut alors formuler le résultat clef de cet article;  afin de simplifier quelque peu la formule  obtenue, nous supposerons qu'on a  $D=2D_1$  \emph{uniquement dans cette introduction}.

  \begin{theoreme} \label{thm:intro-endo} (Cas $D=2D_1$; pour un résultat plus général, cf. théorème \ref{thm:comparaison comptage})
On suppose \eqref{eq:intro-hyp-p}.    Pour tout entier $m\geq 1$, on a
    \begin{align}\label{eq:intro-endo}
      | M_n^\be(C,D)(\Fqm)|=\frac{1}{|P^0(\Fqm)|} \sum_{\chi\in Q_m^1} \frac{q^{m r_\chi\deg(D)}}{d_\chi} \chi(\be)^{-1} |M^e_{n_\chi}(C_\chi,D_\chi)(\Fqm)|
    \end{align}
  où   $r_\chi =n_\chi^2 d_\chi(d_\chi-1)/2=n(n-n_\chi)/2$.
  \end{theoreme}

  On renvoie à la seconde partie de l'introduction pour une explication sur l'origine de la restriction sur la caractéristique qui est due à la méthode de la preuve. On conjecture que la relation  \eqref{eq:intro-endo} vaut encore sans cette restriction (ou du moins sous la restriction plus faible $\mathrm{car}(\Fq)>n$). Lorsqu'on prend $D$ le diviseur canonique, on conjecture que la relation \eqref{eq:intro-endo} devrait être encore correcte.
  \end{paragr}

  \begin{paragr} Le théorème \ref{thm:intro-endo} permet d'exprimer le comptage en terme des polynômes de  Mozgovoy et  O'Gorman. Formulons pour l'introduction notre résultat ainsi:
    
\begin{theoreme} \label{thm:intro-comptage} (Cas $D=2D_1$.) On suppose \eqref{eq:intro-hyp-p}.    Pour tout entier $m\geq 1$, on a
    \begin{align}\label{eq:intro-compt}
      | N_n^\be(C,D)(\Fqm)|=\frac{q^{-m(\deg(D)+1-g))}}{|P^0(\Fqm)|} \sum_{\chi\in Q_m^1} \frac{q^{m r_\chi\deg(D)}}{d_\chi} \chi(\be)^{-1}  \mathbf{H}_{g_\chi,n_\chi,p_\chi}(\la_{1,\chi},\ldots,\la_{g_\chi,\chi}, q^m) 
    \end{align}
    où   l'on introduit, en sus de celles  du théorème \ref{thm:intro-endo}, les notations suivantes:
    \begin{itemize}
    \item $g_\chi=d_\chi(g-1)+1$ est le genre de la courbe $C_\chi$;
    \item $ p_\chi= d_\chi(\deg(D)-(2g-2))=\deg(D_\chi)-(2g_\chi-2)$; 
    \item le $2g_\chi$-uplet   $(\la_{1,\chi},\dots,\la_{g_\chi,\chi},q^m \la_{1,\chi}^{-1},\dots,q^m\la_{g_\chi,\chi}^{-1})\in (\Qlb^\times)^{2g_\chi}$  est formé des valeurs propres de $\Fr^m$ agissant sur $H^1(C_\chi,\Qlb)$.
    \end{itemize}
  \end{theoreme}

  Notons qu'on peut analyser un peu plus la formule obtenue: tout d'abord, on peut regrouper les éléments de $Q_m^1$ selon leur restriction à $P^0(\Fqm)$. Celles-ci forment le groupe des caractères d'ordre divisant $n$ de $P^0(\Fqm)$ et ce dernier, via la norme, ne dépend en fait pas de $m$ (il y a une condition à vérifier sur $q$, cf. § \ref{S:n-torsion}). On peut alors dans la formule \eqref{eq:intro-compt} ne faire intervenir qu'un nombre fini, indépendant de $m$, de revêtements cycliques de la courbe $C$ d'ordre divisant $n$. On obtient alors une formule dans l'esprit de celle de Mozgovoy-O'Gorman: le lecteur trouvera sa formulation dans le théorème \ref{thm:comptage tr 0}, qui s'affranchit d'ailleurs de l'hypothèse $D=2D_1$.
  \end{paragr}

  \begin{paragr}
    Il est bien connu que la cohomologie de la variété $M_n^e(C,D)$ est pure: il s'ensuit que le comptage de Mozgovoy-O'Gorman donne le polynôme de Poincaré de la cohomologie $\ell$-adique à support compact de $M_n^e(C,D)$. C'est le polynôme en la variable $u$ donné par la spécialisation du polynôme   $\mathbf{H}_{g,n,p}$, défini en \eqref{eq:intro-MO} pour $p=\deg(D)-(2g-2)$, en $x_1=\cdots=x_g=-u$ et $z=u^2$. Soit $d|n$ et $\tilde e\in \ZZ$. On définit en § \ref{S:specialisation} un polynôme  $\tilde{\mathbf{H}}_{g,n,p,d,\tilde e}^\flat\in \ZZ[ u]$ en divisant par $(1-u)^{2g} u^{2\deg(D)-(2g-2)}$ des spécialisations légèrement différentes du polynôme  $\mathbf{H}_{g',n',dp}$ avec $n'd=n$ et $g'-1=d(g-1)$. On obtient alors la formule suivante (qui vaut sans l'hypothèse $D=2D_1$):

    \begin{theoreme}\label{thm:intro-Poincare} (cf. corollaire \ref{cor:Poincare}) On suppose \eqref{eq:intro-hyp-p}.   Le polynôme de Poincaré de la cohomologie $\ell$-adique à support compact de $N_n^\be(C,D)$ est donné par
  \begin{align*}
  \sum_{d|n} \psi_g(d) u^{   \deg(D)(d-1)n^2/d} \tilde{\mathbf{H}}^\flat_{g,n,p,d,\tilde e}(-u)
  \end{align*}
  où
  \begin{itemize}
  \item $\tilde e=e+\frac{n(n-1)}{2}\deg(D)$;
  \item  $\psi_g(d)$ est le nombre d'éléments de $(\ZZ/d\ZZ)^{2g}$ d'ordre $d$.
  \end{itemize}
\end{theoreme}

\begin{remarque}
  On voit que le polynôme de  Poincaré  dépend uniquement du genre $g$ de $C$, du rang $n$, du degré de $D$ et du degré $e$ de $\be$. Il ne dépend pas du choix de $D$ et de $\be$ quand leur degré est fixé. 
\end{remarque}

On peut aussi calculer explicitement la caractéristique d'Euler-Poincaré:

\begin{corollaire}\label{cor:intro-euler} (pour un énoncé plus général si $\deg(D)$ est impair, cf. corollaire \ref{cor:Euler}) Supposons \eqref{eq:intro-hyp-p} et $\deg(D)$ pair. La caractéristique d'Euler-Poincaré de $N_n^\be(C,D)$ vaut:
  \begin{itemize}
  \item $1$ si $g=1$;
    \item  $\mu(n) n^{4g-3}$  si $g\geq 2$.
    \end{itemize}
  \end{corollaire}

  \begin{remarque}
    Il est bien connu que les formules  obtenues dans le théorème \ref{thm:intro-Poincare} et le corollaire \ref{cor:intro-euler} (ou celles plus générales des corollaires \ref{cor:Poincare} et \ref{cor:Euler}) qui valent en caractéristique assez grande valent encore pour la variété $N_n^\be(C,D)$ et sa cohomologie de Betti à coefficients dans $\CC$ lorsque la courbe $C$ est définie sur le corps $\CC$ des nombres complexes. Mentionnons que ces calculs sur $\CC$  devraient aussi résulter de la théorie de l'endoscopie appliquée à la cohomologie de $N_n^\be(C,D)$, par exemple des résultats de  Maulik et Shen dans \cite{Endoscopic-Maulik-Shen} qui traitent également du cas du diviseur canonique. Comme on va le voir dans la sous-section \ref{ssec:a propos}, nos résultats passent par une étude de la cohomologie relative à valeurs dans un système local de la fibration de Hitchin pour le groupe $\GL(n)$ et évite toute analyse de celle du groupe $\SL(n)$. 
  \end{remarque}
\end{paragr}

\subsection{À propos de la preuve}\label{ssec:a propos}

\begin{paragr} Le théorème \ref{thm:intro-endo} se déduit, via la formule des traces de Grothendieck-Lefschetz, d'un énoncé cohomologique. Soit $\chi\in Q^1_1$ et $\al\in P^1(\Fq)$ tel que $\chi(\al)=1$, cf. § \ref{S:intro-Picard}. Soit $\gc_\chi$ le $\Qlb$-faisceau lisse de rang $1$ sur  $P^0$ qu'on obtient en poussant $\chi^{-1}$ par le torseur de Lang, cf. §  \ref{S:chi}. Soit  $\gc_{\chi,\al}$ le système local sur $P^e$ obtenu lorsqu'on tire en arrière $\gc_\chi$ par le morphisme $P^e\to P^0$ donné par la tensorisation par $\al^{-e}$.   On dispose alors d'un morphisme $M_n^e(C,D)\to P^e$ qui, à un fibré de Hitchin $(\ec,\theta)$, associe  son déterminant. En tirant en arrière  $\gc_{\chi,\al}$ par ce morphisme, on obtient un système local noté $\lc_\chi$. Soit $f:  M_n^e(C,D)\to A$ le morphisme de Hitchin qui, à  un fibré de Hitchin $(\ec,\theta)$, associe le polynôme caractéristique de $\theta$. C'est un morphisme projectif dont la base $A=\oplus_{i=1}^n H^0(C,\oc_C(iD))$ est un espace affine.

  Soit $\rho_\chi: C_\chi \to C$ le revêtement cyclique de groupe de Galois $\Ga_\chi\simeq P(\Fq)/\Ker(\chi)$ et  $D_\chi$ le diviseur sur $C_\chi$ qu'on a définis au § \ref{S:intro-Picard} pour $m=1$. L'ordre de $\Ga_\chi$, qui est aussi l'ordre du caractère $\chi$,  est un diviseur de $n$ noté $d$:  on écrit $n=n'd$.  Comme ci-dessus mais cette fois-ci relativement à la courbe $C_\chi$, au rang $n'$, au degré $e$ et au diviseur $D_\chi$, on dispose d'un morphisme propre $f': M_{n'}^e(C_\chi,D_\chi)\to A'$. On dispose aussi d'un morphisme fini $\iota:A'\to A$ entre les bases, cf. § \ref{S:morph}. Le groupe   $\Ga_\chi$ agit naturellement sur l'espace total et la base de la fibration $f'$. De plus, les morphismes $f'$ et  $\iota$  sont respectivement  équivariant et invariant pour les actions de $\Ga_\chi$. Par composition, le morphisme $\iota\circ f'$ est également  propre et $\Ga_\chi$-invariant: il se descend en un morphisme propre noté  $\bar f': M_{n'}^e(C_\chi,D_\chi) / \Ga_\chi\to A$.   En fait, le groupe $\Ga_\chi$ agit sans point fixe sur la variété  $M_{n'}^e(C_\chi,D_\chi)$ et le quotient    $M_{n'}^e(C_\chi,D_\chi)/\Ga_\chi$ est lisse. On le munit du système local  noté $\mu_{\chi^e}$  associé au caractère  $\chi^e$ de $\Ga_\chi$, cf.  § \ref{S:mu chi}.

  \begin{theoreme}\label{thm:intro-iso-fx} (ici $D=2D_1$, pour un énoncé plus général cf. théorème \ref{thm:iso-fx})

    On suppose \eqref{eq:intro-hyp-p}. Il existe un isomorphisme entre les  semi-simplifications des faisceaux pervers gradués
    \begin{align*}
 \oplus_{i\in \ZZ}   {}^{\mathrm{p}}\mathcal{H}^{i}(Rf_{*} \lc_\chi )[2r](r)      \text{  et  }   \oplus_{i\in \ZZ}   {}^{\mathrm{p}}\mathcal{H}^{i}(R\bar f'_{*} \mu_{\chi^{e}} )
      \end{align*}
      avec  $r= \frac{n(n-n')}{2}\deg(D)$.
    \end{theoreme}

    On a noté ${}^{\mathrm{p}}\mathcal{H}^{i}$ le $i$-ème faisceau pervers de cohomologie. Les complexes $Rf_{*} \lc_\chi $ et $R\bar f'_{*} \mu_{\chi^{e}}$ sont purs et, de ce fait, somme de leurs faisceaux pervers de cohomologie à des décalages près (pour des rappels, cf. § \ref{S:def support}). Le théorème \ref{thm:intro-endo} se déduit alors du théorème \ref{thm:intro-iso-fx} par une application de  la formule des traces de Grothendieck-Lefschetz et une simple inversion de Fourier.
    \end{paragr}

    \begin{paragr}[Théorème de support.] --- Il résulte de \cite{Ast100} que les faisceaux pervers du théorème \ref{thm:intro-iso-fx}, après changement de base à une clôture algébrique de $\Fq$, sont tous semi-simples. Une information souvent cruciale est la détermination des supports des constituants simples. En caractéristique nulle, nous  avions montré avec Laumon dans \cite{Support} (qui faisait suite aux travaux de Ngô \cite{Ngo-LF} et à notre travail \cite{LFPII} avec Laumon) que le seul support de l'image directe du faisceau constant par le morphisme de Hitchin était la base toute entière. Il n'en est pas de même ici pour $Rf_{*} \lc_\chi$: on montre en effet que tous les supports sont inclus dans le fermé, noté $\iota(A')$, qui est l'image du morphisme $\iota$, cf. proposition \ref{prop:borne-supp}. Néanmoins, on montre que $\iota(A')$ est le seul support pour les constituants des faisceaux pervers  du théorème \ref{thm:intro-iso-fx}, cf.  théorèmes \ref{thm:support} et \ref{thm:support2}. L'argument, qui  apparaît déjà pour l'essentiel dans les travaux précédemment mentionnés,  repose de manière cruciale sur une inégalité de Severi, cf. \eqref{eq:Severi}. On vérifie celle-ci dans la sous-section  \ref{ssec:delta cst} lorsque $\mathrm{car}(\Fq)$ vérifie \eqref{eq:intro-hyp-p}. Dès lors, les théorèmes de support ramènent la démonstration du théorème \ref{thm:intro-iso-fx} à sa vérification sur un ouvert dense de $\iota(A')$. On prend l'ouvert dit elliptique de  $\iota(A')$ sur lequel on peut utiliser une variante du lemme fondamental démontré par Ngô.
    \end{paragr}

    \begin{paragr}[Remarques.] --- Lorsqu'on prend pour $D$ le diviseur canonique, on ne s'attend plus à ce que la description des supports soit aussi simple, cf. \cite{DeCataldoHM} pour une étude de l'image directe du faisceau constant dans le cas d'une courbe sur le corps des nombres  complexes. Néanmoins, il devrait être possible d'établir un isomorphisme comme dans le  théorème \ref{thm:intro-iso-fx}. Une possibilité est de réaliser l'espace de Hitchin $M_n^e(C,D-[x])$ comme ensemble critique d'une certaine fonction sur $M_n^e(C,D)$ et d'utiliser le foncteur des cycles évanescents, cf. \cite[section 4]{Endoscopic-Maulik-Shen} pour une telle approche sur $\CC$ et le faisceau constant sur $N_n^\be(C,D)$.

      Signalons encore deux autres approches possibles  au théorème \ref{thm:intro-endo} qui devraient pouvoir s'affranchir de l'hypothèse que $\mathrm{car}(\Fq)$ est grand devant $\deg(D)$ et valoir encore pour le diviseur canonique. Une première approche est d'utiliser notre travail \cite{scfh} dans lequel on fait le lien entre le comptage des points sur les corps finis de $M_n^e(C,D)$ et la formule des traces pour les algèbres de Lie sur le corps de fonctions de $C$ associée à une fonction test très simple liée au diviseur $D$. La théorie de l'induction automorphe suggère que cette  formule des traces pour les algèbres de Lie, convenablement  tordue par le caractère $\chi$,  devrait être égale à la même formule (non tordue) associée à  $M_{n'}^e(C_\chi,D_\chi)$: c'est le pendant du théorème \ref{thm:intro-iso-fx}. Le lemme fondamental qu'on utilise permet de réaliser cette comparaison pour tous les termes associés à des orbites régulières semi-simples elliptiques. Des méthodes automorphes devraient permettre d'étendre cette comparaison à toutes les orbites. Une seconde approche est d'utiliser les méthodes d'intégration $p$-adiques développées dans \cite{GWZ,GWZ2} par Groechenig,  Wyss et  Ziegler  soit via l'endoscopie du groupe $\SL(n)$ soit  peut-être, ce qui serait plus direct, par l'introduction d'intégrales pour le groupe $\GL(n)$  tordues par le caractère $\chi$.  
    \end{paragr}

    \subsection{Organisation de l'article}

    \begin{paragr}
      Dans la section \ref{sec:cbe-spec}, on introduit les courbes spectrales et l'inégalité de Severi. Sous une hypothèse sur la caractéristique du corps de base, on démontre cette inégalité dans la sous-section \ref{ssec:delta cst}. La section \ref{sec:Hitchin} est consacrée à des rappels sur la fibration de Hitchin. La section  \ref{sec:coh} a pour objectif d'énoncer et de démontrer les théorèmes de support mentionnées plus haut. Ceux-ci sont un ingrédient essentiel de la preuve du théorème \ref{thm:intro-iso-fx} de comparaison des cohomologies relatives, preuve à laquelle se consacre  la section \ref{sec:compa}. L'autre point essentiel, à savoir une version du lemme fondamental, est introduit dans la sous-section \ref{ssec:local}. La section \ref{sec:applications} finale donne des applications du théorème \ref{thm:intro-iso-fx}. On introduit tout d'abord dans la sous-section \ref{ssec:prelim} les différents polynômes qui vont intervenir. Dans la sous-section \ref{ssec:det fixe}, on exploite le  théorème \ref{thm:intro-iso-fx}  pour établir le théorème \ref{thm:comptage det fix} qui donne une formule explicite pour le comptage des fibrés de Hitchin de déterminant fixé. On en déduit dans la sous-section \ref{ssec:trace nulle}, une formule  explicite pour le comptage des fibrés de Hitchin de déterminant fixé et de trace nulle, cf. théorème \ref{thm:comptage tr 0}. On obtient  le calcul d'invariants topologiques dans la sous-section \ref{ssec:Poincare}.
    \end{paragr}

    \begin{paragr}[Remerciements.] --- Je remercie l'Institut universitaire de France (IUF) pour m'avoir offert d'excellentes conditions de travail durant la finalisation de cet article.
          \end{paragr}

    \section{Courbes spectrales}\label{sec:cbe-spec}

\subsection{Définition}\label{ssec:def-cbe-spec}

\begin{paragr}
    Soit $C$ une courbe projective,  lisse et  géométriquement connexe sur un corps $k$. Soit $g_C$ son  genre. Soit $D=\sum_{x}D_{x}[x]$ un diviseur sur $C$. 
  \end{paragr}

  \begin{paragr}
  Soit $n\geq 1$ un entier. On introduit le schéma affine 
  \begin{align}\label{eq:Ac}
    A_{n}=\bigoplus_{i=1}^{n}H^{0}(C,\mathcal{O}_{C}(iD)).
  \end{align}
  Lorsque le contexte est clair, on omet l'indice $n$ et on note simplement $A$ ce schéma. 
  Lorsqu'on a
  \begin{align}\label{eq:degD}
    \deg(D)>2g_C-2,
  \end{align}
condition qu'on supposera dans les principaux résultats de cet article, la dimension de $A$ est donnée par la formule de Riemann-Roch:
  \begin{align*}
    d_{A}=n(1-g_C)+\frac{n(n+1)}{2}\deg(D).
  \end{align*}

\end{paragr}

\begin{paragr}\label{S:esp-total}
  Soit
  \begin{align}\label{eq:esp-total}
    p:\Sigma=\mathbb{V}(\mathcal{O}_{C}(-D))\rightarrow C
  \end{align}
  l'espace total du fibré en droites $\mathcal{O}_{C}(D)$. Soit $u$ la section universelle de $p^{\ast}\mathcal{O}_{C}(D)$. Pour tout $a\in A$ et $1\leq i\leq n$, on dispose d'une section $a_i\in H^0(C,\oc_C(iD))$ et de  la section globale
\begin{align}\label{eq:a(u)}
  a(u)=u^{n}+p^{\ast}(a_{1})u^{n-1}+\cdots+p^{\ast}(a_{n})
\end{align}
de $p^{\ast}\mathcal{O}_{C}(nD))$. La courbe spectrale universelle $X \to  A$ est le diviseur de Cartier relatif dans $A \times_{k}\Sigma /A$ défini par l'équation $a(u)=0$. Pour $a\in A(k)$, on obtient un diagramme commutatif:
  \begin{align}\label{eq:diag}
    \xymatrix{ X_a  \ar@{^{(}->}[r]^{\rho_a}    \ar[rd]_{\pi_a} & \Sigma \ar[d]^{p}  \\      & C}
  \end{align}
  où la courbe $\rho_a : X_a\hookrightarrow \Sigma$ est définie par  le faisceau d'idéaux $\Ic_a$ qui est l'image de l'homomorphisme $p^*\oc_C(-nD)\to \oc_\Sigma$ induit par la section globale $a(u)$ et $\pi_{a}:X_{a}\rightarrow C$ est la  restriction  de $p$ à $X_{a}$.  La courbe $X_{a}$ est connexe  mais elle n'est ni réduite ni irréductible en général. La restriction $\pi_{a}:X_{a}\rightarrow C$ de $p$ à $X_{a}$ est un morphisme fini et plat de degré $n$, et on a
$$
\pi_{a,\ast}\mathcal{O}_{X_{a}}=\mathcal{O}_{C}\oplus\mathcal{O}_{C}(-D)\oplus\cdots\oplus\mathcal{O}_{C}(-(n-1)D),
$$
En particulier le genre arithmétique de $X_a$ est donné par
\begin{align}
\label{eq:qa}  q_a&=1-\chi(X_a, \mathcal{O}_{X_{a}})\\
\nonumber  &= 1+n(g_C-1)+\frac{n(n-1)}{2}\deg(D).
\end{align}
On utilisera parfois la formule:
\begin{align}
  \label{eq:dA-qa}
  d_A-q_a=n(\deg(D)-2g_C +2)-1
\end{align}
\end{paragr}

\begin{paragr}[Lieu elliptique.] --- \label{S:ell}  Soit  $r\geq 1$ un entier, $\underline m=(m_i)_{1\leq i\leq  r}$ et $\underline n=(n_i)_{1\leq i\leq  r}$ deux familles d'entiers $\geq 1$ telles que
  \begin{align}
    \label{eq:summi}
    \sum_{i=1}^r m_in_i=n.
  \end{align}
 Soit
  \begin{align*}
    \iota_{\underline n,\underline m}: \prod_{i=1}^r A_{n_i}\to A=A_n
  \end{align*}
  où l'image $a$ de $(a_1,\ldots,a_r)$ est donnée par
  \begin{align*}
    a(u)=\prod_{i=1}^r a_i^{m_i}(u),
  \end{align*}
  cf. \eqref{eq:a(u)}.  Le morphisme $\iota_{\underline n,\underline m}$ est fini. On note $A_{\underline n,\underline m}$ son image. C'est un fermé de $A$. Soit $A^{\el}$ l'ouvert complémentaire de la réunion des fermés $A_{\underline n,\underline m}$ lorsque $\underline n$ et $\underline m$ parcourent tous les $r$-uplets  possibles vérifiant \eqref{eq:summi}  exceptés $(\underline n,\underline m)=((n),(1))$ pour $r=1$. C'est encore le lieu des $a\in A$ pour lesquels la courbe spectrale $X_a$ est géométriquement intègre. Pour des raisons de dimension, l'ouvert $A^{\el}$ est non vide.
\end{paragr}

\begin{paragr} \label{S:deformation}  La théorie des déformations du diagramme \eqref{eq:diag}  est contrôlée par le complexe
  \begin{align*}
    R\Hom_{\oc_{X_a}}(  L_{X_a/\Sigma_D},\oc_{X_a})
  \end{align*}
où l'on introduit le complexe   cotangent:
  \begin{align*}
    L_{X_a/\Sigma_D}=[\rho_a^* \Om^1_{\Sigma/C} \to \Om^1_{X_a/C} ],
  \end{align*}
  complexe placé en degré $-1$ et $0$. Ce dernier s'identifie au faisceau  localement libre $\Ic_a/\Ic_a^2$ placé en degré $-1$. Par ailleurs, on a
  \begin{align}\label{eq:I-I2}
    \Ic_a/\Ic_a^2=\rho_a^*\Ic_a\simeq \rho_a^*p^* \oc_C(-nD)=\pi_a^*\oc_C(-nD).
  \end{align}
   En utilisant \eqref{eq:I-I2}, on voit que, pour $i=0,1$, on a 
  \begin{align*}
    \Ext^{i+1}_{\oc_{X_a}}(L_{X_a/\Sigma_D},   \oc_{X_a}) & =\Ext^i_{\oc_{X_a}}( \Ic_a/\Ic_a^2      ,\oc_{X_a})  \\
                                                                                                                    &=\oplus_{i=1}^n  \Ext^i_{\oc_{C}}( \oc_C     ,   \oc_C(iD)).
  \end{align*}
Sous l'hypothèse  \eqref{eq:degD}, on a $\Ext^{2}_{\oc_{X_a}}(L_{X_a/\Sigma_D},   \oc_{X_a})=0$: il n'y a donc  pas d'obstruction à déformer \eqref{eq:diag} et  l'espace $A$ est lisse  d'espace tangent en $a$
  \begin{align*}
     \Ext^{1}_{\oc_{X_a}}(L_{X_a/\Sigma_D}   ,\oc_{X_a})
  \end{align*}
  de dimension $d_{A}$.
\end{paragr}

\subsection{Strate à $\delta$ constant}\label{ssec:delta cst}

\begin{paragr}
  On continue avec les notations précédentes. Soit $p$ la caractéristique du corps $k$ si celle-ci est non-nulle. Si la caractéristique de $k$ est nulle, on pose $p=+\infty$ afin d'avoir des énoncés uniformes. On suppose qu'on a
  \begin{align}\label{eq:hyp-p}
    (p-1)/2> 1+n(g_C-1)+\frac{n(n-1)}{2}\deg(D)
  \end{align}
  c'est-à-dire
    \begin{align*}
      p> 3+n(2g_C-2)+n(n-1)\deg(D).
    \end{align*}
    Rappelons que le membre de droite dans \eqref{eq:hyp-p} est le genre arithmétique $q_a$ commun aux courbes spectrales $X_a$ pour $a\in A$ et calculé en \eqref{eq:qa}.

    \begin{proposition}\label{prop:lissite}
    Sous l'hypothèse  \eqref{eq:hyp-p}, pour tout $a\in A^{\el}$, la normalisée $\tilde X_a$ de $X_a$ est lisse.
    \end{proposition}

    \begin{preuve}   Il s'agit de prouver que $\tilde X_a$ est géométriquement régulière. Il suffit pour cela que la courbe $\tilde X_a$ soit géométriquement normale. Soit  $K'/K$ une extension du corps résiduel $K$ de $a$. Il s'agit de montrer que $Y=\tilde X_a\times_K K'$ est normale.  On introduit $\nu:\tilde Y\to Y$ la normalisation de $Y$. Notons que $X_a$ est géométriquement intègre vu que $a\in A^{\el}$. Il en est de même de $\tilde X_a$ et $\tilde Y$ \cite[lemme 1.1]{Schroer}. D'après un théorème de Tate revisité par Schröer, voir \cite{Tate} et \cite[corollaire 2.3]{Schroer}, la différence
      \begin{align}\label{eq:diff}
        \dim_K(H^1(\tilde X_a,\oc_{\tilde X_a}))-\dim_{K'}(H^1(\tilde Y, \oc_{\tilde Y})),
      \end{align}
      est divisible par $(p-1)/2$.  
      En utilisant la suite exacte longue en cohomologie associée à la suite courte
      \begin{align*}
       0 \to \oc_Y \to \nu_* \oc_{\tilde Y}\to   \nu_* \oc_{\tilde Y}/\oc_Y \to 0
      \end{align*}
      et le fait que  $\nu_* \oc_{\tilde Y}/\oc_Y $ est un faisceau de torsion, on obtient la suite exacte
      \begin{align}\label{eq:suite-Hi}
         0 \to H^0(Y,\oc_Y) \to H^0(Y,\nu_* \oc_{\tilde Y})\to H^0(Y,  \nu_* \oc_{\tilde Y}/\oc_Y )    \to H^1(Y,\oc_Y) \to H^1(Y,\nu_* \oc_{\tilde Y})   \to 0.
      \end{align}
      On en déduit qu'on a
         \begin{align}\label{eq:ineg-H1}
         \dim_K(H^1(\tilde X_a,\oc_{\tilde X_a}))=\dim_{K'}(H^1(Y, \oc_{Y}))   \geq \dim_{K'}(H^1(Y,\nu_* \oc_{\tilde Y}) )=\dim_{K'}(H^1(\tilde Y, \oc_{\tilde Y})).
         \end{align}
         De même, on voit qu'on a
         \begin{align*}
         \dim_K(H^1(X_a,\oc_{X_a}))\geq \dim_K(H^1(\tilde X_a,\oc_{\tilde X_a})).
         \end{align*}
         Si la différence \eqref{eq:diff} était non nulle, on aurait
         \begin{align*}
            q_a= \dim_K(H^1(X_a,\oc_{X_a}))\geq \dim_K(H^1(\tilde X_a,\oc_{\tilde X_a}))-\dim_{K'}(H^1(\tilde Y, \oc_{\tilde Y}))\geq (p-1)/2
         \end{align*}
ce qui contredirait l'inégalité \eqref{eq:hyp-p}. L'expression \eqref{eq:diff} est donc nulle et l'inégalité \eqref{eq:ineg-H1} est une égalité d'où l'on tire $\dim_{K'}(H^1(Y, \oc_{Y}))  =\dim_{K'}(H^1(Y,\nu_* \oc_{\tilde Y}) )$. Comme on a par ailleurs $\dim_{K'}(H^0(Y, \oc_{Y}))  =\dim_{K'}(H^0(Y,\nu_* \oc_{\tilde Y}) )=1$, la suite exacte \eqref{eq:suite-Hi} entraîne qu'on a
\begin{align*}
  H^0(Y,  \nu_* \oc_{\tilde Y}/\oc_Y )=(0)
\end{align*}
et donc $\oc_Y =\nu_* \oc_{\tilde Y}$. Ainsi $Y$ est une courbe normale ce qu'il fallait voir.

    \end{preuve}
  
\end{paragr}

\begin{paragr}[Inégalité de Severi.] --- \label{S:delta} Pour tout $a\in A^{\el}$, on définit l'invariant $\delta_a$ par
  \begin{align*}
    \delta_a=\dim_{K_a}(H^0(X_a, \nu_*\oc_{\tilde X_a}/\oc_{X_a}))
  \end{align*}
  où $K_a$ est le corps résiduel de $a$ et $\nu:\tilde X_a\to X_a$ est le morphisme de normalisation. On obtient ainsi une fonction $A^{\el}\to \NN$ qui est constructible et semi-continue supérieurement. Pour tout $\delta\in \NN$, on définit   $A^{\el, \geq \delta}\subset A^{\el}$ comme le sous-schéma fermé réduit défini par la condition $\delta_a\geq \delta$. On appelle \emph{inégalité de Severi} la famille d'inégalités suivantes, où $q_a$ est le genre arithmétique  commun à toutes les courbes spectrales $X_a$  et ne dépend donc pas de $a\in A$:
  \begin{align}\label{eq:Severi}
         \dim(A^{\el, \geq \delta})\leq \dim(A^{\el})-\delta   \text{ pour tout }   0\leq \delta \leq q_a.
    \end{align}

  \begin{proposition}\label{prop:Severi} Soit $p$ un premier qui vérifie l'inégalité \eqref{eq:hyp-p}.  Alors l'inégalité de Severi \eqref{eq:Severi} est vérifiée lorsque $k$ est la clôture algébrique d'un corps fini de caractéristique $p$.   
  \end{proposition}

  \begin{preuve}
    Comme  $A^{\el, \geq \delta}$ est un schéma réduit sur un corps parfait, il existe un ouvert dense $B\subset A^{\el, \geq \delta}$ qui est lisse. Quitte à restreindre $B$, on peut et on va également supposer que l'application $a\mapsto \delta_a$ est constante de valeur $\delta$ sur $B$. Soit $\Xgo=X_{|B}$ la restriction à $B$ de la courbe spectrale universelle. Soit $\nu:\tilde\Xgo\to \Xgo$ le morphisme de normalisation. D'après un résultat de Teissier, voir par exemple \cite[preuve de la proposition A.2.1]{Lau-jac}, le morphisme obtenu par composition $\tilde\Xgo\to B$ est plat et, pour tout $a\in B$, le morphisme $\nu_a:\tilde\Xgo_a\to \Xgo_a=X_a$ est la normalisation de $X_a$. Il résulte alors de la  proposition \ref{prop:lissite} que les fibres du morphisme composé
    \begin{align*}
      \xi:\tilde\Xgo\to B
    \end{align*}
    sont lisses et que ce morphisme plat est donc lisse.

    Pour tout  $a\in B$, soit $\varphi_a=\rho_a\circ \nu_a:\tilde \Xgo_a\to \Sigma$, cf. le diagramme \eqref{eq:diag}. Ces morphismes s'organisent en un morphisme de $B$-schémas
    \begin{align*}
      \varphi:\tilde \Xgo\to B\times \Sigma.
    \end{align*} 
    Soit
    \begin{align*}
    d\varphi : T\tilde\Xgo  \to \varphi^* T(B\times\Sigma)
  \end{align*}
  la différentielle de $\varphi$. Soit $\nc$ le  conoyau de  $d\varphi$.   Soit $p_1: B\times \Sigma\to B$ la première projection. On a $\xi=p_1\circ \varphi$. Soit $\tilde \kappa$ défini comme la composition
    \begin{align*}
      \tilde \kappa: \xi^* TB \hookrightarrow \varphi^* T(B\times \Sigma)\to \nc.
    \end{align*}
    On pose
    \begin{align*}
       \kappa:   TB \to  \xi_*\xi^* TB  \to_{\xi_* \tilde\kappa}   \xi_*\nc.
    \end{align*}
    Finalement, pour $a\in B$, soit  $\nc_a=\nc\otimes\oc_{\tilde \Xgo_a}$. On obtient, par composition de  $T_aB \to  (\xi_*\nc)_a$ avec l'application naturelle  $(\xi_*\nc)_a\to H^0(\tilde \Xgo_a,\nc_a)$ une application, dite  d'Horikawa et notée $\kappa_a$,
    \begin{align*}
      \kappa_a: T_aB \to H^0(\tilde \Xgo_a,\nc_a).
    \end{align*}
   Le faisceau $\nc_a$ contient un sous-faisceau de torsion $\nc_{a,tors}$ supporté sur le lieu où   $d\varphi_a$ n'est pas injective. On pose $\nc_a'=\nc_a/\nc_{a,tors}$.

   Le lemme suivant est une variation sur une observation d'Arbarello et Cornalba, cf.  \cite[p. 26]{Petri} ou \cite[lemme 2.3]{Capo}.
   
  \begin{lemme}\label{lem:injectivite}
    L'application induite
    \begin{align*}
         \kappa_a': T_aB \to H^0(\tilde \Xgo_a,\nc_a')
    \end{align*}
    est injective.
  \end{lemme}

  \begin{preuve}
    On a  un diagramme commutatif à lignes exactes:
    \begin{align}\label{eq:diag-2l}
      \xymatrix{   0   \ar[r] & T_a\tilde \Xgo  \ar[r]\ar[d]& \varphi_a^* T\Sigma \ar@{=}[d]\ar[r]& \nc_a \ar[d]\ar[r]& 0\\
            0   \ar[r]&   \hc om_{\oc_{X_a}}(  \Om^1_{X_a/k},   \oc_{\tilde \Xgo_a})  \ar[r]&   \varphi_a^* T\Sigma \ar[r]&    \hc om_{\oc_{X_a}}( \Ic_a/\Ic_a^2 , \oc_{\tilde \Xgo_a})     &}
    \end{align}
    la ligne inférieure étant obtenue par application du foncteur   $\hc om_{\oc_{X_a}}(  \cdot,   \oc_{\tilde \Xgo_a}) $ à la suite exacte courte

    \begin{align}\label{eq:courte}
     \xymatrix{ 0   \ar[r]& \Ic_a/\Ic_a^2     \ar[r]   &\rho_a^* \Om^1_{\Sigma/C}   \ar[r] &\Om^1_{X_a/C}   \ar[r] &0}.
    \end{align}
    Observons que le faisceau  $\hc om_{\oc_{X_a}}( \Ic_a/\Ic_a^2 , \oc_{\tilde \Xgo_a})$ est localement libre. La flèche verticale à droite se factorise par un morphisme injectif:
    \begin{align*}
      \nc_a'\hookrightarrow  \hc om_{\oc_{X_a}}( \Ic_a/\Ic_a^2 , \oc_{\tilde \Xgo_a}).
    \end{align*}
    On obtient ainsi une application injective
    \begin{align*}
      H^0(\tilde \Xgo_a,\nc_a') \hookrightarrow H^0(\tilde \Xgo_a,\hc om_{\oc_{X_a}}( \Ic_a/\Ic_a^2 , \oc_{\tilde \Xgo_a}))=H^0(C, \tilde\pi_{a,*}\oc_{\tilde \Xgo_a} \otimes_{\oc_C} \oc(nD)),
    \end{align*}
    où l'on pose $\tilde\pi_a=\pi_a\circ \nu_a$.

  Par composition avec $\kappa'$, on obtient une application qui se factorise ainsi:
  \begin{align*}
    T_aB \to H^0(C, \pi_{a,*}\oc_{ X_a} \otimes_{\oc_C} \oc_C(nD)) \to H^0(C, \tilde\pi_{a,*}\oc_{\tilde \Xgo_a} \otimes_{\oc_C} \oc_C(nD)).
  \end{align*}
  Cette application est injective puisque la première flèche s'identifie à $ T_aB \to T_a A$, cf. § \ref{S:deformation}, et que la seconde est évidemment injective. L'injectivité cherchée est alors claire.
\end{preuve}

Il résulte alors du lemme \ref{lem:injectivite} qu'on a
\begin{align*}
  \dim(B)=\dim(T_aB)\leq \dim_K(H^0(\tilde \Xgo_a,\nc_a'))
\end{align*}
où $K$ est le corps résiduel de $a$.
On conclut par le lemme suivant qui est analogue au \cite[corollaire 2.4]{Capo}.

\begin{lemme}
  On a
  \begin{align*}
     \dim_K(H^0(\tilde \Xgo_a,\nc_a'))\leq \dim(A^{\el})-\delta.
  \end{align*}
\end{lemme}

\begin{preuve}
  Soit $Z$ le diviseur de ramification de $\varphi_a$. Il résulte de la première ligne du diagramme \eqref{eq:diag-2l} qu'on a une suite exacte de faisceau localement libre sur $\tilde \Xgo_a$.
      \begin{align*}
      \xymatrix{   0   \ar[r] & (T_a\tilde\Xgo)(Z)  \ar[r]& \varphi_a^* T\Sigma \ar[r]& \nc_a' \ar[r]& 0}.
      \end{align*}
      On en déduit qu'on a un isomorphisme de fibrés en droites
      \begin{align*}
        \nc_a' \simeq \Om^1_{\tilde\Xgo_a/K}(-Z)\otimes_{\oc_{\tilde\Xgo_a}} \varphi_a^*(\om^1_{\Sigma/K})^{-1}
      \end{align*}
      où $\om^1_{\Sigma/K}$ est le carré extérieur du faisceau des différentielles $ \Om^1_{\Sigma/K}$.     Ainsi,   $\nc_a' $ s'identifie à un sous-faisceau du fibré en droites  $\Om^1_{\tilde\Xgo_a/K}\otimes_{\oc_{\tilde\Xgo_a}} \varphi_a^*(\om^1_{\Sigma/K})^{-1}$ ce qui entraîne qu'on a
      \begin{align*}
                \dim_K(H^0(\tilde \Xgo_a,\nc_a'))\leq   \dim_K(H^0(\tilde \Xgo_a, \Om^1_{\tilde\Xgo_a/K}\otimes_{\oc_{\tilde\Xgo_a}} \varphi_a^*(\om^1_{\Sigma/K})^{-1}))
      \end{align*}
      Pour conclure, on va calculer la dimension du membre de droite à l'aide de la formule de Riemann-Roch. On commence par obtenir  à l'aide du lemme \ref{lem:deg-Om} ci-dessous et de l'inégalité \eqref{eq:degD}

      \begin{align*}
        \deg(\Om^1_{\tilde\Xgo_a/K}\otimes_{\oc_{\tilde\Xgo_a}} (\om^1_{\Sigma/K})^{-1})&=  \deg(\Om^1_{\tilde\Xgo_a/K})-\deg( \varphi_a^*\Om^1_{\Sigma/K})\\
                                                                                        &=2g_{ \tilde\Xgo_a }-2+n(\deg(D)-2g_C +2)\\
        &> 2g_{ \tilde\Xgo_a }-2.
      \end{align*}
On a donc par la formule de Riemann-Roch  et  l'égalité \eqref{eq:dA-qa}
\begin{align*}
  \dim_K(H^0(\tilde \Xgo_a, \Om^1_{\tilde\Xgo_a/K}\otimes_{\oc_{\tilde\Xgo_a}} (\om^1_{\Sigma/K})^{-1}))&=\chi(\tilde \Xgo_a, \Om^1_{\tilde\Xgo_a/K}\otimes_{\oc_{\tilde\Xgo_a}} (\om^1_{\Sigma/K})^{-1}))\\
                                                                                                        &= \deg(\Om^1_{\tilde\Xgo_a/K}\otimes_{\oc_{\tilde\Xgo_a}} (\om^1_{\Sigma/K})^{-1})+1-g_{ \tilde\Xgo_a }\\
                                                                                                        &=g_{ \tilde\Xgo_a }-1+n(\deg(D)-2g_C +2)\\
                                                                                                        &= q_a-\delta_a -1+n(\deg(D)-2g_C +2)\\
  &=d_A-\delta_a.
\end{align*}
\end{preuve}

\begin{lemme}
        \label{lem:deg-Om}
        On a
        \begin{align*}
          \deg( \varphi_a^*\Om^1_{\Sigma/K})=n(2g_C - 2 - \deg(D)).
        \end{align*}
      \end{lemme}
      \begin{preuve}
       Le faisceau $\varphi_a^*\Om^1_{\Sigma/K}$ s'incrit dans une suite exacte
   \begin{align*}
     0\to \varphi_a^*p^*\Om^1_C \to \varphi_a^*\Om^1_{\Sigma/K} \to \varphi_a^*\Om^1_{\Sigma/C} \to 0,
   \end{align*}
   l'exactitude à gauche résulte de ce que $\varphi_a^*p^*\Om^1_C$ est localement libre. On a $\Om^1_{\Sigma/C} =p^{*}\oc_C(-D)$. Le résultat est alors une conséquence immédiate du fait que pour tout  fibré en droites $\lc$ sur $C$ on a 
    \begin{align*}
      \deg( \varphi_a^*p^*\lc)=n\deg(\lc).
    \end{align*}
      \end{preuve}

\end{preuve}
  
\end{paragr}

\section{Fibration de Hitchin}\label{sec:Hitchin}

\subsection{Cadre}\label{ssec:cadre-Hitchin}

\begin{paragr} On continue avec les notations de la section \ref{sec:cbe-spec}.
\end{paragr}

\begin{paragr} Par fibré de Hitchin sur $C$, on entend la donnée d'un couple $(\ec,\theta)$  formé d'un $\oc_C$-module $\ec$ localement libre de rang fini   et d'un endomorphisme tordu 
  \begin{align*}
    \theta :\mathcal{E}\rightarrow \mathcal{E}(D)=\mathcal{E}\otimes_{\mathcal{O}_{C}}\mathcal{O}_{C}(D).
  \end{align*}
  Le degré et le rang d'un fibré de Hitchin  $(\mathcal{E},\theta)$ sont le degré et le rang, respectivement notés $\deg(\ec)$ et  $\rang(\ec)$, du fibré vectoriel sous-jacent $\ec$.
\end{paragr}

  \begin{paragr}\label{S:stabilite} Un fibré de Hitchin $(\ec,\theta)$ est dit stable si pour tout sous-fibré $(0)\not=\mathcal{F}\subsetneq\mathcal{E}$ tel que $\theta(\mathcal{F})\subset\mathcal{F}(D)$, on a
$$
\mu(\mathcal{F})<\mu(\mathcal{E})
$$
où l'on introduit la pente
\begin{align*}
  \mu(\fc)=\frac{\deg(\fc)}{\rang(\fc)}.
\end{align*}
\end{paragr}

\begin{paragr}\label{S:Hitchin}  Rappelons que depuis la section  \ref{sec:cbe-spec}, on  a fixé un entier $n\geq 1$. Pour tout  $e\in \ZZ$, soit $M^e$ le schéma de modules grossier qui paramètre les classes d'isomorphisme de fibrés de Hitchin sur $C$ qui sont stables, de degré $e$ et rang $n$. Lorsque le contexte sera clair, on pourra omettre l'exposant $e$.  D'après \cite[proposition 7.4]{Nitsure}, c'est un schéma quasi-projectif et   lisse sur $k$, connexe, purement de dimension
  \begin{align*}
    d_{M}=n^{2}\deg(D)+1.
  \end{align*}
    La fibration de Hitchin est le morphisme de schémas
\begin{equation*}
f:M\rightarrow A
\end{equation*}
qui envoie $(\mathcal{E},\theta)$ sur le polynôme caractéristique de $\theta$
$$
a_\theta(X)=X^n-\mathrm{tr}(\theta)X+\ldots +(-1)^{n}\mathrm{det}(\theta).
$$
Lorsque $e$ est premier au rang $n$, ce qui sera désormais notre hypothèse, le  morphisme $f$  est projectif, plat,   à fibres connexes  purement de dimension
\begin{align}\label{eq:df}
d_{f}=  n(g_C-1)+\frac{n(n-1)}{2}\deg(D)+1,
\end{align}
voir \cite[théorème 6.1]{Nitsure} et \cite[corollaire 8.2]{Support}. Pour tout $a\in A$, soit $M_a$ la fibre de $f$ au-dessus de $a$. 
\end{paragr}

\subsection{Action du schéma de Picard}\label{ssec:actionJa}

\begin{paragr}\label{S:Ja} Soit $J_{X/  A}$ l'espace grossier de la composante neutre du champ de Picard relatif de $X/  A$.  C'est un schéma  en groupes, lisse et de type fini sur $ A$, cf. \cite[section 5]{Support} pour quelques rappels. Soit $a\in A$. La fibre $J_{X_a}$, noté simplement $J_a$, est le schéma en groupes des classes d'isomorphisme de fibrés inversibles sur $X_{a}$ dont la restriction à chaque composante irréducible de $X_{a}$ est de degré $0$.  

  Par une correspondance due à Hitchin, Beauville-Narasimhan-Ramanan, cf. \cite{Hitchin,BLR,Schaub} par ordre croissant de généralité, tout fibré de Hitchin $(\ec,\theta)$ de rang $n$ et de polynôme caractéristique $a$ correspond de manière biunivoque à un module $\fc$ sans torsion et  de rang $1$ sur la courbe spectrale $X_{a}$. Rappelons simplement que le fibré de Hitchin associé à un tel $\fc$ est le fibré vectoriel $\ec=\pi_{a,*}\fc$ muni de l'endomorphisme induit par la section universelle $u$. Dans cette correspondance, on a
   \begin{align*}
      \deg(\fc)&:=\chi(X_{a},\fc)-\chi(X_a,\oc_{X_a})\\
      &=\deg(\ec)+\frac{n(n-1)}{2}\deg(D).
   \end{align*}
   On dit que $\fc$ est stable si le fibré de Hitchin correspondant est stable.  Pour une traduction de la condition de stabilité purement  en terme du module $\fc$, on renvoie le lecteur à \cite[section 4]{Support}.   Pour tout module inversible $\lc$ sur $X_a$ et tout module $\fc$ sans torsion et  de rang $1$ sur $X_a$, le produit tensoriel $\lc\otimes_{\oc_{X_a}}\fc$ est encore sans torsion et  de rang $1$. De plus, si le degré de $\lc$ est  $0$ sur chaque composante irréducible de $X_a$, alors $\lc\otimes_{\oc_{X_a}}\fc$ est, d'une part, de même degré que $\fc$ et, d'autre part,  stable si et seulement si $\fc$ l'est.  Via la correspondance évoquée plus haut, on obtient une action de  $J_{X/  A}$ sur $M$ et cette action se fait fibre à fibre, le schéma en groupes  $J_{a}$ agissant sur $M_a$ pour tout $a\in A$.

   Pour $a\in A^{\el}$, la courbe spectrale est intègre : le schéma en groupes $J_a$ agit simplement transitivement sur l'ouvert dense $M_A^{\reg}$ de $M_a$ formé de fibrés de Hitchin correspondant à des fibrés en droites sur $X_a$. La dimension du complémentaire $M_a\setminus M_a^{\reg}$ est de dimension $\leq d_f-1$, cf. \cite[Proposition 4.16.1]{Ngo-LF} et le résultat  \cite{AIK} de Altman, Iarrobino et Kleiman sur la densité de la jacobienne dans la jacobienne compactifiée.
 \end{paragr}

 \begin{paragr}[La norme.] ---   \label{S:Norme}   Soit $P_C$ (ou simplement $P$) le schéma de Picard de $C$ qui paramètre les classes d'isomorphisme des $\oc_C$-modules inversibles. C'est un  $k$-schéma en groupes abéliens pour le produit tensoriel des $\oc_C$-modules inversibles. L'application degré identifie son groupe des composantes connexes à $\ZZ$. Pour tout $i\in \ZZ$, soit $P^i_C$ la composante connexe formée des $\oc_C$-modules inversibles de degré $i$. Soit  $J_C$ (ou simplement $J$) la jacobienne de $C$ c'est-à-dire la composante neutre $P_C^0$.

   Pour tout $\oc_C$-module $\ec$ localement libre de rang $r$, on note $\det(\ec)=\La^r\ec$ le fibré en droites donné la puissance extérieure maximale de $\ec$.
    Soit $a\in  A$. On a  un  morphisme de groupes (voir \cite[section 6.5]{EGA2}) induit par la norme
  \begin{align*}
    N_{X_a/C}: J_{a}\to J.
  \end{align*}
  Pour tout module inversible $\lc$ sur $X_a$ et tout module $\fc$ sans torsion et  de rang $1$ sur $X_a$, on a (voir \cite[proposition 3.10]{Hausel-Pauly})
    \begin{align*}
      \det(\pi_{a,*}(\lc\otimes_{\oc_{X_a}}\fc))= \det(\pi_{a,*}\fc)\otimes_\oc N_{X_a/C}(\lc).
    \end{align*}

    Vue comme diviseur effectif sur $\Sigma$, la courbe spectrale $X_a$ s'écrit
    \begin{align*}
      \sum_{i=1}^r m_i X_{a_i}
    \end{align*}
où $m_i\geq 1$ sont les multiplicités et $X_{a_i}$ est la courbe spectrale intègre associée au facteur irréductible $a_i\in A_{n_i}$ avec $n_i\geq 1$ et $\sum_{i=1}^r m_in_i=n$.
Soit $\nu_i:\tilde X_{a_i}\to X_{a_i}$ la normalisation  de $X_{a_i}$. Soit $\tilde J_{a_i}$ la jacobienne de $\tilde X_{a_i}$. On a un homomorphisme $\nu_i^*:J_{a_i}\to \tilde J_{a_i}$.   On a également un homomorphisme norme
          \begin{align*}
            N_{\tilde X_{a_i}/C}: \tilde J_{a_i}\to J.
          \end{align*}

          \begin{proposition} \label{prop:devissage}L'application $\lc\mapsto(\lc_i)_{1\leq i\leq r}$ donnée par 
\begin{align}\label{eq:lci}
 \lc_i=\nu_i^*(\lc\otimes_{\oc_{X_a}} \oc_{X_{a_i}} ).
\end{align}
induit un homomorphisme
            \begin{align}\label{eq:Ja-Jai}
                J_a\to \prod_{i=1}^{r}\tilde J_{a_i}
            \end{align}
 surjectif, à noyau un groupe algébrique affine connexe noté $J_a^{\aff}$.
          \end{proposition}
          
  \begin{proposition}\label{prop:norme}
    Pour tout fibré inversible $\lc$ sur $X_a$, on  a
    \begin{align*}
      N_{X_a/C}(\lc)=\otimes_{i=1}^{r}  N_{\tilde X_{a_i}/C}(\lc_i)^{m_i}
    \end{align*}
    où $\lc_i$ est défini par \eqref{eq:lci}.    
\end{proposition}

\begin{preuve}
 On renvoie le lecteur aux preuves des \cite[lemmes 3.4 à 3.6]{Hausel-Pauly}.
\end{preuve}
\end{paragr}

\section{Théorèmes de support}\label{sec:coh}

\subsection{Conventions} \label{ssec:convention}

\begin{paragr}
  On suppose désormais que $k$ est la clôture algébrique d'un corps fini  $\Fq$ de cardinal $q$. Soit $\Fr$ l'automorphisme de Frobenius géométrique, c'est-à-dire le générateur topologique du groupe de Galois $\Gal(k/\Fq)$ donné par l'inverse de $x\mapsto x^q$. 
\end{paragr}

\begin{paragr} Soit $X_0$ un  $\Fq$-schéma de type fini. On ôte l'indice $0$ pour désigner le changement de base à $k$ c'est-à-dire on note $X=X_0\times_{\Fq} k$.  Cette convention vaut pour les $\Fq$-schémas et les $\Fq$-morphismes qui seront affublés d'un indice $0$. On munit $X_0$  de l'endomorphisme de Frobenius $F_{X_0}:X_0\to X_0$, qui sera souvent noté  simplement $F$: c'est le  morphisme fini qui est l'identité sur l'espace topologique sous-jacent et qui est donné par l'élévation à la puissance $q$ sur le faisceau structural. Par changement de base à $k$, on obtient l'endomorphisme de Frobenius (\og relatif\fg{}) $F_X=F_{X_0}\times \Id_{k}$ qu'on note parfois simplement $F$.  Pour tout entier $n\geq 1$, l'ensemble les points fixes de $F_X^n=F_X\circ\cdots\circ F_X$ ($n$ fois) sur $X(k)$ est $X_0(\Fqm)$. On dispose également d'une action du Frobenius géométrique sur $X$ par $\Id_{X_0}\times \Fr$. L'action de $\Fr$ qui s'en déduit sur $X(k)$ admet  $X_0(\Fq)$ pour ensemble de points fixes.
\end{paragr}

\begin{paragr} Soit $\ell$ un nombre premier distinct de la caractéristique de $k$ et $\Qlb$ une clôture algébrique du corps $\Ql$ des nombres $\ell$-adiques. Soit $\Dbc(X,\Qlb)$ la catégorie dérivée des faisceaux $\ell$-adiques sur $X$ et $\Perv(X,\Qlb)$  la sous-catégorie pleine des faisceaux pervers définie par Beilinson, Bernstein, Deligne et Gabber \cite{Ast100}. Cette sous-catégorie est abélienne et préservée par la dualité de Verdier notée simplement $\D$. Pour tout $i\in \ZZ$, on dispose d'un foncteur de cohomologie ${}^{\mathrm{p}}\mathcal{H}^{i}:\Dbc(X,\Qlb)\to \Perv(X,\Qlb)$. Les faisceaux pervers simples sont de la forme
        \begin{align}\label{eq:pervers-simple}
          \IC_{Z}(\fc)[d_Z]
        \end{align}
        où $Z\subset X$ est un sous-schéma intègre de dimension $d_Z$ et  $\fc$ un système local $\ell$-adique irréductible sur un ouvert dense de $Z$. On a noté $ \IC_{Z}(\fc)$ l'image directe dans $X$ du complexe d'intersection de $Z$ à valeurs dans $\fc$. Les objets de $\Perv(X,\Qlb)$ sont de longueur finie et les faisceaux pervers semi-simples sont des sommes finies de faisceaux pervers simples.
      \end{paragr}
      
\begin{paragr}
  On définit comme ci-dessus $\Dbc(X_0,\Qlb)$ et la sous-catégorie  $\Perv(X_0,\Qlb)$. Soit $\gc_0\in  \Dbc(X_0,\Qlb)$. Soit $\gc\in  \Dbc(X,\Qlb)$ déduit de $\gc_0$ par changement de base. Alors $\gc$ est muni d'un isomorphisme $\Fr^*:F_X^*\gc \simeq \gc$. En particulier, pour tout point $x\in X_0(\Fq)$, on dispose d'un automorphisme $\Fr^*$ de $\gc_x$ et on note:
  \begin{align*}
    \trace(\Fr,\gc_x)=\sum_{i\in \ZZ} (-1)^i \trace(\Fr^*, \hc^i(\gc)_x).
  \end{align*}
\end{paragr}

\subsection{Revêtements cycliques}\label{ssec:revet}

\begin{paragr}
  Soit $J_0$ la jacobienne de $C_0$ et $J=J_0\times_{\Fq} k$. 
   
  L'isogénie de Lang 
\begin{align}\label{eq:Lang}
   \Lgo: J_0\to J_0
\end{align}
donnée par $\lc\mapsto F^*\lc\otimes\lc^{-1}$ fait de   $J_0$ un   $J_0(\Fq)$-torseur au-dessus de  $J_0$.
\end{paragr}

\begin{paragr}[Application d'Abel-Jacobi. ]\label{S:torseur} ---   Soit $P_0$ le schéma de Picard de $C_0$. Pour tout $i\in \ZZ$, soit $P_0^i$ la composante  des fibrés de degré $i$. On dispose de l'application d'Abel-Jacobi  $AJ:C_0\to P^1_0$. Soit $\al\in P_0^1(\Fq)$ et  $AJ_\al:C_0\to J_0$ la composition de $AJ$ avec le morphisme  $P_0^1\to J_0$ donné par $\lc\mapsto\lc\otimes \al^{-1}$. 
  En tirant en arrière le torseur  de Lang \eqref{eq:Lang}  par l'application $AJ_\al$,  on obtient un revêtement $\tilde\rho_\al:\tilde C_\al\to C_0$ fini, étale, galoisien  de groupe de Galois $J_0(\Fq)$. Soit $\be\in P_0^1(\Fq)$ et $\lc_0\in J_0(k)$ tel que $\Lgo(\lc_0)=\al^{-1}\be$.  On a alors un diagramme commutatif

\begin{align*}
    \xymatrix{ \tilde C_\al \ar[r]^{}  \ar[d]^{\tilde\rho_\al} & J_0 \ar[r]^{\otimes \lc_0}  \ar[d]^{\Lgo} & J_0 \ar[d]^{\Lgo}  \\
 C_0 \ar[r]^{AJ_\al} & J_0\ar[r]^{\otimes \al^{-1}\be} & J_0}
\end{align*}
  qui induit un isomorphisme $\phi:\tilde C_{\al}\times_{\Fq} k\to \tilde C_{\be}\times_{\Fq} k$. On prendra garde que ce morphisme ne se descend en général pas à $\Fq$: en fait,  on a $ \phi\circ \Fr= \al\be^{-1}\Fr\circ \phi$.
\end{paragr}

\begin{paragr}[Caractère $\chi$.] --- \label{S:chi}   Soit $\ell$ un nombre premier inversible dans $k$. Soit 
  \begin{align}\label{eq:chi-J0}
    \chi: J_0(\Fq)\to \Qlb^\times
  \end{align}
  un caractère et $d$ son ordre.   Soit $\gc_\chi$ le $\Qlb$-faisceau lisse de rang $1$ sur  $J_0$  obtenu lorsqu'on pousse $\chi^{-1}$ par le torseur défini en § \ref{S:torseur}.

Soit
  \begin{align*}
    \Ga_\chi=J_0(\Fq)/ \Ker(\chi).
  \end{align*}
  Le caractère $\chi$ identifie $\Gamma_\chi$ au  sous-groupe des racines $d$-ièmes de l'unité de $\Qlb^\times$.  Pour tout $m\geq 1$ et tout $\beta\in J_0(\Fqm)$ on a $\trace(\Fr^m,\gc_{\chi,\be})=\chi(N_m(\be))$ où $N_m(\be)=\be F(\be)\cdots F^{m-1}(\be)\in P_0(\Fq)$.
\end{paragr}

\begin{paragr} \label{S:C chi al}En prenant le quotient par l'action du sous-groupe $\Ker(\chi)$, on déduit du revêtement  $\tilde\rho_\al:\tilde C_\al\to C_0$  un revêtement $\rho_{\chi,\al}: C_{\chi,\al} \to C_0$ fini, étale, galoisien  de  groupe de Galois $\Ga_\chi$.  Notons que  $C_{\chi,\al}$ ne dépend en fait que du noyau de $\chi$ et de l'image de $\al$ dans le quotient $P_0^1(\Fq)/\Ker(\chi)$.
\end{paragr}

\begin{paragr}[Variante sur $\Fqm$.] --- \label{S:variante}Soit $m\geq 1$ et $\al\in P_0^1(\Fqm)$. On utilisera la variante suivante: soit
  \begin{align*}
    \tilde\rho_{\al,m}:\tilde C_{\al,m}\to C_0\times_{\Fq}\Fqm
  \end{align*}
 le revêtement galoisien de groupe $J_0(\Fqm)$ qu'on obtient en tirant en arrière le torseur de Lang $\lc\mapsto F^{m,*}\lc\otimes\lc^{-1}$ par l'application
  \begin{align*}
    AJ_\al:  C_0\times_{\Fq}\Fqm\to J_0\times_{\Fq}\Fqm.
  \end{align*}
 En prenant le quotient par le sous-groupe $\Ker(\chi\circ N_m)$, on obtient    un revêtement $\rho_{\chi,\al,m}: C_{\chi,\al,m} \to C_0$ fini, étale, galoisien  de  groupe de Galois $ J_0(\Fqm)/\Ker(\chi\circ N_m)$.  Ce groupe s'identifie à $\Ga_\chi$ via la norme $N_m$. 
 Pour tout caractère $\mu:J_0(\Fqm)\to \Qlb^\times$, soit $\fc_{\mu,\al}$ le $\Qlb$-faisceau lisse de rang $1$ sur $C_0\times_{\Fq}\Fqm$ qu'on obtient en poussant le revêtement $\tilde\rho_{\al,m}$  par $\mu^{-1}$.  Par abus, on note de la même façon son tiré en arrière sur $C$. On a 
 \begin{align*}
    H^1( C_{\chi,\al,m}\times_{\Fqm}k, \Qlb)=\bigoplus_{i=0}^{d-1}   H^1( C, \fc_{\chi^i,\al}).
 \end{align*}
 Pour $i=0$, on a simplement $\fc_{\chi^0,\al}=\Qlb$.
 Pour tous $\al,\be\in P_0^1(\Fqm)$, les revêtements  $\rho_{\chi,\al,m}$ et  $\rho_{\chi,\be,m}$ sont isomorphes si et seulement si $\al\be^{-1}\in \Ker(\chi\circ N_m)$. En tout cas,  ils le sont  après changement de base à $\mathbb{F}_{q^{md}}$. En particulier,  on peut identifier  les groupes de cohomologie mais l'action du Frobenius $\Fr^m$ dont ils sont équipés  dépend du choix de $\al$. Plus exactement, l'action de $\Fr^m$ sur $ H^1( C, \fc_{\chi^i,\be}) $ est celle sur $ H^1( C, \fc_{\chi^i,\al}) $  multipliée par $(\chi^i\circ N_m)(\al^{-1}\be)$.

 Supposons de plus $\al\in P_0^1(\Fq)$. On peut faire le lien avec les constructions du § \ref{S:torseur}. On a un morphisme naturel $\tilde C_{\al,m}\to \tilde C_\al \times_{\Fq} \Fqm$ au-dessus de  $C_0\times_{\Fq} \Fqm$. Il est équivariant pour l'action de $J_0(\Fq)$ qui agit sur  $\tilde C_\al \times_{\Fq} \Fqm$ via la norme $N_m$. Cela donne une identification des quotients $C_{\chi,\al,m}\simeq C_{\chi,\al} \times_{\Fq} \Fqm$.
\end{paragr}

\subsection{Morphisme entre les bases}\label{ssec:morphisme}

\begin{paragr}\label{S:C'} 
Dans cette sous-section, on fixe des entiers naturels $n,n',d$ de sorte que $n=dn'$. On fixe aussi  un élément  $\al\in P^1(\Fq)$ ainsi qu'un caractère $\chi: J_0(\Fq)\to \Qlb^\times$ d'ordre $d$. On note simplement $\Ga=\Ga_\chi$, cf. § \ref{S:chi}, et $\rho_0:C_0'\to C_0$ le revêtement  $\rho_{\chi,\al}: C_{\chi,\al} \to C_0$ de groupe de Galois $\Ga$ construit au § \ref{S:C chi al}. Par changement de base à $k$, on obtient un revêtement galoisien $\rho:C'\to C$ de groupe $\Gamma$ où  $C'=C_0'\times_{\Fq}k$ est une courbe  projective, lisse et connexe.
\end{paragr}

\begin{paragr}\label{S:morph}   Soit $A_0$ l'espace affine \eqref{eq:Ac} relatif à la courbe $C_0$,  à un diviseur $D_0$ fixé sur $C_0$ et à l'entier $n$.  Soit $A'_0$ l'espace affine \eqref{eq:Ac} relatif à la courbe $C'_0$, au diviseur $D'_0=\rho^*D_0$ et à l'entier $n'$. Le groupe $\Ga$ agit naturellement  sur $A_0'$. Soit
  \begin{align*}
    p'_0:\Sigma'_0=\mathbb{V}(\mathcal{O}_{C'_0}(-D'_0))\rightarrow C'_0
        \end{align*}
        l'espace total de $\mathcal{O}_{C'_0}(D'_0)$ qu'on obtient en tirant en arrière \eqref{eq:esp-total} par $\rho_0$.   La section universelle $u\in H^0(\Sigma_0, {p}_0^{\ast}\oc_{C_0}(D_0))$ induit la section universelle $u'\in H^0(\Sigma'_0, {p'}_0^{\ast}\oc_{C'_0}(D'_0))$. Pour tout $a'\in A'_0$,   on dispose de la section  ${a'}(u')\in H^0(\Sigma'_0, {p'}_0^{\ast}\oc_{C'_0}(n'D'_0))$. La section
          \begin{align*}
          \prod_{\sigma\in \Gamma} \sigma({a'}(u'))\in H^0(\Sigma'_0, {p'}_0^{\ast}\oc_{C'_0}(n D'_0))
        \end{align*}
        est fixe sous l'action naturelle de $\Gamma$ et s'identifie naturellement  à $a(u)$ pour un unique élément $a$ de $A_0$. On obtient ainsi un morphisme 
    \begin{align}\label{eq:iota}
            \iota_0:A'_0\to A_0
    \end{align}
    qui est $\Ga$-invariant  et fini sur son image, notée  $A_{\Ga,0}$, qui est un sous-schéma fermé de $A_0$.

    \begin{remarque}
      La notation $A_{\Ga,0}$ est un peu abusive en ce sens que $A_{\Ga,0}$ dépend \emph{a priori} aussi du choix de $\al\in P_0^1(\Fq)$. En revanche,  $A_{\Ga}=A_{\Ga,0}\times_{\Fq} k$ ne dépend que de $\Ga$.
    \end{remarque}

        \end{paragr}

 \begin{paragr}\label{S:img-iota} On définit comme dans la section \ref{sec:cbe-spec} et  la sous-section \ref{ssec:actionJa} la courbe spectrale $\pi_a:X_{a} \subset C_0$ et le schéma de Picard $J_a$. La proposition suivante donne une caractérisation des points de $A_{\Ga,0}$.

  \begin{proposition}\label{prop:img-iota}
Soit $a\in A_0(\Fq)$.  Les deux assertions sont équivalentes:
          \begin{enumerate}
          \item Il existe $a' \in A'_0(k)$ tel que $\iota_0(a')=a$.
          \item Le caractère $\chi \circ   N_{X_{a}/C_0} $ est trivial sur $J_{a}(\Fq)$ où le morphisme norme  $N_{X_{a}/C_0}$ est introduit au § \ref{S:Norme}.
          \end{enumerate}
        \end{proposition}

        \begin{preuve}  En utilisant le fait que la norme  $N_m:J_{a}(\Fqm) \to J_{a}(\Fq)$ est surjective (lemme de Lang), on voit que, quitte à remplacer $q$ par $q^m$, changer de base à $\Fqm$ tous les objets et remplacer $\chi$ par $\chi\circ N_m$, on peut (et on va) faire l'hypothèse suivante:
          \begin{itemize}
          \item   On a une décomposition en composantes irréductibles  $X_{a}=\sum_{i=1}^r m_i X_{a_{i}}$, cf § \ref{S:Norme}, avec $a_{i}\in A_{n_i,0}^{\el}(\Fq)$ et $\sum_{i=1}^r m_in_i=n$. En particulier, les courbes   spectrales $X_{a_{i}}$  sont géométriquement irréductibles.
          \end{itemize}

          Soit  $\tilde J_{a_i}$ la jacobienne de $\tilde X_{a_i}$. En utilisant le lemme de Lang et la proposition \ref{prop:devissage}, on a un morphisme surjectif $J_a(\Fq)\to \prod_{i=1}^r\tilde J_{a_i}(\Fq)$  induit par un morphisme défini comme en \eqref{eq:Ja-Jai}. Par la proposition \ref{prop:norme}, on voit que la condition 2 ci-dessus est équivalente à ce que la condition suivnate soit vérifiée pour tout $1\leq i\leq r$:
          
           \begin{itemize}
          \item[3.]           le caractère $\chi^{m_i} \circ   N_{\tilde X_{a_i}/C_0} $ est trivial sur $\tilde P_{a_i}(\Fq)$.
          \end{itemize}
\medskip

On a obtenu la courbe $C_0'$ comme quotient de $\tilde C_\al$ par $\Ker(\chi)$, cf. §§ \ref{S:C chi al} et  \ref{S:C'}. En quotientant $\tilde C_\al$ par le sous-groupe $\Ker(\chi^{m_i})$, on obtient  un revêtement fini étale cyclique $\rho_{i,0}:C_{i,0}'\to C_0$ de degré   $d_i=d/(m_i,d)$ de $F$ où $(m_i,d)$ est le pgcd de $m_i$ et $d$.  Le groupe de Galois de $C_{i,0}'/C_0$ est un quotient noté $\Ga_i$ de $\Ga$. À l'aide de la théorie du corps de classes, on  montre que la condition 3 équivaut à l'existence de $b_i$ dans l'espace $A'_{C_{i,0}',n_i/d_i}$ des polynômes caractéristiques sur  $C_{i,0}'$ de rang $n_i/d_i$ relatifs au diviseur $\rho_{i,0}^*D_0$ tel que la propriété suivante soit vérifiée:\\

     \begin{itemize}
     \item[4.]  Il existe une décomposition
       \begin{align*}
   a_i=\prod_{\sigma\in \Ga_i} \sigma(b_i).
 \end{align*}
\end{itemize}

\medskip

En posant $m_i'=m_i/(d,m_i)$ on a 
 \begin{align*}
   a=\prod_{i=1}^r a_i^{m_i}&= \prod_{i=1}^r \prod_{\sigma\in \Ga}  \sigma(b_i)^{m_i'}\\
   &= \prod_{\sigma\in \Ga}   \sigma( \prod_{i=1}^r b_i^{m_i'}).
    \end{align*}
    Le polynôme $b=\prod_{i=1}^r b_i^{m_i'}$ est alors de degré $\sum_{i=1}^r n_im_i'/d_i= \sum_{i=1}^r n_im_i/d=n'$. C'est un élément de $A'_0(\Fq)$ d'image $a$ par $\iota$. On a donc prouvé que 2 implique 1.

    Traitons la réciproque. Pour tout entier $l$, soit  $A'_{l,0}$ l'espace affine \eqref{eq:Ac}  relatif à $C_0'$, $D_0'$ et l'entier $l$. Quitte à prendre une extension plus grande du corps fini, on peut et on va supposer qu'il existe $a'\in A_0'(\Fq)$ d'image $a$ avec une décomposition $a'=\prod_{j=1}^{r'} c_j^{k_j}$ et les $c_j\in (A'_{l_j,0})^{\el}(\Fq)$ sont des polynômes caractéristiques deux à deux distincts pour $C_0'$ de degré respectif $l_j$ avec $\sum_{j=1}^{r'}k_jl_j=n'$. Soit $\Ga'_j$ le quotient de $\Ga$ par le stabilisateur de $c_j$. Soit $\delta_j$ l'ordre de $\Ga_j'$. On peut donc écrire:
    \begin{align*}
      a=\prod_{\sigma\in \Ga }   \sigma(a')= \prod_{j=1}^{r'}    ( \prod_{ \sigma\in \Ga'_j } \sigma(c_j))^{k_j d \delta_j^{-1}}.
    \end{align*}
    Les polynômes $\prod_{ \sigma\in \Ga'_j } \sigma(c_j)$ pour $1\leq j \leq r'$ sont deux à deux distincts  et appartiennent respectivement à $A_{\delta_j l_j,0}(\Fq)$. On a donc $r'=r$ et, quitte à réindexer, pour $1\leq i\leq r$ on a   $a_i=\prod_{ \sigma\in \Ga'_i  } \sigma(c_i)$, et  $m_i=k_id \delta_i^{-1} $ et $n_i=l_i\delta_i$. Il s'ensuit que $d|m_i\delta_i$ et donc l'ordre de $\Ga_i$ qui est $d/(m_i,d)$ divise $\delta_i$. Ainsi le groupe $\Ga_i$ est un quotient de $\Ga_i'$. On obtient alors la décomposition 4 ci-dessus pour $b_i=\prod_{\sigma\in \Ker(\Ga_i'\to \Ga_i)}\sigma(c_i)$. L'assertion 2 résulte alors de l'équivalence de 3 et 4.
  \end{preuve}
\end{paragr}

\begin{paragr} \label{S:A chi ell} Soit $A_{\Ga,0}^{\el}=A_{\Ga,0}\cap A_0^{\el}$. C'est un ouvert de $ A_{\Ga,0}$.

    \begin{proposition}\label{lem:gamma_a}
    \begin{enumerate}
    \item La partie localement fermée $A_{\Ga,0}^{\el}$ de $A_0$ est non vide.
    \item L'action de $\Gamma$ sur $\iota_0^{-1}(A_{\Ga,0}^{\el})$ est sans point fixe.
    \item Pour $a\in A_{\Ga,0}^{\el}(\Fq)$, l'ensemble des $a'\in A'(k)$ tel que $\iota_0(a')=a$ est un espace principal homogène sous l'action de $\Ga$. En particulier, il existe un unique $\gamma_{a,\Fr} \in \Gamma$ tel que pour tout $a'\in A'(k)$ tel que $\iota_0(a')=a$ on a $\gamma_{a,\Fr} \Fr(a')=a'$. 
    \end{enumerate}
  \end{proposition}

  \begin{remarque}
    Si  $a\in A_{\Ga,0}^{\el}(\Fqm)$ pour $m\geq 1$, la proposition  donne  un élément $\gamma_{a,\Fr^m}$.  Si, de plus, $a\in A_{\Ga,0}^{\el}(\Fq)$, on a la relation  $\gamma_{a,\Fr^m}=\gamma_{a,\Fr}^m$.
  \end{remarque}
  
  \begin{preuve} Soit $\Om\subset A_0'$  le plus grand ouvert sur lequel $\Ga$ agit sans point fixe. Pour des raisons de dimension, cet ouvert est non vide.   De même, l'ouvert elliptique $(A_{0}')^{\el}$ est non vide. Il est clair qu'on a $\Om\cap (A_{0}')^{\el}=\iota_0^{-1}(A_{\Ga,0}^{\el})$ ce qui vérifie 1 et 2.
Supposons que $a=\iota_0(a')$ appartienne à $A_{\Ga,0}^{\el}(\Fq)$.  Dans la décomposition de $\rho_0^*a$,  les facteurs irréductibles sont les $\sigma(a')$ pour $\sigma\in \Gamma$. La transitivité de l'action de $\Gamma$ sur la fibre $\iota_0^{-1}(a)$ s'en déduit. Ces facteurs  sont aussi permutés par l'action du Frobenius $\Fr$. Le 3 s'ensuit puisque les actions de $\Gamma$ et $\Fr$ commutent.
      \end{preuve}
  
\end{paragr}

\subsection{Un système local}\label{ssec:syst-loc}

\begin{paragr}   \label{S:al} On continue avec les notations de la sous-section précédente. On fixe un entier $e\in \ZZ$ premier à l'entier $n$. Soit $M_0$ l'espace de modules des fibrés de Hitchin de rang $n $ et de degré $e$ définis relativement au diviseur $D_0$.  Soit
  \begin{align*}
    \det: M_0 \to P_0^e
  \end{align*}
  le morphisme lisse qui, à $(\ec,\theta)$, associe la classe du fibré en droites $\La^n\ec$ donné par la puissance extérieure maximale de $\ec$.  Soit $\gc_{\chi,\al}$  le $\Qlb$-faisceau lisse de rang $1$ sur  $P_0^e$ qu'on obtient en tirant en arrière le faisceau $\gc_\chi$ du § \ref{S:chi} par l'application $P_0^e\to J_0$ donnée par la multiplication par $\al^{-e}$.

  Soit
    \begin{align}\label{eq:L chi al}
      \lc_{\chi,\al}={\det}^{*}\gc_{\chi,\al}.
    \end{align}
C'est un $\Qlb$-faisceau lisse de rang $1$ sur  $M_0$. Lorsque le contexte est clair, on omet l'indice $\al$. Par abus, on note  encore $\lc_\chi$  le $\Qlb$-faisceau lisse de rang $1$ sur  $M$ obtenu par extension  des scalaires.

        \begin{proposition}\label{prop:borne-supp}
          Pour tout $a\in A(k)$, on  a
          \begin{align*}
           (R f_{*} \lc_{\chi})_a=0
          \end{align*}
          sauf si $a$ appartient à $A_\Ga$.
        \end{proposition}

        \begin{preuve}
          Sans perte de généralité, on suppose que $a$ s'identifie à un élément de $a\in  A_0(\Fqm)$ pour $m\geq 1$ assez divisible. Un élément $x\in J_a(\Fqm)$ agit sur $(M_a,\lc_\chi)$ et l'action qui en résulte sur la cohomologie $H^\bullet(M_a,\lc_\chi)$ est donnée par la multiplication  par $\chi^{-1}(N_{X_a/C_0}(x))$. Par ailleurs, comme $J_a$ est connexe, un argument d'homotopie, cf.  \cite[exposé 6, lemme 2.8]{SGA4.5}, montre que cette action ne dépend pas de $x$. Autrement dit la multiplication par $\chi^{-1}(N_{X_a/C_0}(x))-1$ donne $0$ sur $H^\bullet(M_a,\lc_\chi)$. Donc, si on a $H^\bullet(M_a,\lc_\chi)\not=0$, le caractère  $\chi\circ N_{X_a/C_0}$ est trivial sur $J_a(\Fqm)$ ce qui implique que $a$ appartient à l'image du morphisme $\iota$ d'après la proposition \ref{prop:img-iota}.
        \end{preuve}
      \end{paragr}

      \subsection{Un théorème de support}\label{ssec:thm support}

      \begin{paragr}\label{S:def support}
        Le dual de Verdier du complexe de faisceaux  $\ell$-adiques $\lc_{\chi}[d_M]$ sur $M_0$ est  $\lc_{\chi^{-1}}[d_M](d_M)$. Il s'ensuit que $\lc_{\chi}[d_M]$ est un complexe pur de poids $d_M$. Le morphisme $f_0$ étant propre, par \cite[proposition 5.1.14]{Ast100} conséquence de \cite{Weil2}, le complexe de faisceaux $Rf_{0,*} \lc_\chi [d_M]$ sur $A_0$ est aussi pur de poids $d_M$. D'après \cite[théorème 5.4.5]{Ast100}, on a
        \begin{align}\label{eq:sum-fxpervers}
       Rf_{*} \lc_\chi [d_M]   \cong\bigoplus_{i\in \ZZ}{}^{\mathrm{p}}\mathcal{H}^{i}(Rf_{*} \lc_\chi [d_M])[-i].
        \end{align}
        En fait, d'après \cite[note 94 p. 177]{Ast100}, \cite[proposition 2.1]{Cataldo-IMRN} ou \cite[corollaire 3.2.5]{Sun-Zheng}, on a aussi une décomposition
        
 \begin{align}\label{eq:sum-fxpervers0}
       Rf_{0,*} \lc_\chi [d_M]   \cong\bigoplus_{i\in \ZZ}{}^{\mathrm{p}}\mathcal{H}^{i}(Rf_{0*} \lc_\chi [d_M])[-i].
        \end{align}
   Les  faisceaux pervers  ${}^{\mathrm{p}}\mathcal{H}^{i}(Rf_{0,*} \lc_\chi [d_M])$ sont purs d'après \cite[corollaire 5.4.4]{Ast100}. Il résulte alors de \cite[théorème 5.3.8]{Ast100} que les faisceaux pervers ${}^{\mathrm{p}}\mathcal{H}^{i}(Rf_{*} \lc_\chi [d_M])$ sont semi-simples et qu'on a 
        \begin{align}\label{eq:thm-dec}
          Rf_{*} \lc_\chi [d_M]=\bigoplus_{Z\in \Sc_\chi} \bigoplus_{i\in \ZZ}  \IC_Z(\fc_{Z,\chi}^i)[d_Z-i]      
        \end{align}
        où $\Sc_{\chi}$ est un ensemble fini  de sous-schéma intègres de $A$ (non vides) tels que, pour tout   $Z\in \Sc_{\chi}$,
        \begin{itemize}
        \item  pour tout $i\in \ZZ$,  le  système local $\ell$-adique $\fc_{Z,\chi}^i$ est semi-simple sur un ouvert dense de $Z$;
        \item  l'ensemble $I_{Z,\chi}=\{i\in \ZZ\mid \fc_{Z,\chi}^i\not=0\}$ est fini et non vide. 
        \end{itemize}
Un élément $Z$ de $\Sc_\chi$ est appelé un support de $Rf_{*} \lc_\chi [d_M]$.  Le but de cette sous-section est de prouver le théorème suivant.

\begin{theoreme}\label{thm:support}
    Supposons que la caractéristique de $\Fq$ est $>n$ et satisfait l'inégalité \eqref{eq:hyp-p}.  L'ensemble des supports  $\Sc_\chi$ est  le singleton formé du fermé $A_\Ga$.
        \end{theoreme}

        \begin{preuve}
          Elle occupe les §§ \ref{S:majoration} à \ref{S:conclusion}.
        \end{preuve}
        
      \end{paragr}

      \begin{paragr}[Majoration de l'amplitude.] ---  \label{S:majoration}On définit l'amplitude d'un support $Z$ dans $\Sc_\chi$ par:
        \begin{align*}
          \Amp(Z,\chi)=\max(I_{Z,\chi})-\min(I_{Z,\chi}).
        \end{align*}

        \begin{lemme}\label{lem:maj}
           Pour tout $Z\in \Sc_\chi$, on  a
          \begin{align*}
              \Amp(Z,\chi)\leq 2(d_Z+2d_f-d_M).
          \end{align*}
                  \end{lemme}

                  \begin{preuve}
                    Soit $Z\in \Sc_\chi$. Soit $i\in \ZZ$ tel que $\fc_{Z,\chi}^i\not=0$. Le système local  $\fc_Z^i$ apparaît comme facteur en degré $i-d_Z$ dans  $Rf_{*} \lc_\chi [d_M]$. Comme ce dernier a sa cohomologie concentrée en degrés compris entre $-d_M$ et $2d_f-d_M$, on obtient   $i\leq d_Z+2d_f-d_M$ et 
                    \begin{align}\label{eq:max}
                      \max(I_{Z,\chi})\leq d_Z+2d_f-d_M.
                    \end{align}
La dualité de Verdier commute à l'image directe vu que $f$ est propre;  ainsi on a
\begin{align*}
  D Rf_{*} \lc_\chi [d_M] = Rf_{*} \lc_{\chi^{-1}} [d_M].
\end{align*}
On obtient alors:
 \begin{align*}
          Rf_{*} \lc_{\chi^{-1}} [d_M]=\bigoplus_{Z\in \Sc_\chi} \bigoplus_{i\in \ZZ}  \IC_Z(\fc_{Z,\chi}^{i,\vee})[d_Z+i]
        \end{align*}
        Il s'ensuit que $\fc_{Z,\chi}^{i,\vee}$ apparaît en degré $-i-d_Z$ dans $Rf_{*} \lc_{\chi^{-1}} [d_M]$. Comme ci-dessus on a $-i\leq d_Z+2d_f-d_M$ et
            \begin{align}\label{eq:min}
                      -     \min(I_{Z,\chi}) \leq d_Z+2d_f-d_M.
                    \end{align}
Il suffit alors de combiner \eqref{eq:max} et \eqref{eq:min}.
      \end{preuve}
    \end{paragr}
    
    \begin{paragr}[Minoration de l'amplitude.] --- \label{S:minoration}Pour tout point fermé $a$ de $A$, soit $\delta^{\aff}_a$ la dimension de la partie affine du groupe algébrique $J_a$ défini au §  \ref{S:Ja}. L'application $a\mapsto \delta^{\aff}_a$ est semi-continue supérieurement. Pour tout sous-schéma intègre $Z\subset A$, on note $\delta^{\aff}_Z$ la valeur minimale de $ \delta^{\aff}_a$ pour $a\in Z$. On aura besoin de la variante suivante de \cite[proposition 7.3.2]{Ngo-LF}:

      \begin{proposition}\label{prop:minor}
           Pour tout $Z\in \Sc_\chi$, on  a
          \begin{align*}
              \Amp(Z,\chi)\geq 2(d_f-\delta^{\aff}_Z).
          \end{align*}
        \end{proposition}

        \begin{preuve}
          Il s'agit d'une variante de \cite[proposition 1.3]{Maulik-Shen-chi}, elle-même variante de \cite[proposition 7.3.2]{Ngo-LF}. Mais comme nous ne sommes pas exactement dans le cadre de  \cite{Maulik-Shen-chi} quelques mots sont dus au lecteur. Tout d'abord le triplet $(M,J_{X/A},A)$ est bien une fibration faiblement abélienne (\og\emph{weak abelian fibration}\fg) au sens  de  \cite[sec. 1.1]{Maulik-Shen-chi}: essentiellement   $J_{X/A}$ est un schéma en groupes lisses dont les fibres sont de dimension $d_f$, l'action de $J_{X/A}$ sur $M$ a des stabilisateurs affines et le module de Tate $T_{\Qlb}(J_{X/A})$ est polarisable (cf. \cite[sec. 4.12]{Ngo-LF} et \cite[section 3]{Cataldo-SLn} pour des explications détaillées). Notons que les hypothèses \cite[théorème 1.1  (a) et (c)]{Maulik-Shen-chi} sont bien satisfaites pour le faisceau $\kc=\lc_\chi$, cf. en particulier \eqref{eq:sum-fxpervers} et \eqref{eq:sum-fxpervers0}. Cependant il ne vérifie pas tout-à-fait  \cite[théorème 1.1, (b)]{Maulik-Shen-chi}. En fait,  cette hypothèse ne joue pas  de rôle dans la preuve de \cite[proposition 1.3]{Maulik-Shen-chi} de sorte que celle-ci s'applique \emph{verbatim} au faisceau $\lc_\chi$.
        \end{preuve}
      \end{paragr}

      \begin{paragr}[Conclusion.] --- \label{S:conclusion} En combinant la majoration du lemme \ref{lem:maj} et la minoration de la proposition \ref{prop:minor}, on voit que pour tout $Z\in S_\chi$, on doit avoir
        \begin{align}
          \label{eq:maj-delta}
          \delta_Z^{\aff}\leq d_A-d_Z.
        \end{align}
        En toute généralité, il existe $\underline m=(m_i)_{1\leq i\leq  r}$ et $\underline n=(n_i)_{1\leq i\leq  r}$ satisfaisant \eqref{eq:summi} tels que la strate $A_{\underline n,\underline m}$  définie au §  \ref{S:ell} rencontre $Z$ selon un ouvert. On peut alors trouver des ouverts $Z_i\subset A_{n_i}^{\el}$ de sorte que l'image de l'ouvert $Z_1\times \ldots\times Z_r$ par   $\iota_{\underline n,\underline m}$ soit incluse dans l'ouvert de $Z$ où $\delta_a^{\aff}$ prend sa valeur minimale. En suivant §  \ref{S:delta}, on a une fonction $a_i\mapsto \delta_{a_i}$ sur $Z_i$ où la fonction $\delta_{a_i}$ est définie relativement à $A_{n_i}$ et aux courbes spectrales de degré $n_i$.         Soit $\delta_i$  la valeur minimale prise par la fonction $a_i\mapsto \delta_{a_i}$ sur $Z_i$. On a donc
        \begin{align*}
          Z_i \subset A_{n_i}^{\el, \geq \delta_{i}}
        \end{align*}
        et par  l'hypothèse \eqref{eq:hyp-p} et la proposition \ref{prop:Severi} on a
        \begin{align*}
          d_{Z_i}\leq d_{A_{n_i}}-\delta_i.
        \end{align*}
        Par commodité, notons $d_i$ la dimension des fibres de la fibration de Hitchin en rang $n_i$ (donnée par la formule \eqref{eq:df} où l'on doit substituer $n_i$ à $n$). Soit $a_i\in Z_i$ et $a\in Z$ l'image de $(a_1,\ldots,a_r)$. En utilisant la proposition \ref{prop:devissage}, on a, en reprenant  ses notations,  $d_i=\dim(J_{a_i})$ et $d_i-\delta_i=\dim(\tilde J_{a_i})$: il vient alors 
                \begin{align*}
          \delta_Z^{\aff}=d_f+ \delta_{1}+\ldots+\delta_{r}-(d_1+\ldots+d_r).
        \end{align*}
On obtient ensuite 
 \begin{align*}
          d_A-d_Z &\geq d_A-(d_{Z_1}+\ldots+d_{Z_r})\\
                  &\geq (\delta_{1}+\ldots+\delta_{r})+ d_A-(d_{A_{n_1}}+\ldots+d_{A_{n_r}})\\
                  &\geq  \delta_Z^{\aff}+(d_A-d_f) -((d_{A_{n_1}}-d_1)+\ldots+(d_{A_{n_r}}-d_r))\\
   &\geq \delta_Z^{\aff}+ r-1+(n-(n_1+\ldots+n_r))(\deg(D)+2g_C-2)
 \end{align*}
 en vertu de \eqref{eq:dA-qa}. Cette dernière égalité n'est compatible avec \eqref{eq:degD} et \eqref{eq:maj-delta} que si $r=1$ et qu'on a égalité dans \eqref{eq:maj-delta} et donc dans le lemme \ref{lem:maj}. Cela implique d'une part que  $A^{\el}$  rencontre le support $Z$  et d'autre part que le faisceau de plus haut degré $R^{2d_f}f_* \lc_\chi$ admet sur l'ouvert elliptique un facteur direct un faisceau $\ell$-adique non trivial supporté sur  $Z\cap A^{\el}$. Pour conclure, il nous reste à déterminer le faisceau de plus haut degré sur l'ouvert elliptique. Pour cela, on note $f^{\el}$ la restriction de la fibration de Hitchin à l'ouvert elliptique $A^{\el}$. Soit $g$ la restriction de $f^{\el}$ à l'ouvert elliptique régulier $M^{\el,\reg}$ (c'est le lieu où les fibrés de Hitchin s'identifient par la correspondance rappelée au § \ref{S:Ja} à des fibrés en droites sur les courbes spectrales).  Comme dans \cite[proposition 6.5.1]{Ngo-LF}, on a, pour des raisons de dimension, un isomorphisme
 \begin{align*}
   R^{2d_f}g_! \lc_\chi \to R^{2d_f}f^{\el}_* \lc_\chi.
 \end{align*}
Puisque l'action de $J_{X/A^{\el}}$ est, fibre à fibre, transitive sur $M^{\el,\reg}$, on voit que le support du faisceau $R^{2d_f}f^{\el}_* \lc_\chi$ est  inclus dans $A_\Ga\cap A^{\el}$  et que sur $A_\Ga\cap A^{\el}$ ce faisceau est le faisceau constant $\Qlb$. On en déduit qu'on  a nécessairement $Z\cap A^{\el}=A_\Ga\cap A^{\el}$. Par conséquent,  on  a $A_\Ga\subset Z$. Mais, par ailleurs, on a nécessairement  $Z\subset A_\Ga$ par  la proposition \ref{prop:borne-supp}. On a donc bien $Z=A_\Ga$.
\end{paragr}

\subsection{Un second théorème de support}\label{ssec: second thm}

\begin{paragr}\label{S:C prime} On reprend les notations et les hypothèses des  sous-sections \ref{ssec:morphisme} et \ref{ssec:syst-loc}.  Soit $f'_0:M'_0\to A'_0$ la fibration de Hitchin relative à la courbe $C'_0$, au diviseur $D'_0=\rho_0^*D_0$, au rang $n'$ et au degré $e$ (cf. §§  \ref{S:C'} et \ref{S:morph}). Par la formule d'Hurwitz, le genre $g_{C'_0}$ de la courbe $C'_0$ est donné par 
  \begin{align}
    \label{eq:Hurwitz}
    2g_{C'_0}-2=d(2g_{C_0}-2).
  \end{align}
  On a également
  \begin{align}
    \label{eq:degD prime}\deg(D'_0)=d \deg(D_0).
  \end{align}
  Rappelons qu'on enlève l'indice $0$ pour noter les mêmes objets obtenus par changement de base à  $k$. Les égalités ci-dessus valent pour les objets sur $k$ (donc sans indice $0$).
\end{paragr}

\begin{paragr}
  Le groupe de Galois $\Gamma$ du revêtement $\rho_0:C_0'\to C_0$ agit naturellement sur $M_0'$.  Il résulte de la coprimalité de $n$ et $e$ que cette action est sans point fixe et que le quotient $M_0'/\Ga$ est lisse. Ce dernier s'identifie d'ailleurs par $(\ec,\theta)\mapsto  (\rho_{0,*}\ec,\rho_{0,*}\theta)$ à la sous-variété de $M_0$ des points fixes sous l'action par tensorisation d'un fibré en droite $\lc$ sur $C$ tel que
  \begin{align*}
    \rho_{0,*}\oc_{C_0'}=\oc_{C_0}\oplus \lc\oplus\lc^{\otimes 2}\oplus\ldots \oplus \lc^{\otimes d-1} \ \text{ et  }  \   \lc^{\otimes d}= \oc_{C_0}.
  \end{align*}
  Sur ces points le lecteur pourra consulter \cite[§ 2 et propositions 3.1 et  3.3]{Narasim-Ram}, points  qui sont d'ailleurs repris, dans notre  contexte,  dans  \cite[proposition 7.1]{Hausel-Tha}.
\end{paragr}

\begin{paragr} \label{S:mu chi} Le morphisme de Hitchin $f_0':M_0'\to A_0'$ est propre et $\Gamma$-équivariant;   le morphisme $\iota_0:A_0'\to A_0$ est fini et $\Ga$-invariant. Il s'ensuit que le composé $\iota_0\circ f_0'$ est propre et $\Ga$-invariant et se descend en un morphisme propre, noté $\bar f_0':M_0'/\Gamma \to A_0$.  Soit $\mu_\chi$ le système local sur $M_0'/\Ga$ qu'on obtient en poussant le revêtement $M_0'\to M_0'/\Ga$ par le  caractère  $\chi$ de $\Ga$.  Par abus, on note  encore   $\mu_\chi$  le $\Qlb$-faisceau lisse de rang $1$ sur  $M'/\Ga$ obtenu par extension  des scalaires.  Comme au § \ref{S:def support}, on a 

   \begin{align}\label{eq:thm-dec f'}
          R\bar f_{*}' \mu_\chi [d_{M'}]=\bigoplus_{Z\in \Sc_\chi} \bigoplus_{i\in \ZZ}  \IC_Z(\fc_{Z,\chi}^i)[d_Z-i]      
   \end{align}
   où les notations sont celles utilisées pour le second membre de \eqref{eq:thm-dec}.
\end{paragr}

\begin{paragr} Voici le second théorème de support que nous démontrerons.

  \begin{theoreme}\label{thm:support2}
 Supposons que la caractéristique de $\Fq$ est $>n$ et satisfait l'inégalité \eqref{eq:hyp-p}. L'ensemble $\Sc_\chi$  dans la décomposition \eqref{eq:thm-dec f'} est le singleton formé du fermé $A_\Ga$.
        \end{theoreme}

        \begin{preuve} 
   En utilisant la formule \eqref{eq:Hurwitz} d'Hurwitz et \eqref{eq:degD prime}, on voit que l'inégalité \eqref{eq:hyp-p} entraîne
         \begin{align*}
           p-1&>2+n'(2g_{C'}-2)+n'(n-1)\deg(D')\\
           &>2+n'(2g_{C'}-2)+n'(n'-1)\deg(D').
         \end{align*}
         Autrement dit on a l'égalité \eqref{eq:hyp-p} relative à $C',n'$ et $D'$.         On peut donc appliquer le théorème \ref{thm:support} à la fibration de Hitchin $f':M'\to A'$ et au système local trivial:   on en déduit que l'image directe $Rf'_*\Qlb[D_{M'}]$ a pour seul support la base $A'$ toute entière dans le théorème de décomposition. Le morphisme $\iota:A'\to A_\Ga$ est fini, surjectif entre variétés irréductibles: il  s'ensuit que  $R(\iota\circ f')_*\Qlb$  est également, à un décalage près, une somme de faisceau pervers semi-simples de support $A_\Ga$ tout entier, cf. par exemple \cite[lemme 3.5]{Maulik-Shen-IC} pour un argument. Pour conclure, il suffit d'observer que  $R\bar f_{*}' \mu_\chi [d_{M'}]$  est un facteur direct de $R(\iota\circ f')_*\Qlb$.
         
        \end{preuve}
        
\end{paragr}

\section{Comparaison de cohomologies relatives}\label{sec:compa}

\subsection{Calculs de traces de Frobenius}\label{ssec:traces}

\begin{paragr}  On conserve les notations et les hypothèses utilisées dans toute  la section \ref{sec:coh}. 
\end{paragr}

\begin{paragr} \label{S:laD}  Soit $F$ le corps de fonctions de la courbe $C_0$. Soit $|C_0|$ l'ensemble des points fermés de $C_0$. Pour tout $v\in |C_0|$ soit $F_v$ le complété de $F$ en $v$ et $\oc_v\subset F_v$ l'anneau des entiers. On fixe $\varpi_v\in \oc_v$ une uniformisante. On écrit $D_0=\sum_{v\in |C_0| }D_v[v]$. Soit $\AAA$ l'anneau des adèles de $F$. Soit  $\oc=\prod_{v\in |C_0|}\oc_v\subset \AAA$ le sous-anneau compact maximal. On définit le morphisme degré
  \begin{align*}
    \deg:\AAA^\times\to \ZZ
  \end{align*}
  par
  \begin{align*}
    \deg(x)=\sum_{v\in |C_0| } \deg(v) v(x)
  \end{align*}
  où l'on identifie $v$ à  la  valuation normalisée sur $F_v$. Soit $\AAA^0$ le noyau du morphisme $\deg$. Les groupes $J_0(\Fq)$ et  $P_0(\Fq)$ s'identifient respectivement aux groupes $F^\times\back \AAA^0/ \oc^\times$ et  $F^\times\back \AAA^\times/ \oc^\times$. Dans cette identification qui respecte les degrés, la classe du fibré  $\oc_{C_0}(D_0)$ correspond à  la classe de l'élément $\varpi_D\in  \AAA^\times$ défini par  $\varpi_D=(\varpi_v^{D_v})_{v\in |C_0|}$ et l'élément $\al\in P^1_0(\Fq)$ correspond à une classe encore notée $\al$. On identifie alors  le  caractère $\chi$ d'ordre $d$ de $J_0(\Fq)$ défini au § \ref{S:chi} à un caractère d'ordre $d$, encore noté $\chi$, du groupe  $ \AAA^0$.

  Soit $G=\GL(n)$. Pour tout $i\in \ZZ$, soit
\begin{align*}
    G(\AAA)^i=\{g  \in G(\AAA)\mid  \deg(\det(g))=i\}.
\end{align*}
Plus généralement pour un sous-groupe $R\subset G(\AAA)$ on notera $R^i=R\cap G(\AAA)^i$. Soit  $\ggo=\mathfrak{gl}(n)$ l'algèbre de Lie de $G$. Soit $\ggo(\oc)=\prod_{v\in |C_0|}\ggo(\oc_v)$ et $\ggo_{D_0}(\oc)=\varpi_D^{-1} \ggo(\oc)$. Soit $\mathbf{1}_{ \ggo_{D_0}(\oc)}$ la fonction  sur $\ggo(\AAA)$ caractéristique du sous-ensemble $\ggo_{D_0}(\oc)$.

\end{paragr}

\begin{paragr}[Un premier calcul de trace de Frobenius.] --- \label{S:chi a} Dans toute la suite, on fixe un élément  $a\in A_{\Ga,0}^{\el}(\Fq)$ où la partie localement fermée   $A_{\Ga,0}^{\el}\subset A_0^{\el}$ est définie au  § \ref{S:A chi ell}. Soit $F_a$   le corps de fonctions de la  courbe  spectrale $X_a$ et $\AAA_{F_a}$ l'anneau des adèles de $F_a$. Soit $\oc_{F_a}\subset \AAA_{F_a}$ le sous-anneau compact maximal. On a un morphisme norme $N_{F_a/F}: \AAA_{F_a}^\times \to \AAA^\times$. Le morphisme degré $\AAA_{F_a}^\times \to \ZZ$ est donné par la composition $\deg\circ N_{F_a/F}$. Soit $\AAA_{F_a}^0$ son noyau.  Vu qu'on a $a\in A_{\Ga,0}^{\el}(\Fq)$, la courbe $X_a$ est géométriquement connexe, cf. § \ref{S:ell}, et le morphisme degré est surjectif.   Soit  $x_a\in  \AAA_{F_a}^\times $ tel que $\deg(N_{F_a/F}(x_a))=1$. On étend alors le caractère $\chi$ de $\AAA^0$ en un caractère $\chi_a$ de $\AAA^\times$ d'ordre $d$  en posant
  \begin{align}
    \label{eq:chi a}
    \chi_a(y)=\chi(y N_{F_a/F}(x_a)^{-\deg(y)}) \ \ \ \forall y\in \AAA^\times.
  \end{align}
  Soit $\tilde J_a$ est la jacobienne de la normalisée $\tilde X_a$ de $ X_a$. Comme au § \ref{S:laD}, on a une identification entre les groupes $\tilde J_a(\Fq)$ et $F_a^\times\back \AAA_{F_a}^0/ \oc_{F_a}^\times$. On a  $a\in A_{\Ga,0}(\Fq)$; d'après la proposition \ref{prop:img-iota}, le  caractère $\chi\circ N_{X_a/C_0}$ est trivial sur $J_a(\Fq)$. À l'aide de la proposition \ref{prop:norme},  on en déduit que le  caractère $\chi\circ N_{\tilde X_a/C_0}$ est trivial sur $\tilde J_a(\Fq)$. Ce dernier s'identifie au caractère $\chi\circ N_{F_a/F}$ sur $\AAA_{F_a}^0$ qui est donc aussi trivial. Il s'ensuit le caractère $\chi_a$ sur $\AAA^\times$ d'une part ne dépend pas  du choix de $x_a$ et d'autre part est trivial sur $N_{F_a/F}(\AAA_{F_a}^\times)$.  On munit le  groupe  $\AAA_{F_a}^0$  de la mesure de Haar qui donne le volume $1$ au sous-groupe compact maximal $\oc_{F_a}^\times$. Pour la mesure quotient par la mesure de comptage sur $F_a^\times$, on a
  \begin{align}\label{eq:vol Fa}
    \vol(F_a^\times\back \AAA_{F_a}^0)=(q-1)^{-1}|\tilde J_a(\Fq)|.
  \end{align}
  On munit le  groupe $\AAA_{F_a}^\times $  de la mesure de Haar qui donne la mesure de comptage sur le quotient  $\AAA_{F_a}^\times /\AAA_{F_a}^0$.

  Soit  $U\in \ggo(F)$ un élément dont le polynôme caractéristique s'identifie à $a$.  Observons que $U$ est semi-simple régulier et qu'il est uniquement défini à $G(F)$-conjugaison près. Soit $G_U$ le centralisateur de $U$ dans $G$: c'est un sous-tore maximal. L'identification de  $U$ à un endomorphisme de $F^n$  munit $F^n$ d'une structure de $F_a$-module monogène. Il s'ensuit que le groupe $G_U(F)$ s'identifie au groupe multiplicatif $F_a^\times$. Dans l'identification $ G_U(\AAA) \simeq \AAA_{F_a}^\times $ qui en résulte, la restriction du morphisme déterminant  $\det:G(\AAA)\to \AAA^\times$ à $G_U(\AAA)$ s'identifie au morphisme norme $N_{F_a/F}$ sur  $\AAA_{F_a}^\times $. Il s'ensuit que  le caractère  $\chi_a\circ \det$ est trivial sur $G_U(\AAA)$.  Le noyau $G_U(\AAA)^0$ du morphisme $\deg\circ\det $ s'identifie alors à  $\AAA_{F_a}^0$.

On munit le groupe  $G(\AAA)$ de la mesure de Haar qui donne le volume $1$ au sous-groupe compact maximal $G(\oc)$. On transporte sur $G_U(\AAA)$  la mesure de Haar qu'on a fixée sur $\AAA_{F_a}^\times $. Vu \eqref{eq:vol Fa}, on a, pour la mesure quotient par la mesure de comptage sur $G_U(F)$,
  \begin{align}\label{eq:vol GX}
    \vol(G_U(F)\back G_U(\AAA)^0)=(q-1)^{-1}|\tilde J_a(\Fq)|.
  \end{align}

  Rappelons qu'on a fixé un degré $e\in \ZZ$ premier à $n$ pour définir l'espace $M$, cf. § \ref{S:Hitchin}. On peut alors formuler un premier calcul de trace de Frobenius.

  \begin{proposition}\label{prop:trace}
    On a
    \begin{align*}
 \trace(\Fr,(Rf_{*} \lc_\chi )_a)   =|\tilde J_a(\Fq)|\chi(\gamma_{a,\Fr})^e \int_{G_U(\AAA)\back G(\AAA)} \mathbf{1}_{\ggo_{D_0}(\oc)}(g^{-1}Ug)\chi_a(\det(g)) \, dg
    \end{align*}
    où l'on note $dg$ la mesure quotient et où  l'élément $\gamma_{a,\Fr}\in \Ga$ est défini  dans la proposition \ref{lem:gamma_a}. 
  \end{proposition}

 \begin{preuve} En utilisant le dictionnaire de Weil entre faisceaux localement libres sur $C_0$ et adèles, cf. \cite{scfh} pour une référence commode, on constate que l'ensemble  $M_{a,0}(\Fq)$ des points $\Fq$-rationnels de la fibre de Hitchin $M_{a,0}$ au-dessus de $a$ s'identifie à l'ensemble des orbites du groupe $G_U(F)$ agissant par translations à gauche sur l'ensemble
  \begin{align*}
  \Xgo_a=  \{ g\in G(\AAA)^e/G(\oc) \mid g^{-1}Ug \in  \ggo_{D_0}(\oc)\}.
  \end{align*}
D'après la formule des traces de Grothendieck-Lefschetz, cf. \cite[théorème 1.9]{SGA4.5}, on obtient:
   \begin{align*}
      \trace(\Fr,(Rf_{0,*} \lc_\chi )_a)=      \sum_{g\in G_U(F)\back \Xgo_a} \chi(\det(g)\al^{-e}).
   \end{align*}
Soit $g\in \Xgo_a$ et  $(\ec,\theta)$ le fibré de  Hitchin correspondant. Comme $a\in A_0^{\el}(\Fq)$, ce fibré est automatiquement stable. Il s'ensuit que le groupe $G_U(F)\cap g^{-1}G(\oc) g$, qui n'est autre que le groupe des automorphismes du fibré de Hitchin, s'identifie à $\Fq^\times$. On a donc en introduisant le caractère $\chi_a$ et en utilisant le lemme \ref{lem:x a}
     \begin{align*}
       \trace(\Fr,(Rf_{*} \lc_\chi )_a)&=    (q-1)  \sum_{g\in G_U(F)\back \Xgo_a} \frac{\chi(\det(g)\al^{-e})}{| G_U(F)\cap g^{-1}G(\oc) g   | }\\
                                         &=  (q-1)  \int_{G_U(F)\back G(\AAA)^e} \mathbf{1}_{ \ggo_{D_0}(\oc)}(g^{-1}Ug)\chi(\det(g)\al^{-e}) \, dg
     \end{align*}
     Comme on l'a observé plus haut, le caractère $\chi\circ\det$ est trivial sur $     G_U(\AAA)^0 $. L'expression ci-dessus s'écrit donc encore
     \begin{align*}
       &  |\tilde J_a(\Fq)| \int_{G_U(\AAA)^0 \back G(\AAA)^e} \mathbf{1}_{ \ggo_{D_0}(\oc)}(g^{-1}Ug)\chi(\det(g)\al^{-e}) \, dg\\
        &= \chi(N_{F_a/F}(x_a)\al^{-1})^e |\tilde J_a(\Fq)| \int_{G_U(\AAA)^0 \back G(\AAA)^e} \mathbf{1}_{ \ggo_{D_0}(\oc)}(g^{-1}Ug)\chi_a(\det(g)) \, dg\\
                                       &= \chi(N_{F_a/F}(x_a)\al^{-1})^e  |\tilde J_a(\Fq)| \int_{G_U(\AAA) \back G(\AAA)} \mathbf{1}_{ \ggo_{D_0}(\oc)}(g^{-1}Ug)\chi_a(\det(g)) \, dg.
     \end{align*}
    Par abus, on a noté de la même façon, à savoir $dg$, les diverses  mesures quotients qui apparaissent, le groupe $G_U(F)$ étant muni de la mesure de comptage. On conclut alors avec le lemme \ref{lem:x a} ci-dessous.

  \end{preuve}

  \begin{lemme}
    \label{lem:x a}
    On a
    \begin{align*}
     \chi(N_{F_a/F}(x_a)\al^{-1})= \chi(\gamma_{a,\Fr})
    \end{align*}
    où  $\gamma_{a,\Fr}\in \Ga$ est défini  dans la proposition \ref{lem:gamma_a}. 
  \end{lemme}

  \begin{preuve} Notons $\chi_\al$ le caractère d'ordre $d$ de $\AAA_F^\times$ donné par $\chi_\al(x)=\chi(x\al^{-\deg(x)})$ pour $x\in \AAA^\times$.  Soit $\tilde E$ le  corps de fonctions de  $\tilde C=\tilde C_{\al} \otimes_{\Fq}k$, cf. § \ref{S:torseur}. Soit $F\subset E, E_a\subset \tilde E$ les extensions de $F$ cycliques de degré $d$ tels que $\Ker(\chi_a)=F^\times N_{E_a/F}(\AAA_{E_a}^\times)$ et $\Ker(\chi)=F^\times N_{E/F}(\AAA_{E}^\times)$. Puisque $\chi_\al$ et $\chi_a$ coïncident sur les éléments de degré divisible par $d$, il en résulte qu'on a  $E_a\otimes_{\Fq}\Fqd= E\otimes_{\Fq}\Fqd$. Par le morphisme de réciprocité d'Artin, l'élément $N_{F_a/F}(x_a)\in \AAA^\times$ détermine un $F$-automorphisme $\Phi$ de $E_a\otimes_{\Fq}\Fqd= E\otimes_{\Fq}\Fqd$. Puisque $N_{F_a/F}(x_a)$ est de degré $1$, on a $\Phi=\delta\otimes  \Fr$ avec $\delta\in \Ga=\Gal(E/F)$. Ici $\delta$ est l'image de la classe de $N_{F_a/F}(x_a)$ dans $\AAA_F^\times/F^\times N_{E/F}(\AAA_E^\times)$ par l'application de réciprocité. On a donc $\chi_\al(\delta)=\chi_\al(N_{F_a/F}(x_a))$. On a une autre interprétation de $\delta$. En effet, on a une décomposition $a=\prod_{\sigma\in \Gal(E_a/F)} \sigma(b)$ avec  $b\in E_a[X]\subset (E\otimes_{\Fq} \Fqd)[X]$. La restriction de $\Phi$ à $E_a$ est triviale puisque par construction $\chi_a(N_{F_a/F}(x_a))=1$ et que la classe de $N_{F_a/F}(x_a)$ dans $\AAA_F^\times/F^\times N_{E_a/F}(\AAA_{E_a}^\times)$ est ainsi triviale. Il s'ensuit qu'on a $\Phi(b)=b$ dont on déduit $\delta(\Fr(b))=b$. Il est clair alors que $\delta$ s'identifie à l'élément  $\gamma_{a,\Fr}\in \Ga$ de la proposition \ref{lem:gamma_a}. 
  \end{preuve}
\end{paragr}

 \begin{paragr}[Un second  calcul de trace de Frobenius.] --- \label{S:Ca} On continue avec $a\in A_{\Ga,0}^{\el}(\Fq)$. Soit $i\in \ZZ$. On cherche à calculer:
   \begin{align*}
     \trace(\Fr,(R\bar f'_{*} \mu_{\chi^i} )_a).
   \end{align*}
D'après la formule des traces de Grothendieck-Lefschetz, on a :
    \begin{align}\label{eq:trace mu chi}
      \trace(\Fr,(R\bar f'_{*} \mu_{\chi^i} )_a)=\frac1{|\Ga|} \sum_{\ga\in\Ga} \sum_{a'}\chi(\gamma)^i |(M'_{a'}(k))^{\gamma\Fr}|
    \end{align}
    où la somme porte sur l'ensemble (fini) des éléments $a'\in A'(k)$ tels que $\iota_0(a')=a$ et $(M'_{a'}(k))^{\gamma\Fr}$ désigne l'ensemble des points dans $M'_a(k)$ fixes sous l'action du Frobenius tordu par $\gamma$. Il résulte de la proposition \ref{lem:gamma_a} que $(M'_{a'}(k))^{\gamma\Fr}$  est vide sauf si $\gamma=\gamma_{a,\Fr}$. D'après la proposition \ref{lem:gamma_a}  assertion 3, le groupe  $\Gamma$ agit simplement transitivement sur l'ensemble des éléments $a\in A'(k)$ tels que $\iota_0(a')=a$  et  $\gamma\in \Ga$ envoie bijectivement $(M'_{a'}(k))^{\gamma_{a,\Fr}\Fr}$ sur $(M'_{\gamma(a')}(k))^{\gamma_{a,\Fr}\Fr}$. On en déduit qu'on a 
    \begin{align}\label{eq:trace mu chi 2}
      \trace(\Fr,(R\bar f'_{*} \mu_{\chi^i} )_a)= \chi(\gamma_{a,\Fr})^i |(M'_{a'}(k))^{\gamma_{a,\Fr}\Fr}|
    \end{align}
    pour un élément quelconque $a'\in A'(k)$ tels que $\iota_0(a')=a$. Pour aller plus loin,   pour tout $\delta\in P_0^1(\Fq)$, on utilise  le revêtement  $\rho_{\chi,\delta}:C_{\chi,\delta}\to C_0$ étale, cyclique de groupe $\Ga$ qu'on a introduit au § \ref{S:C chi al}.  Rappelons que la construction ne dépend que de l'image de $\delta$ dans le quotient $P_0^1(\Fq)/\Ker(\chi)$. Le cas intéressant pour nous est $\delta=\gamma_{a,\Fr}\al$ pour lequel on pose $C_a=C_{\chi,\delta}$ et $\rho_a=\rho_{\chi,\delta}$. Alors $a'$ s'identifie à un élément $a''$ de $A_0''(\Fq)$ et  $(M'_{a'}(k))^{\gamma_{a,\Fr}\Fr}$ à $M''_{a'',0}(\Fq)$ où $M''_0$ est l'espace des fibrés de Hitchin définis relativement à la courbe $C_a$, au rang $n'$, au degré $e$ et au diviseur $D_0''=\rho_a^*D_0$ et $A_0''$ est la base correspondante.

    Comme dans la preuve de la proposition \ref{prop:trace}, on peut écrire ce comptage sous forme d'une intégrale adélique. Pour cela, on introduit $E_a$  le corps de fonctions de $C_{a}$,    l'anneau des adèles $\AAA_{E_a}$ de $E_a$ et $\oc_{E_a}\subset \AAA_{E_a}$ son sous-anneau compact maximal.

    \begin{lemme}
      \label{lem:sur chia}
      Le caractère $\chi_a$ de $\AAA^\times$ défini en \eqref{eq:chi a} a pour noyau $F^\times N_{E_a/F}(\AAA_{E_a}^\times)$ où $N_{E_a/F}$ désigne la norme.
    \end{lemme}

    \begin{preuve} Par la théorie du corps de classes, le sous-groupe $F^\times N_{E_a/F}(\AAA_{E_a}^\times)$ est d'indice $d$ dans $\AAA^\times$. Par ailleurs, le  caractère $\chi_a$ est d'ordre $d$. Il suffit donc de prouver que $F^\times N_{E_a/F}(\AAA_{E_a}^\times)\subset \Ker(\chi_a)$.  La restriction du  caractère $\chi_a$ à  $\AAA^0$ est le caractère $\chi$, qui a, par construction,  pour noyau $F^\times N_{E_a/F}(\AAA_{E_a}^0)$ où $\AAA_{E_a}^0\subset \AAA_{E_a}^\times$ est le sous-groupe des éléments de degré $0$. Soit $y\in \AAA^\times$ de degré $1$ dont la classe dans $\AAA^\times/ \Ker(\chi)$ correspond à  l'élément $\gamma_{a,\Fr}\al \in  P_0^1(\Fq)/\Ker(\chi)$. Comme $F^\times N_{E_a/F}(\AAA_{E_a}^\times)$ est engendré par $y$ et $F^\times N_{E_a/F}(\AAA_{E_a}^0)$ il suffit de prouver qu'on a $\chi_a(y)=1$. Mais c'est ce que l'on peut vérifier:
      \begin{align*}
        \chi_a(y)=\chi(y N_{F_a/F}(x_a)^{-1})= \chi(\gamma_{a,\Fr}\al N_{F_a/F}(x_a)^{-1})=1
      \end{align*}
      par le lemme \ref{lem:x a}.
    \end{preuve}

    Soit $H_a=\Res_{E_a/F}(\GL(n'))$ avec $n'd=n$ et $\hgo_a$ l'algèbre de Lie de $H_a$. Le groupe $H_a(\AAA)=\GL(n',\AAA_{E_a})$ est muni de la  mesure de Haar qui donne le volume $1$ au sous-groupe compact maximal $\GL(n',\oc_{E_a})$. Soit $V\in \hgo_a(F)$ de polynôme caractéristique $a''$ et  $H_{V}\subset H_a$ le centralisateur de $V$ dans $H_a$. Alors $H_V$ est un sous-tore maximal de $H_a$. On munit le groupe $H_{V}(\AAA)$  de la mesure de Haar qui donne le volume $1$ à l'unique sous-groupe compact maximal.     La courbe spectrale qui est le  revêtement de $C_a$ de degré $n'$ associée à $a''$, s'identifie à la courbe spectrale $X_a$, revêtement de $C$ de degré $n$ associé à $a$. On en déduit qu'on a 
 \begin{align}\label{eq:vol HV}
    \vol(H_V(F)\back H_V(\AAA)^0)=(q-1)^{-1}|\tilde J_a(\Fq)|
  \end{align}
où $H_V(\AAA)^0\subset H_V(\AAA)$ est le sous-groupe ouvert défini comme le noyau du morphisme $\deg\circ \det$.

    Comme au § \ref{S:laD}, on définit la fonction caractéristique   $\mathbf{1}_{\hgo_{a,D_0''}(\oc)}$ (ici relativement à la courbe $C_a$, au diviseur $D_0''$ et à l'entier $n'$).   En identifiant $\varpi_D$ à un élément de $\AAA_{E_a}^\times$ via le plongement $\AAA\hookrightarrow \AAA_{E_a}$,  on voit que la fonction $\mathbf{1}_{ \hgo_{a,D_0''}(\oc)  }$ s'identifie alors à la fonction sur  $\mathfrak{gl}(n',\AAA_{E_a})$ caractéristique de $\varpi_D^{-1}\mathfrak{gl}(n',\oc_{E_{a}})$. On obtient alors l'énoncé suivant:
    
 \begin{proposition}\label{prop:trace 2}
   On a
    \begin{align*}
      \trace(\Fr,(R\bar f'_{*} \mu_{\chi^i} )_a)   =|\tilde J_a(\Fq)|  \chi(\gamma_{a,\Fr})^i \int_{H_{V}(\AAA)\back H_a(\AAA)} \mathbf{1}_{\hgo_{a,D_0''}(\oc)}(h^{-1}Vh)\, dh
    \end{align*}
    où $dh$ est la mesure quotient.
  \end{proposition}

  \begin{remarque}
    L'espace $M'$ dépend du degré $e$ mais le second membre de l'égalité de la proposition \ref{prop:trace 2} n'en dépend visiblement pas. Dans la sous-section \ref{ssec: second thm}, on fait l'hypothèse $e$ premier au rang $n$ : elle n'est pas nécessaire pour la validité de la proposition \ref{prop:trace 2}.
  \end{remarque}
\end{paragr}

\subsection{Un lemme fondamental}\label{ssec:local}

\begin{paragr} On continue avec les notations de la section \ref{ssec:traces}. 
  On fixe aussi  un point fermé $v$ de $C_0$. Soit  $\chi_{a,v}$ le caractère de $F_v^\times$ obtenu comme  la composante en $v$ du caractère $\chi_a$ du  § \ref{S:chi a}. Notons que ce caractère est trivial sur le sous-groupe $\oc_v^\times$.
\end{paragr}

\begin{paragr} On a une décomposition $E_a\otimes_F F_v=\prod_{w|v} E_w$ indexée par les points fermés $w$ de $C_a$ au-dessus de $v$; le corps $E_w$ est le complété de $E$ en $w$. Il s'ensuit que $H_a(F_v)$ s'identifie au produit $ \prod_{w|v}\GL(n',E_w)$.  Soit  $\oc_{E_w}\subset E_w$ l'anneau de valuation. 

  Soit  $\ggo(\oc_{v})=\mathfrak{gl}(n,\oc_v)$ et $\hgo_a(\oc_v)=\prod_{w|v}\mathfrak{gl}(n',\oc_{E_w})$. Soit $\mathbf{1}_{\ggo(\oc_{v})}$ et $\mathbf{1}_{\hgo_a(\oc_{v})}$ les fonctions sur $\ggo(F_v)$ et $\hgo_a(F_v)$ caractéristiques des sous-algèbres  $\ggo(\oc_{v})$ et $\hgo_a(\oc_v)$ respectivement.
\end{paragr}

\begin{paragr} \label{S:image} Soit ${U_v}\in \ggo(F_v)$ un élément semi-simple régulier de $\ggo(F)$; par définition son polynôme caractéristique, qui est un élément de l'anneau de polynômes $F_v[t]$, est séparable. Soit ${V_v}\in \hgo_a(F_v)=\mathfrak{gl}(n',E_a\otimes_F F_v)$ et $P_{V_v}$ le polynôme caractéristique de ${V_v}$: c'est un polynôme unitaire de degré $n'$ à coefficients dans $E_a\otimes_F F_v$. Le groupe de Galois $\Ga$ agit sur $E_a\otimes_F F_v$ via son action sur $E_a$. On a un plongement $F_v \to E_a\otimes_F F_v$ donné par $1\otimes \Id_{F_v}$. De la sorte, on peut voir  $P_{U_v}$ comme un élément de $(E_a\otimes_F F_v)[t]$.   On dit que ${V_v}$ est une image de ${U_v}$ si le polynôme caractéristique de ${U_v}$ vérifie l'égalité suivante
  \begin{align*}
    P_{U_v}=\prod_{\sigma\in \Ga} \sigma(P_{V_v}).
  \end{align*}
Supposons que ${V_v}$ est une image de ${U_v}$. Alors  ${V_v}$ est aussi une image de $g^{-1}{U_v}g$ pour tout $g\in G(F_v)$.  Soit $G_{U_v}$ et $H_{V_v}$ les centralisateurs respectifs de ${U_v}$ et ${V_v}$ dans $G\times_F F_v$ et $H_a\times_F F_v$. Ce sont des sous-tores maximaux. On a des identifications naturelles
\begin{align*}
  H_{V_v}(F_v)\simeq ((E_a\otimes_F F_v)[t]/(P_{V_v}))^\times \simeq (F_v[t]/(P_{U_v}))^\times \simeq G_{U_v}(F_v).
\end{align*}
On pose
\begin{align}\label{eq:FT}
    \chi_{a,v}({U_v})=\chi_{a,v}^{-1}(\det(x,{U_v}x,\ldots,{U_v}^{n-1}x))
  \end{align}
  où, dans le membre de droite, le déterminant $\det$ est pris dans la base canonique de $F_v^n$ et $x\in F_v^n$ est tel que $(x,{U_v}x,\ldots,{U_v}^{n-1}x)$ soit une base de $F_v^n$. Comme ${U_v}$ est régulier semi-simple, de tels $x$ existent et forment une unique orbite sous l'action de $G_{U_v}(F_v)$. Le caractère $\chi_{a,v}$ est  trivial sur $N_{E_a/F}( (E_a\otimes_F F_v)^\times)$, cf. lemme \ref{lem:sur chia}. Il s'ensuit que le caractère $\chi_{a,v} \circ \det$ est trivial sur $G_{U_v}(F_v)$ et que $\chi_{a,v}({U_v})$ ne dépend pas du choix de $x$.  On vérifie qu'on a 
  \begin{align}\label{eq:invariance chia}
    \chi_{a,v}(g^{-1}{U_v}g)=\chi_{a,v}(\det(g))  \chi_{a,v}({U_v}).
  \end{align}
  
  On introduit aussi les discriminants
\begin{align*}
&  D^G_v({U_v})=\det(\ad({U_v})|\ggo(F_v)/\ggo_{U_v}(F_v))\in F_v^\times \\
  &D^{H_a}_v({V_v})=\det(\ad({V_v})|\hgo_a(F_v)/\hgo_{V_v}(F_v))\in F_v^\times,
\end{align*}
où $\ad$ désigne l'action adjointe. On note $|\cdot|_v$ la valeur absolue normalisée sur $F_v$. 

On munit $G(F_v), H_a(F_v)$, $G_{U_v}(F_v)$ et $H_{V_v}(F_v)$ des mesures de Haar qui donnent le volume $1$ à tout sous-groupe compact maximal. 
\end{paragr}

\begin{paragr}     Dans la section suivante, on utilisera le théorème suivant.
  
  \begin{theoreme} \label{thm:LF} Supposons $p>n$. Pour tout ${U_v}\in \ggo(F_v)$,  semi-simple régulier dans $\ggo(F)$, et ${V_v}\in \hgo_a(F_v)$ qui est une image de ${U_v}$, on a l'égalité
    \begin{align}
\label{eq:LF}      |D^G_v({U_v})|_v^{1/2} \int_{G_{U_v}(F_v)\back G(F_v)}\mathbf{1}_{\ggo(\oc_{v})}(g^{-1}{U_v}g)\chi_{a,v}(g^{-1}{U_v}g) \, dg\\
\nonumber      =  |D_v^{H_a}({V_v})|_v^{1/2}\int_{H_{V_v}(F_v)\back H_a(F_v)} \mathbf{1}_{\hgo_a(\oc_{v})}(h^{-1}{V_v}h)\, dh.
    \end{align}
où l'on note $dg$ et $dh$ les mesures quotients.
  \end{theoreme}

  \begin{preuve} On fixe une place  $w|v$. Soit $d'=[E_w:F_v]$. D'après le lemme \ref{lem:sur chia}, le caractère $\chi_{a,v}$ a pour noyau $N_{E_w/F_v}(E_w^\times)$: il est donc d'ordre $d'$. Soit $\Ga_w\subset \Ga$ le stabilisateur de $w$ dans $\Ga$: il s'identifie au groupe de Galois de l'extension $E_w/F_v$. Pour tout $w'|v$, on fixe $\sigma\in \Ga$ tel que $w'=\sigma(w)$ et on identifie $E_{w'}$ à $E_w$ par $\sigma^{-1}$. Soit $G'=\Res_{E_w/F_v} (\GL( \frac{n}{d'}))$. On identifie alors $H_a$ d'abord à  $\Res_{E_w/F_v} (\GL(n'))^{d/d'}$ puis à un sous-groupe de Levi noté $M'$ de $G'$. Soit $\mgo'$ l'algèbre de Lie de $M'$. L'élément $V_v$ s'identifie à un élément noté $V' \in \mgo'(F_v)$. Soit $P_{V'}$ le polynôme caractéristique de $V'$ vu comme endomorphisme de $E_w^{n/d'}$. On a donc $P_U=\prod_{\sigma\in \Ga_w} \sigma(P_{V'})$. De plus, par un argument standard de descente parabolique,  le membre de droite de l'égalité \eqref{eq:LF} est égal à l'expression
    \begin{align}\label{eq:IO G'}
      |D^{G'}({V'})|_v^{1/2}\int_{G'_{V'}(F_v)\back G'(F_v)} \mathbf{1}_{\ggo'(\oc_{v})}(h^{-1}{V'}h)\, dh,
    \end{align}
    où $G'_{V'}$ est le centralisateur de $V'$ dans $G'$ (c'est un sous-tore maximal) d'algèbre de Lie notée $\ggo'_{V'}$ et
    \begin{align*}
 D^{G'}({V'})=\det(\ad({V'})|\ggo'(F_v)/\ggo'_{V'}(F_v))\in F_v^\times .
    \end{align*}
    Dès lors, l'énoncé est une variante pour les algèbres de Lie du lemme fondamental pour l'induction automorphe. Il a été prouvé, pour les groupes sur les corps $p$-adiques,  par Waldspurger dans \cite{Wald-LF-induction}. La définition du facteur $\chi_{a,v}({U_v})$ est ici empruntée à  \cite[p.2]{Courtes}. On laisse le lecteur  le soin de vérifier la compatibilité des deux formulations.

    Esquissons une preuve dans notre contexte.     Tout d'abord, on se ramène aisément au cas où $U_v$ est de trace nulle. L'énoncé se ramène  alors à celui du lemme fondamental pour l'algèbre de Lie dans l'endoscopie ordinaire du groupe $\mathrm{SL}(n,F_v)$, énoncé démontré par Ngô dans \cite[théorème 1]{Ngo-LF}. En effet, le membre de gauche de  l'égalité \eqref{eq:LF} s'interprète comme une \og $\kappa$-intégrale orbitale\fg{} pour ce groupe. L'expression \eqref{eq:IO G'} quant à elle s'interprète comme une \og intégrale orbitale stable\fg{}  pour le groupe des éléments de  $\mathrm{SL}(n/d',E_w)$ dont la norme relative à $E_w/F_v$ de leur déterminant vaut $1$. Ce dernier groupe est un groupe endoscopique de $\mathrm{SL}(n,F_v)$.
    
  \end{preuve}
\end{paragr}

\subsection{Identification des cohomologies relatives}

\begin{paragr}
  On reprend les notations et les hypothèses de la section \ref{sec:coh}, cf. §§ \ref{S:C'} et \ref{S:al} en particulier. Soit $f':M'\to A'$ la fibration de Hitchin relative à la courbe $C'$, au diviseur $D'$ au rang $n'$ et au degré $e$ (premier à $n$) (cf. §§ \ref{S:C'} et \ref{S:morph}). On rappelle toutefois que tous les objets considérés proviennent par changement de base d'objets définis sur $\Fq$ et notés de la même façon avec un indice $0$, que $C_0'$ dépend du choix de $\al\in P^1(\Fq)$ et d'un caractère $\chi$ de $J_0(\Fq)$ d'ordre $d=n/n'$. Rappelons également qu'on a défini un élément $\varpi_D\in \AAA^\times$ au § \ref{S:laD}. Observons qu'on a
  \begin{align}\label{eq:signe}
    \chi_{}(\varpi_D \al^{-\deg(D)})^{n(n-1)/2}=\pm1. 
  \end{align}
Ce signe  est même $+1$ sous l'une des conditions suivantes :
\begin{itemize}
\item $n$ est impair;
\item  $n$ est pair et $2d|n $ où $d$, rappelons-le, est l'ordre de $\chi$ ;
\item   $D=2D_1$.
\end{itemize}
En particulier, on voit que le signe \eqref{eq:signe} est $+1$ si $d$ est impair.
On note $\eta_D$ le caractère quadratique de $\Gal(k/\Fq)$ qui envoie $\Fr$ sur  $\chi_{}(\varpi_D\al^{-\deg(D)})^{n(n-1)/2}$. On note par $(\eta_D)$ la torsion par ce caractère.  Le but de cette section est de prouver le théorème suivant:

  \begin{theoreme}\label{thm:iso-fx} Supposons que la caractéristique de $\Fq$ est $>n$ et satisfait l'inégalité \eqref{eq:hyp-p}. Il existe un isomorphisme entre les  semi-simplifications de faisceaux pervers gradués
    \begin{align*}
 \oplus_{i\in \ZZ}   {}^{\mathrm{p}}\mathcal{H}^{i}(Rf_{0,*} \lc_\chi )[2r](r)(\eta_D)  
    \end{align*}
    et
      \begin{align*}
   \oplus_{i\in \ZZ}   {}^{\mathrm{p}}\mathcal{H}^{i}(R\bar f'_{0,*} \mu_{\chi^{\tilde e}} )
      \end{align*}
      avec  $r=\dim(A_0)-\dim(A_0')$ et $\tilde e=e+\frac{n(n-1)}{2}\deg(D)$.
    \end{theoreme}

    \begin{remarque}
        En utilisant la formule pour les dimensions de $A_0$ et $A_0'$ ainsi que \eqref{eq:Hurwitz} et \eqref{eq:degD prime} on vérifie qu'on a 
    \begin{align*}
        r&=\frac{n(n-n')}{2}\deg(D)=(n')^2  \frac{d(d-1)}{2}\deg(D).
      \end{align*}   
    \end{remarque}
    
    Le théorème \ref{thm:iso-fx} est en fait équivalent au théorème suivant comme il résulte,  via le dictionnaire faisceaux-fonctions, de \eqref{eq:sum-fxpervers0} et de la pureté des faisceaux de cohomologie perverse.

    \begin{theoreme}\label{thm:egalite traces} Supposons que la caractéristique de $\Fq$ est $>n$ satisfait l'inégalité \eqref{eq:hyp-p}.       Pour tout $m\geq 1$ et tout $a\in A_0(\Fqm)$, on a 
 \begin{align}\label{eq:egalite traces}
   \trace(\Fr^m,(Rf_{*} \lc_\chi [2r](r)(\eta_D))_a)=\trace(\Fr^m, (R\bar f'_{*} \mu_{\chi^{\tilde e}})_a). 
    \end{align}      
           \end{theoreme}

L'hypothèse \eqref{eq:hyp-p}  sur la caractéristique entraîne la validité des théorèmes de support \ref{thm:support} et \ref{thm:support2}: la principale conséquence est la proposition suivante.

  \begin{proposition}\label{prop:utilisation support}
  Supposons que la caractéristique de $\Fq$ est $>n$ et satisfait l'inégalité \eqref{eq:hyp-p}.  Soit $\Om\subset A_{\Ga,0}$ un ouvert dense défini sur $\Fq$.  Supposons que l'égalité  \eqref{eq:egalite traces} du théorème \ref{thm:egalite traces} vaille pour tout $m\geq 1$ et tout $a\in \Om(\Fqm)$. Alors l'assertion du  théorème  \ref{thm:iso-fx} est vraie.
  \end{proposition}

  \begin{preuve}
    La proposition est une conséquence des théorèmes de support \ref{thm:support} et \ref{thm:support2}, cf. \cite[§  8.5]{Ngo-LF}, le lecteur pourra consulter aussi \cite[corollaire A.4]{Mig-SupportII}.
  \end{preuve}
\end{paragr}

\begin{paragr}[Fin de la preuve du théorème \ref{thm:iso-fx}.] --- Compte tenu de la proposition \ref{prop:utilisation support}, on va se contenter de prouver  l'égalité  \eqref{eq:egalite traces} du théorème \ref{thm:egalite traces} pour tout $m\geq 1$ et tout $a\in A_{\Ga,0}^{\el}(\Fqm)$. Pour alléger les notations, on se contente de rédiger  les égalités cherchées pour $m=1$.

  Soit $a\in A_{\Ga,0}^{\el}(\Fq)$ et $U\in \ggo(F)$ de polynôme caractéristique  $a$. Rappelons qu'on a défini un caractère $\chi_a$ de $\AAA^\times$ en \eqref{eq:chi a}. On utilisera l'égalité suivante, conséquence du lemme \ref{lem:x a}:
  \begin{align*}
    \chi_{a}(\varpi_D)=\chi_{\al}(\varpi_D ) \chi( \gamma_{a,\Fr})^{-\deg(D)}
  \end{align*}
où l'on pose  $\chi_\al(x)=\chi(x \al^{-\deg(x)})$ pour tout $x\in \AAA^\times$.

  En utilisant cette formule, la formule du produit ainsi que \eqref{eq:FT} et \eqref{eq:invariance chia}, on voit que, pour tout  $g=(g_v)_{v\in |C_0|}\in G(\AAA)$, on  a
  \begin{align*}
    &q^{-\deg(D)n(n-1)/2} \chi_\al(\varpi_D)^{n(n-1)/2}\chi_a(\det(g))\\
    &=  \chi( \gamma_{a,\Fr})^{\deg(D)n(n-1)/2}    \prod_{v\in |C_0|} |D^G_v(\varpi_v^{D_v}U)|_v^{1/2}\chi_{a,v}(\varpi_v^{D_v}g_v^{-1}Ug_v).
  \end{align*}
Il résulte de l'égalité ci-dessus et de la proposition \ref{prop:trace} que $\trace(\Fr,(Rf_{*} \lc_\chi )_a[2r](r)(\eta_D))$ est égale à 
  
    \begin{align}
    \label{eq:tr IO1}  q^{\deg(D)n(n'-1)/2} |\tilde J_a(\Fq)|\chi(\gamma_{a,\Fr})^{\tilde e}  \times  \\
\nonumber     \prod_{v\in |C_0|}  |D^G_v(\varpi_v^{D_v}U)|_v^{1/2}\int_{G_U(F_v)\back G(F_v)}\mathbf{1}_{\ggo(\oc_{v})}(g_v^{-1}\varpi_v^{D_v}Ug_v)\chi_{a,v}(g^{-1}\varpi_v^{D_v} Ug_v) \, dg_v.
    \end{align}

    Au § \ref{S:Ca}, on a introduit un groupe $H_a$, un élément  $V\in \hgo_a(F)$ et son centralisateur $H_V\subset H_a$. Pour tout $v\in |C_0|$, on  observe que $\varpi_v^{D_v}V$ est une image de $\varpi_v^{D_v}U$ au sens du § \ref{S:image}.  Il résulte du théorème \ref{thm:LF} que l'expression \eqref{eq:tr IO1} est égale à 
   \begin{align*}
  q^{\deg(D)n(n'-1)/2} |\tilde J_a(\Fq)|\chi(\gamma_{a,\Fr})^{\tilde e}  \times  \\
     \prod_{v\in |C_0|}  |D^{H_a}_v(\varpi_v^{D_v}V)|_v^{1/2}\int_{H_{V}(F_v)\back H_a(F_v)}\mathbf{1}_{\hgo_a(\oc_{v})}(h_v^{-1}\varpi_v^{D_v}Vh_v)\, dh_v.
    \end{align*}
En utilisant l'égalité
    
  \begin{align*}
    q^{-\deg(D)n(n'-1)/2}=  \prod_{v\in |C_0|} |D^{H_a}_v(\varpi_v^{D_v}V)|_v^{1/2},
  \end{align*}
  on voit que l'expression \eqref{eq:tr IO1} est aussi égale à
 \begin{align*}
 |\tilde J_a(\Fq)|\chi(\gamma_{a,\Fr})^{\tilde e}    \prod_{v\in |C_0|}  \int_{H_{a,V}(F_v)\back H_a(F_v)}\mathbf{1}_{\hgo_a(\oc_{v})}(h_v^{-1}\varpi_v^{D_v}Vh_v)\, dh_v
 \end{align*}
qui n'est autre que  $\trace(\Fr, R\bar f'_{*} \mu_{\chi^{\tilde e}})$ par  la proposition \ref{prop:trace 2}.
\end{paragr}

\section{Applications}\label{sec:applications}

\subsection{Préliminaires}\label{ssec:prelim}

\begin{paragr}
  Pour tout entier $g\geq 1$,  le corps  $\QQ(x_1,\ldots,x_g,z,t)$, où $x_1,\ldots,x_g,z,t$ sont des indéterminées, est muni:
  \begin{itemize}
  \item d'opérateurs d'Adams $\psi_n$, indexés par les entiers $n\geq 1$ et  définis par 
    \begin{align*}
      \psi_n(f(x_1,\ldots,x_g,z,t))=f(x_1^n,\ldots,x_g^n,z^n,t^n)
    \end{align*}
pour tout $f\in \QQ(x_1,\ldots,x_g,z,t)$;
\item d'une action du produit semi-direct $(\ZZ/2\ZZ)^g\rtimes \SG_g$: le groupe symétrique $\SG_g$ d'ordre $g!$ agit naturellement sur $(\ZZ/2\ZZ)^g$  et permute les variables $x_1,\ldots,x_g$ et le $i$-ième générateur canonique de  $(\ZZ/2\ZZ)^g$ envoie $x_i$ sur $zx_i^{-1}$. 
  \end{itemize}
  Un élément de $\QQ(x_1,\ldots,x_g,z,t)$ invariant par $(\ZZ/2\ZZ)^g\rtimes \SG_g$ est la fraction suivante:
   \begin{align*}
    Z_g=\frac{\prod_{i=1}^g(1-x_it)(1-x_i^{-1}zt)} {(1-t)(1-zt)}.
  \end{align*}
\end{paragr}

\begin{paragr}[Exponentielle pléthystique.] --- On étend l'opérateur $\psi_n$ à l'anneau des séries formelles $  \QQ(x_1,\ldots,x_g,z,t) [[T]]$ par la formule
  \begin{align*}
      \psi_n(\sum_{i=0}^\infty{a_i}T^i)=\sum_{i=0}^\infty\psi_n(a_i)T^{ni}.
  \end{align*}

        On définit alors une bijection
$$\mathrm{Exp} : T\QQ(x_1,\ldots,x_g,z,t) [[T]] \to  1+  T\QQ(x_1,\ldots,x_g,z,t) [[T]]$$
par $\mathrm{Exp}=\exp\circ \sum_{r\geq 1} \psi_r/r$ d'inverse $\mathrm{Log}= \sum_{r\geq 1}\frac{\mu(r)}{r} \psi_r \circ  \log$ où $\mu$ est la fonction de Möbius. Ici $\exp$ et $\log$ désignent les séries formelles usuelles exponentielle et logarithme.
  \end{paragr}

  \begin{paragr}[Partitions.] --- Soit $\La$ l'ensemble des partitions d'entiers: c'est l'ensemble des suites $\la=(\la_n)_{n\geq 1}$  à valeurs dans $\NN$, presque nulles et décroissantes. Pour tout $\la\in \La$, on pose $|\la|=\sum_{i=1}^\infty \la_i$. Pour tout entier $n$ soit $\La_n\subset\La$ le sous-ensemble des partitions $\la$ telles que $|\la|=n$. Pour tous $\la,\mu \in \La$, on définit $\la+\mu\in \La$ comme le seul élément de $\La$ tel que pour tout entier $i>0$ on a $|(\la+\mu)^{-1}(\{i\})   |=|\la^{-1}(\{i\})   | + |\mu^{-1}(\{i\}) |$. Pour tout entier $n\geq 1$ et $\la\in \La$, on pose $n\la=\la+\ldots+\la$ ($n$ fois). 
    On assimile un élément $\la\in \La$ à un diagramme de Young à $|\la|$ cases, où la première ligne contient $\la_1$ cases, la deuxième $\la_2$ cases, etc. La notation $s\in \la$ désigne une case $s$ du diagramme de Young de $\la$. Soit $s\in \la$. On note  $l(s)$ et $a(s)$  le nombre de cases qui sont respectivement strictement en-dessous et strictement à droite de  $s$. Pour la $j$-ième case de la $i$-ième ligne (avec par conséquent $j\leq \la_i$) on a $a=\la_i-j$ et $l$ est le nombre d'entiers $k>j$ tels que $\la_k\geq j$. On pose  $h(s)=1+a(s)+l(s)$. Notons que pour $s$ la dernière case de la dernière ligne on  a $a(s)=l(s)=0$ et $h(s)=1$. 
  \end{paragr}

  \begin{paragr}\label{S:Log}  Pour tous  $\la\in \La$ et $p\geq 0$, soit $\zc_{g,p,\la}\in \QQ(x_1,\ldots,x_g,z,t) $ défini par
             \begin{align*}
               \zc_{g,p,\la}=\prod_{s\in\la} (-t^{a(s)-l(s)}z^{a(s)})^p t^{(1-g)(2l(s)+1)}Z_g(x_1,\ldots,x_g,z,t^{h(s)}z^{a(s)}).
             \end{align*}

             Comme on l'a vu, il y a au moins $s\in \la$ tel que $a(s)=0$, $l(s)=0$ et $h(s)=1$. Le facteur correspondant à un tel $s$ est:
             \begin{align*}
               (-1)^p t^{1-g}Z_g.
             \end{align*}
             
             Soit     $\zc_{g,p}\in  \QQ(x_1,\ldots,x_g,z,t) [[T]]$ défini par
       \begin{align*}
                \zc_{g,p}=1+ \sum_{n\geq 1}(\sum_{\la\in \La_n}\zc_{g,p,\la} )T^{n} 
       \end{align*}
       On définit alors $\hc_{g,n,p}\in \QQ(x_1,\ldots,x_g,z,t) $ par l'égalité
       \begin{align*}
         \sum_{n\geq 1}  \hc_{g,n,p}T^n=(1-t)(1-zt)\Log\zc_{g,p}.
       \end{align*}
       Explicitement, on a
       \begin{align}
         \label{eq:sum hc la}
       \hc_{g,n,p}=\sum_{\la\in \La_n} \hc_{g,p,\la}
       \end{align}
       où, pour $\la_0\in \La$, on pose
\begin{align}
  \label{eq:Hofexplicit}
 \hc_{g,p,\la_0}= (1-t)(1-zt)\sum_{(m,r) }  \frac{\mu(r)}{r}  (-1)^{\sigma(m)} \sigma(m)! \prod_{\la\in \La} \frac{\psi_r(\zc_{g,p,\la})^{m_\la}}{m_\la !}
\end{align}
où la somme est prise sur l'ensemble (fini) des couples $(m,r)$ formés d'un élément $m$ de l'ensemble $\NN^{(\La)}$ des applications presque nulles de $\la$ dans $\NN$  et d'un entier naturel non nul $r$ tels que 
\begin{align}\label{eq:relation m r}
  r\sum_{\la \in \La}  m_\la \la=\la_0
\end{align}
et où l'on pose
\begin{align*}
  \sigma(m)=-1+\sum_{\la\in \La} m_{\la}.
\end{align*}
\emph{A priori},  $\hc_{g,n,p}$ est une fraction rationnelle dont le numérateur et le dénominateur sont tous deux dans  $\ZZ[x_1^{\pm1},\dots,x_g^{\pm1},z,t^{\pm1}]$ (et supposés sans facteur commun dans cet anneau factoriel). On voit facilement que $\prod_{i=1}^g(1-x_it)(1-x_i^{-1}zt)$ n'a pas de facteur commun avec le dénominateur de $\hc_{g,n,p}$ et divise son numérateur. Il résulte du  \cite[théorème 4.6]{MoOgor} qu'on a en fait $ \hc_{g,n,p}\in \ZZ[x_1^{\pm1},\dots,x_g^{\pm1},z,t^{\pm1}]$ et
\begin{align}
  \label{eq:une fraction}
  \frac{\hc_{g,n,p}}{\prod_{i=1}^g(1-x_it)(1-x_i^{-1}zt)}\in \ZZ[x_1^{\pm1},\dots,x_g^{\pm1},z,t^{\pm1}].
\end{align}
À la suite de Mozgovoy-O'Gorman,  on définit alors  $\mathbf{H}_{g,n,p}$ comme l'élément de $\ZZ[x_1^{\pm1},\dots,x_g^{\pm1},z]$  obtenu par évaluation en $t=1$ 
       \begin{align*}
                (-1)^{pn} z^{(g-1)n^2+p n(n+1)/2}  \hc_{g,n,p}.
       \end{align*}
 On renvoie le lecteur au théorème \ref{thm:MO} dû à  Mozgovoy-O'Gorman pour motiver l'introduction de ce polynôme.
       Il est clair que $\mathbf{H}_{g,n,p}$ est invariant sous l'action de $(\ZZ/2\ZZ)^g\rtimes \SG_g$. On a aussi que
       \begin{align*}
         z^{(g-1)n^2+p n(n+1)/2}\prod_{i=1}^g(1-x_i)(1-x_i^{-1}z)
       \end{align*}
       divise $ \mathbf{H}_{g,n,p}$ dans l'anneau $\ZZ[x_1^{\pm1},\dots,x_g^{\pm1},z]$.
     \end{paragr}

     \begin{paragr} \label{eq:notations n d}Soit $n,n',d\geq 1$ des entiers de sorte que $n=n'd$. Soit $g\geq 1$ et $g'=d(g-1)+1=g+(d-1)g$. On a donc $g'\geq 1$.
     \end{paragr}

     \begin{paragr}   \label{S:action mud}   Soit ${\mu}_d\subset \Qb$ le groupe des racines $d$-ièmes de l'unité. On définit une action du groupe $\mu_d$ sur l'anneau $\Qb[x_1^{\pm 1},\ldots,x_{g'}^{\pm 1},z]$
 \begin{align*}
 (\zeta,P)\in    \mu_d \times \Qb[x_1^{\pm 1},\ldots,x_{g'}^{\pm 1},z]\mapsto \zeta\cdot P \in \Qb[x_1^{\pm 1},\ldots,x_{g'}^{\pm 1},z]
  \end{align*}
       de la façon suivante:
       \begin{itemize}
\item  $\mu_d$ agit trivialement sur les variables $z,x_1,\ldots,x_g$;
  \item $\zeta\in \mu_d$  agit par $\zeta^i$ sur les variables $x_{i(g-1)+2},\ldots,x_{(i+1)(g-1)+1}$ pour $1\leq i <d$.
  \end{itemize}
 Pour tout $e\in \ZZ$,  soit   $\mathbf{H}_{g,n,p,d,e}\in \ZZ[x_1^{\pm 1},\ldots,x_{g'}^{\pm 1},z]$ défini par
  \begin{align}\label{eq:gnpd}
    \mathbf{H}_{g,n,p,d,e}=\frac1d\sum_{\zeta\in \mu_d} \zeta^e (\zeta\cdot \mathbf{H}_{g',n', dp}) .
  \end{align}
Ce polynôme ne dépend que de la classe de $e$ modulo $d$. Si $d=1$ on retrouve le polynôme   $\mathbf{H}_{g,n,p}$ défini ci-dessus.

Le polynôme   
  \begin{align*}
    z^{p+g-1} \prod_{i=1}^g (1-x_i)(1-zx_i^{-1})
  \end{align*}
  divise $\mathbf{H}_{g',n', dp}$ et est fixe sous l'action de $\mu_d$. Il s'ensuit qu'il existe $\tilde{\mathbf{H}}_{g,n,p,d,e}\in \ZZ[x_1^{\pm 1},\ldots,x_{g'}^{\pm 1},z]$   tel que
  \begin{align*}
    \mathbf{H}_{g,n,p,d,e}=z^{p+g-1} \prod_{i=1}^g (1-x_i)(1-zx_i^{-1})\tilde{\mathbf{H}}_{g,n,p,d,e}.
  \end{align*}
     \end{paragr}

     \begin{paragr}[Spécialisations.] --- \label{S:specialisation}On étudie certaines spécialisations des   polynômes précédents. On travaille sur le corps $\QQ(\xi, u)$.  Soit $\mathbf{H}_{g,n,p,d}^\flat\in \ZZ[\xi^{\pm 1},u]$ défini comme la spécialisation de $\mathbf{H}_{g',n', dp}$ en $z=u^2$, $x_i=u$ pour $1\leq i \leq g$ et $x_j=\xi^i u$ pour $1\leq i \leq d-1$ et $i(g-1)+2\leq j\leq (i+1)(g-1)+1$. Notons que l'invariance de $\mathbf{H}_{g',n',pd}$ est invariant sous l'action de $(\ZZ/2\ZZ)^{g'}\rtimes \SG_{g'}$ entraîne bien qu'on a $\mathbf{H}_{g,n,p,d}^\flat\in \ZZ[\xi^{\pm 1},u]$. Soit   $\tilde{\mathbf{H}}_{g,n,p,d,e}^\flat\in \ZZ[ u]$ la spécialisation de $\tilde{\mathbf{H}}_{g,n,p,d,e}$ en $z=u^2$ et $x_i=u$ pour tout $1\leq i\leq g'$. Alors on a
\begin{align}\label{eq:identite clef}
        \frac1d  \sum_{\zeta\in \mu_d}  \zeta^e \mathbf{H}_{g,n,p,d}^\flat(\zeta,u)=u^{2(p+g-1)} (1-u)^{2g} \tilde{\mathbf{H}}^\flat_{g,n,p,d,e}.
         \end{align}       
     \end{paragr}

       \begin{paragr}                Dans la suite, on utilisera le lemme  suivant:

         \begin{lemme}\label{lem:calcul H sur 1-u} Sous les notations de § \ref{eq:notations n d}, on a:
           \begin{enumerate}
           \item Le polynôme $\mathbf{H}_{g,n,p,d}^\flat(\xi,u) $ est divisible par $u^{2(p+g-1)}(1-u)^{2g} $ dans l'anneau $\ZZ[\xi^{\pm1},u]$.
           \item La valeur en $u=1$ de
             \begin{align*}
               \frac{ \mathbf{H}_{g,n,p,d}^\flat(\xi,u) }{(1-u)^{2g} }
             \end{align*}
             est égale, 
\begin{enumerate}
\item  si $g\geq 2,$ à
  \begin{align*}
    (-1)^{p(n-d) }\mu(n') (n')^{2g-3}     \prod_{i=1}^{d-1}((1-\xi^{in'} )(1-\xi^{-in'} ))^{g-1}
      \end{align*}
\item  si $g=1,$ à
  \begin{align*}
    &  2 \text{ si } pd  \text{ est impair et si } n' \equiv 2\mod  4\  ;\\
    &   1 \text{ sinon }.
  \end{align*}
  \end{enumerate}
\end{enumerate}
\end{lemme}

\begin{preuve}
  1. On a vu que le polynôme $z^{p+g-1} \prod_{i=1}^g (1-x_i)(1-zx_i^{-1})$ divisait $\mathbf{H}_{g',n', dp}$ dans l'anneau $\ZZ[x_1^{\pm 1},\ldots,x_{g'}^{\pm 1},z]$.  Par spécialisation, on obtient la divisibilité cherchée.

  2. La valeur cherchée est la spécialisation du polynôme, cf. \eqref{eq:une fraction},
  \begin{align*}
    (-1)^{pn} z^{(g'-1)(n')^2+ pn(n'+1)/2}  \frac{\hc_{g',n',dp}}{\prod_{i=1}^{g}(1-x_it)(1-x_i^{-1}zt)} \in \ZZ[x_1^{\pm 1},\ldots,x_{g'}^{\pm 1},z,t^{\pm 1}]
  \end{align*}
 en $t=1$, $z=1$, $x_i=1$ pour $1\leq i \leq g$ et $x_j=\xi^i $ pour $1\leq i \leq d-1$ et $i(g-1)+2\leq j\leq (i+1)(g-1)+1$. 

  Soit $\hc^\flat_{g',n',pd}$, $Z_{g'}^\flat$  et   $\zc_{g',pd,\la}^\flat$ pour $\la\in \La$ les éléments de $\QQ(\xi,t)$ obtenus comme les spécialisations respectives  de $\hc_{g',n',pd}$,  $Z_{g'}$ et  $\zc_{g',pd,\la}$ en $z=1$, $x_i=1$ pour $1\leq i \leq g$ et $x_j=\xi^i $ pour $1\leq i \leq d-1$ et $i(g-1)+2\leq j\leq (i+1)(g-1)+1$. Ainsi, on a 
    \begin{align*}
      Z_{g'}^{\flat}&= (1-t)^{2g -2} \prod_{i=1}^{d-1}((1-\xi^i t)(1-\xi^{-i} t))^{g-1}.
    \end{align*}

  L'expression qu'on cherche à calculer est alors la valeur en  $t=1$ de l'expression
    \begin{align}\label{eq:le rapport}
(-1)^{pn}     \frac{\hc_{g',n',pd}^\flat}{ (1-t)^{2g} }.
    \end{align}
     Soit une partition $\la_0\in \La_{n'}$. On regarde le terme associé à $\la_0$ dans la somme  qui définit $\hc_{g',n',pd}^\flat$, cette somme étant obtenue après  spécialisation de \eqref{eq:sum hc la}. On cherche alors à évaluer en $t=1$ l'expression
       \begin{align}\label{eq:Hofexplicit-bis} 
      \frac{\hc_{g',pd,\la_0}^\flat}{ (1-t)^{2g} }=  \sum_{(m,r) }  \frac{\mu(r)}{r}  \frac{(-1)^{\sigma(m)} \sigma(m)! }{(1-t)^{2g-2}}\prod_{\la\in \La} \frac{\psi_r(\zc_{g',pd,\la})^{m_\la}}{m_\la !},
       \end{align}
       obtenue par  spécialisation de \eqref{eq:Hofexplicit}.

         Supposons d'abord $g\geq 2$.    Il s'ensuit que dans   \eqref{eq:Hofexplicit-bis} tous les termes s'annulent en $t=1$ sauf celui associé au couple $(m,n')$ où $m$ est l'application qui vaut $1$ sur l'unique  partition de $1$ et $0$ ailleurs.  Ce couple n'apparaît que pour la partition $(1,\ldots,1)$ de $n'$.   On conclut que la valeur  $t=1$ du rapport \eqref{eq:le rapport}  est la valeur en $t=1$ de 

      \begin{align*}
  (-1)^{p(n-d) }   \frac{\mu(n')}{n'}\cdot t^{n'(1-g')} \frac{ (1-t^{n'})^{2g-2 }}{ (1-t)^{2g-2 } }   \prod_{i=1}^{d-1}((1-\xi^{in'} t^{n'})(1-\xi^{-in'} t^{n'}))^{g-1}
      \end{align*}
      c'est-à-dire
      \begin{align*}
        (-1)^{p(n-d) }\mu(n') (n')^{2g-3}     \prod_{i=1}^{d-1}((1-\xi^{in'} )(1-\xi^{-in'} ))^{g-1}
      \end{align*}
      On en déduit 2.a.
      Montrons maintenant 2.b. On suppose $g=1$; on a donc $g'=1$. On a alors $Z^\flat_{g'}=1$. Pour tout $\la\in \La$ et $p\geq 0$, on a
      \begin{align*}
        \zc_{1,pd,\la}=(-1)^{pd|\la|} \prod_{s\in\la} t^{(a(s)-l(s))pd} \\
             \end{align*}
      On en déduit que pour tout $\la_0\in \La$,  l'évaluation en $t=1$  de
      \begin{align}\label{eq:rapport hc1}
       \frac{\hc_{1,pd,\la_0}^\flat}{ (1-t)^{2} }
      \end{align}
      est égale à 
\begin{align}
  \label{eq:Hofexplicit-ter}
\sum_{(m,r) }  \frac{\mu(r)}{r}  (-1)^{\sigma(m)} \sigma(m)!  (-1)^{pd|\la_0|/r}\prod_{\la\in \La} \frac{1}{m_\la !}
\end{align}
où les sommes sont comme en \eqref{eq:Hofexplicit}.  Par conséquent,  en reprenant les définitions du § \ref{S:Log}, on voit que l'évaluation en $t=1$ de 
   \begin{align*}
     (-1)^{pn} \frac{\hc_{1,n',pd}^\flat}{ (1-t)^{2} }
       \end{align*}   
       est  $(-1)^{pn}$ fois le coefficient de degré $n'$ de la série
       \begin{align*}
         \Log(1+\sum_{k\geq 1} |\La_k| (-1)^{kpd} T^k)=\Log( \prod_{k\geq 1}  \frac{1}{1- (-1)^{kpd}T^k}).
       \end{align*}
Pour finir le calcul, on utilise l'identité élémentaire $\Exp(T^i)=(1-T^i)^{-1}$.       Distinguons deux cas suivant la parité de $pd$. Commençons par le cas où $pd$ est pair. Dans ce cas, on a $ (-1)^{pn}=1$. On calcule
       \begin{align*}
         \Log( \prod_{k\geq 1}  \frac{1}{1- T^k})&=\sum_{k\geq 1} \Log(  \frac{1}{1- T^k})\\
         &=\sum_{k\geq 1} T^k.
       \end{align*}
       Supposons ensuite    $pd$ impair.  On a alors  $ (-1)^{pn}=(-1)^{n'}$.
       On calcule
       \begin{align*}
         \Log( \prod_{k\geq 1}  \frac{1}{1- (-T)^k})&=\sum_{k\geq 1} \Log(  \frac{1}{1- T^{2k}})+\sum_{k\geq 0} \Log(  \frac{1}{1+ T^{2k+1}})  \\
                                                    &=\sum_{k\geq 1} T^{2k}+\sum_{k\geq 0} \Log(  \frac{1-T^{2k+1}}{1-T^{4k+2}}) \\
         &=\sum_{k\geq 1}T^{2k}+\sum_{k\geq 0 }T^{4k+2} - \sum_{k\geq 0 }T^{2k+1}.
       \end{align*}
       On conclut facilement.
\end{preuve}

     \end{paragr}

       \subsection{Comptage des fibrés de Hitchin de déterminant fixé}\label{ssec:det fixe}

       \begin{paragr}\label{S:n-torsion}
         On reprend les notations de la section \ref{sec:coh}. Soit un corps fini $\Fq$ ainsi qu'une clôture algébrique $k$ de $\Fq$.  Soit   $C_0$  une courbe projective, lisse, géométriquement connexe sur un corps fini $\Fq$ de genre noté $g$. On suppose dans la suite qu'on a $g\geq 1$. Soit $D_0$ un diviseur sur $C_0$ de degré $>2g-2$. On pose
         \begin{align*}
           p=\deg(D_0)-(2g-2).
         \end{align*}
         Rappelons qu'on suppose $p>0$.          Soit $n\geq 1$ un entier. On suppose que la caractéristique du corps $k$ est supérieure strictement à $n$ et qu'elle vérifie l'inégalité \eqref{eq:hyp-p}. Le sous-groupe $J_0(k)[n]\subset J_0(k)$ des points de $n$-torsion de la jacobienne $J_0$ de $C_0$ est alors isomorphe à $(\ZZ/n\ZZ)^{2g}$. On suppose que les points  de $n$-torsion sont tous définis sur $\Fq$ c'est-à-dire que  $J_0(\Fq)[n]$ est de cardinal $n^{2g}$. Quitte à faire un changement de base à une extension finie de $\Fq$, on peut toujours se ramener à cette situation.

         Pour un groupe fini $G$, on appelle, dans la suite, caractère de $G$ tout  homomorphisme $G \to \Qlb^\times$. Pour tout entier $m\geq 1$, soit $Q_m$ le groupe des caractères de  $J_0(\Fqm)$ d'ordre divisant $n$. C'est encore le groupe des caractères du quotient $J_0(\Fqm)/J_0(\Fqm)^{\otimes n}$ où  $J_0(\Fqm)^{\otimes n}$ est l'image du morphisme $J_0(\Fqm)\to J_0(\Fqm)$ donné par $\delta\mapsto \delta^{\otimes n}$.  On pose $Q=Q_1$ pour $m=1$. L'hypothèse sur les points de $n$-torsion entraîne le  lemme suivant. 

         \begin{lemme}\label{lem:ntorsion}
           Supposons qu'on a  $|J_0(\Fq)[n]|=n^{2g}$. Alors pour tout $m\geq 1$, le groupe  $Q_m$ est isomorphe à $(\ZZ/n\ZZ)^{2g}$ et le morphisme norme $N_m:J_0(\Fqm)\to J_0(\Fq)$ induit un isomorphisme $Q\to Q_m$.
         \end{lemme}

         \begin{preuve}
     Soit   $m\geq 1$.      La suite exacte courte associée à l'isogénie de $J_0$ donnée par $\delta \mapsto \delta^{\otimes n}$
           \begin{align*}
             0\to J_0(k)[n]\to J_0(k) \to J_0(k)\to 0 
           \end{align*}
           donne une suite exacte en cohomogie
           \begin{align*}
              0\to J_0(\Fqm)[n]\to J_0(\Fqm) \to J_0(\Fqm)\to H^1(\Gal(k/\Fqm),J_0(k)[n]) \to H^1(\Gal(k/\Fqm),J_0(k)).
           \end{align*}
           Comme on a  $H^1(\Gal(k/\Fqm),J_0(k))=\{0\}$ (lemme de Lang) et
           \begin{align*}
             H^1(\Gal(k/\Fqm),J_0(k)[n])=\Hom_{\textrm{cont}}(\Gal(k/\Fqm),J_0(k)[n])\simeq J_0(k)[n],
           \end{align*}
           on voit que $Q_m$ est isomorphe au groupe des caractères de $J_0(k)[n]$:  comme groupe il est isomorphe  à $(\ZZ/n\ZZ)^{2g}$ donc d'ordre est $n^{2g}$. La norme induit un morphisme surjectif $N_m:J_0(\Fqm)\to J_0(\Fq)$. On dispose donc dualement d'un morphisme $Q\to Q_m$ injectif donc bijectif pour des raisons de cardinalité.
         \end{preuve}
       \end{paragr}

       \begin{paragr}[Les valeurs propres de Frobenius.] --- \label{S:vp Frob}  Soit $\al\in P_0^1(\Fq)$. Pour  $\chi\in Q$, soit $d_\chi$ son ordre et  $n_\chi$ tel que $n=n_\chi d_\chi$. Comme au § \ref{S:C chi al}, on dispose de revêtements galoisiens $\tilde \rho_\al:\tilde C_\al \to C_0$ et  $\rho_{\chi,\al}:C_{\chi,\al} \to C_0$ de groupes de Galois respectifs $J_0(\Fq)$ et $\Ga_\chi=J_0(\Fq)/\Ker(\chi)$. Soit  $0\leq i\leq d_\chi-1$ et  $\fc_{\chi^i,\al,0}$ le $\Qlb$-faisceau lisse de rang $1$ sur $C_0$ qu'on obtient en poussant le revêtement $\tilde\rho_\al$  par $\chi^{-i}$.  Soit  $\fc_{\chi^i,\al}$ le $\Qlb$-faisceau lisse de rang $1$ sur $C$ qui s'en déduit. On a
         \begin{align*}
           H^1(C_{\chi,\al} \times_{\Fq}k,\Qlb)=\bigoplus_{i=0}^{d-1} H^1(C,\fc_{\chi^i,\al}).
         \end{align*}
  Les groupes de cohomologie  $H^1(C,\fc_{\chi^i,\al})$    sont munis d'une action de $\Fr$ dont on notera ainsi  les valeurs propres $\Fr$:
  \begin{itemize}
  \item  $\la_1,\ldots,\la_g,q\la_1^{-1},\ldots, q\la_g^{-1}$ pour $i=0$;
  \item $\la_{(2i-1)g+3-2i},\ldots,\la_{(2i+1)g  -2i}$ pour $1\leq i <d_\chi/2$;
  \item $\la_{(d_\chi-1)(g-1)+2},\ldots,\la_{d_\chi(g-1)+1},q\la_{(d_\chi-1)(g-1)+2}^{-1},\ldots,q\la_{d_\chi(g-1)+1}^{-1}$  si $d_\chi$ est pair et $i=d_\chi/2$.
  \end{itemize}
Notons que, pour $1\leq i <d_\chi/2$ les valeurs propres de  $\Fr$ agissant sur $H^1(C,\fc_{\chi^{-i},\al})$ sont
\begin{align*}
  q\la_{(2i-1)g+3-2i}^{-1},\ldots,q\la_{(2i+1)g  -2i}^{-1}.
\end{align*}
Soit $g_\chi=d_\chi(g-1)+1$. Pour tout entier $m\geq 1$ et tout $e\in \ZZ$, on note alors
\begin{align*}
  \psi_m(\mathbf{H}_{g,n,p,d_\chi,e})[\chi,\al]
\end{align*}
l'évaluation du polynôme $\psi_m(\mathbf{H}_{g,n,p,d_\chi,e})\in \ZZ[z,x_1^{\pm 1},\ldots,x_{g_\chi}^{\pm 1}]$, cf.  \eqref{eq:gnpd}, en $z=q$, $x_1=\la_1,\ldots, x_{g_\chi}=\la_{g_\chi}$. La même notation vaut pour le polynôme $\tilde{\mathbf{H}}_{g,n,p,d_\chi,e}$.

Si $d_\chi=1$, le caractère $\chi$ est le caractère trivial noté $1$. Dans ce cas,  $\psi_m(\mathbf{H}_{g,n,p,1,e})[1,\al] $ ne dépend que des valeurs propres de $\Fr$ sur $H^1(C,\Qlb)$ et ne dépend donc pas de $\al$ (ni d'ailleurs de $e$). On pourra noter simplement $\psi_m(\mathbf{H}_{g,n,p})[1]= \psi_m(\mathbf{H}_{g,n,p,1,e})[1,\al] $.
\end{paragr}
      
       \begin{paragr}[Comptage des fibrés de Hitchin de degré fixé.] ---    Soit $e\in \ZZ$ premier à $n$.  Soit $M_n^e(C_0,D_0)$  l'espace des fibrés de Hitchin de  rang $n$ et degré $e$ sur la courbe $C_0$ relatif au diviseur $D_0$.  On a le résultat suivant dû à Mozgovoy et  O'Gorman.

\begin{theoreme}(\cite[théorème 1.1]{MoOgor})\label{thm:MO}
  Pour tout entier $m\geq 1$ et tout $e\in \ZZ$ premier à $n$, on a
       \begin{align}\label{eq:Mozgovoy}
         |  M_n^e(C_0,D_0)(\Fqm) |&=\psi_m(\mathbf{H}_{g,n,p})[1].
         \end{align}
\end{theoreme}
       \end{paragr}
       
       \begin{paragr}[Comptage des fibrés de Hitchin de déterminant fixé.] --- Soit $\be\in P_0^e(\Fq)$.   On définit  $M_n^\be(C_0,D_0)$  comme le fermé de $M_n^e(C_0,D_0)$ formé des fibrés de Hitchin $(\ec,\theta)$ tels que le déterminant de $\ec$ soit le fibré en droites $\be$. C'est un $k$-schéma lisse sur $k$. Tout élément $\chi\in Q$, vu comme caractère de $J_0(\Fq)$,  s'étend alors de manière unique en un caractère noté $\chi_\al$ de $P_0(\Fqm)$ qui est trivial sur $\al$. Le théorème suivant fournit une formule explicite pour le comptage $|  M_n^\be(C_0,D_0)(\Fqm)|$ pour tout entier $m\geq 1$. 
         
         \begin{theoreme}\label{thm:comptage det fix}
           Supposons que la caractéristique de $\Fq$ est $>n$ et vérifie l'inégalité \eqref{eq:hyp-p}. Pour tout entier $m\geq 1$ et tout $\be\in P_0(\Fq)$ de degré premier à $n$, on a
      \begin{align*}
        |M_n^\be(C_0,D_0)(\Fqm)|&=  \frac1{|J_0(\Fqm)|} \sum_{\chi\in Q} \left( q^{r_\chi \deg(D)} \chi_\al^{-1}(\be \varpi_D^{n(n-1)/2} ) \right)^m\psi_m(\mathbf{H}_{g,n,p,d_\chi,\tilde e})[\chi,\al] 
      \end{align*}
      où l'on note
      \begin{itemize}
      \item $\varpi_D$ est un élément de $\AAA^\times$ dont la classe dans $F^\times\back \AAA^\times/\oc^\times$ correspond à $D_0$, cf. § \ref{S:laD}.
      \item $D$ est le diviseur sur $C$ obtenu par changement de base de $D_0$; on a d'ailleurs $\deg(D)=\deg(D_0)$;
      \item    $\tilde e=\deg(\be)+\frac{n(n-1)}{2}\deg(D)$        ;
      \item $p=\deg(D)-(2g-2)$;
      \item $n_\chi d_\chi=n$ où $d_\chi$ est l'ordre de $\chi$  et $r_\chi =n_\chi^2 d_\chi(d_\chi-1)/2$.
              \end{itemize}
   \end{theoreme}
   
La démonstration se trouve au § ci-dessous.
         
       \end{paragr}

       \begin{paragr}[Preuve du théorème \ref{thm:comptage det fix}.]  --- Soit $e\in \ZZ$ premier à $n$. Pour tout $\chi\in Q$, soit $\lc_{\chi,\al}$ le système local de rang $1$ sur $M_n^e(C_0,D_0)$ défini en \eqref{eq:L chi al}. Par abus, on note encore $\lc_{\chi,\al}$ le système local qui s'en déduit sur  $M_n^e(C,D)=M_n^e(C_0,D_0)\times_{\Fq}k$. On note $f_{C,D}$ le morphisme de Hitchin sur $M_n^e(C,D)$. Les autres notations sont celles du théorème   \ref{thm:comptage det fix}. Le théorème   \ref{thm:comptage det fix} résulte alors de la combinaison du lemme \ref{lem:step 1} et de la  proposition  \ref{prop:step 2}.
         
       \begin{lemme}\label{lem:step 1}
        Pour tout entier $m\geq 1$,  tout  $\al\in P_0^1(\Fq)$ et tout $\be\in P_0^e(\Fq)$, on a 
    \begin{align*}
    |M_n^\be(C_0,D_0)(\Fqm)|&=  \frac1{|J_0(\Fqm)|} \sum_{\chi\in Q} \chi_\al(\be)^{-m} \sum_{a\in A_0(\Fqm)}\trace(\Fr^{m}, (Rf_{C,D,*} \lc_{\chi,\al})_a) .
    \end{align*}
 
  \end{lemme}

  \begin{preuve}     Soit  $\delta_n\in J_0(\Fqm)$ et  $\delta=\delta_n^{\otimes n}$. Alors  l'application $(\ec,\theta) \mapsto (\ec\otimes \delta_n,\theta)$ induit un isomorphisme de $M_n^\be(C,D)$ sur $M_n^{\be\otimes \delta}(C,D) $ qui se descend à $\Fqm$. On en déduit une bijection entre les ensembles correspondants de $\Fqm$-points: on a donc
    \begin{align}\label{eq:be'}
      |M_n^\be(C_0,D_0)(\Fqm)|=|M_n^{\be\otimes \delta}(C_0,D_0)(\Fqm)|.
    \end{align}
    Par inversion de Fourier, on a 
    \begin{align*}
      \sum_{\delta\in J_0(\Fqm)^{\otimes n}} |M_n^{\be\otimes \delta}(C_0,D_0)(\Fqm)|&=\frac1{|Q_m|} \sum_{\chi\in Q_m}\sum_{(\ec,\theta)} \chi(\det(\ec)\be^{-1})
    \end{align*}
    où, à droite, la somme est prise sur les fibrés de Hitchin $(\ec,\theta) $ dans $  M_n^e(C_0,D_0)(\Fqm)$.
  
  En tenant compte de \eqref{eq:be'}, on aboutit à
  \begin{align*}
    |M_n^\be(C_0,D_0)(\Fqm)|&= \frac1{|J_0(\Fqm)^{\otimes n}|} \sum_{\delta\in J_0(\Fqm)^{\otimes n}} |M_n^{\be\otimes \delta}(C_0,D_0)(\Fqm)|\\
    &= \frac1{|J_0(\Fqm)|} \sum_{\chi\in Q_m}\sum_{(\ec,\theta)} \chi(\det(\ec)\be^{-1})
  \end{align*}
  où la somme sur $(\ec,\theta)$ est comme ci-dessus. En utilisant le lemme \ref{lem:ntorsion}, on en déduit
  \begin{align*}
    |M_n^\be(C_0,D_0)(\Fqm)|&=\frac1{|J_0(\Fqm)|} \sum_{\chi\in Q}\sum_{(\ec,\theta)} \chi(N_m(\det(\ec)\be^{-1}))\\
    & = \frac1{|J_0(\Fqm)|}    \sum_{\chi\in Q}\chi_\al(\be)^{-m} \sum_{(\ec,\theta)} \chi_\al(N_m(\det(\ec)) .
  \end{align*}
 Par la formule des traces de Grothendieck-Lefschetz,  pour tout $a\in A_0(\Fqm)$, on a
  \begin{align*}
\trace(\Fr^{m}, (Rf_{C,D,*} \lc_{\chi,\al})_a)=\sum_{(\ec,\theta)\in f_{C,D}^{-1}(a)(\Fqm)} \chi_\al(N_m(\det(\ec)).
  \end{align*}
  Cela permet de conclure.

\end{preuve}

 \begin{proposition}  \label{prop:step 2}  Sous les hypothèses du théorème \ref{thm:comptage det fix}, pour tout entier $m\geq 1$ et   tout $\chi\in Q$  on a 
    \begin{align*}
 \sum_{ a\in A_0(\Fqm)} \trace(\Fr^{m}, (Rf_{C,D,*} \lc_{\chi,\al})_a) =(q^{n_\chi^2d_\chi (d_\chi-1)\deg(D)/2} \chi_{\al}(\varpi_D)^{-n(n-1)/2})^m \psi_m(\mathbf{H}_{g,n,p,d_\chi,\tilde e})[\chi,\al]
    \end{align*}
    où $n=n_\chi d_\chi$ et $\tilde e=e+\frac{n(n-1)}{2}\deg(D)$.
  \end{proposition}
  
 \begin{preuve}Soit  $f'_0:M'_0\to A'_0$ la fibration de Hitchin relative à la courbe $C_{\chi,\al}$, au diviseur $D_{\chi,\al}=\rho_{\chi,\al}^*D_0$, au rang $n_\chi=n/d_\chi$ et au degré $e$. Le morphisme $f_0'$ se descend en un morphisme  $\bar f_0': M'_0/\Gamma_\chi\to A_0$. Sur $M'_0/\Gamma_\chi$, on dispose du système local $\mu_{\chi^{\tilde e}}$, cf. § \ref{S:mu chi}.
    On pose $r=   n_\chi^2 d_\chi(d_\chi-1)\deg(D)/2$. On applique le théorème \ref{thm:iso-fx} et la formule des traces de Grothendieck; ainsi, pour tout $a\in A_0(\Fqm)$, on a
      \begin{align*}
        \trace(\Fr^{m}, (Rf_{C,D,*} \lc_{\chi,\al})_a)&= q^{mr} \chi_{\al}(\varpi_D)^{-mn(n-1)/2}   \trace(\Fr^{m}, (R\bar f_*' \mu_{\chi^{\tilde e}})_a).
      \end{align*}
      Par conséquent, on obtient:
      \begin{align*}
     &   \sum_{ a\in A_0(\Fqm)} \trace(\Fr^{m}, (Rf_{C,D,*} \lc_{\chi,\al})_a)\\
     &= q^{mr} \chi_{\al}(\varpi_D)^{-mn(n-1)/2} \frac1{|\Ga_\chi|}\sum_{\gamma \in \Ga_\chi}    \chi(\gamma)^{\tilde e}|(M_{n_\chi}^{e}(C_{\chi,\al},  D_{\chi,\al})(k))^{\gamma\Fr^m} |\\
        &=q^{mr} \chi_{\al}(\varpi_D)^{-mn(n-1)/2} \frac1{|\Ga_\chi|}\sum_{\gamma \in J_0(\Fqm)/\Ker(\chi_m)}    \chi_m(\gamma)^{\tilde e}|(M_{n_\chi}^{e}(C_{\chi,\al},  D_{\chi,\al})(k))^{N_m(\ga)\Fr^m} |
      \end{align*}
      où l'on pose $\chi_m=  \chi\circ N_m$. 
      Soit $\ga\in  J_0(\Fqm)$. On a défini un revêtement  $\rho_{\chi,\gamma \al,m}:C_{\chi,\gamma \al,m}\to C_0\times_{\Fq}\Fqm $ au § \ref{S:variante}. 
      
      Soit $D_{m}$ le tiré en arrière de $D_0$ par  $C_0\times_{\Fq}\Fqm \to C_0$ et $D_{\chi,\gamma  \al,m}= \rho_{\chi,\gamma \al,m} ^* D_{m}$. 

      On a alors:
            \begin{align*}
        |(M_{n_\chi}^{e}(C_{\chi,\al},  D_{\chi,\al})(k))^{N_m(\gamma)\Fr^m} |=|M_{n_\chi}^{e}(C_{\chi,\gamma\al,m},  D_{\chi,\gamma \al,m})(\Fqm)|.
            \end{align*}
            Soit $g_\chi=d_\chi(g-1)+1$. Avec les notations du § \ref{S:vp Frob},  les valeurs propres de $\Fr^m$ agissant sur $H^1( C, \fc_{\chi^i,\gamma \al})$ sont données par
 \begin{align*}
\bullet & \la_1^m,\ldots,\la_g^m,q^m\la_1^{-m},\ldots, q^m\la_g^{-m} \text{ pour }  i=0 \ ;\\
\bullet &\chi_m^i(\gamma) \la_{(2i-1)g+3-2i}^m,\ldots,\chi_m^i(\gamma)\la_{(2i+1)g  -2i}^m \text{ pour } 1\leq i <d_\chi/2\ ;\\
   \bullet &\chi_m^{d_\chi/2}(\gamma)\la_{(d_\chi-1)(g-1)+2}^m,\ldots,\chi_m^{d_\chi/2}(\gamma)\la_{g_\chi}^m \text{ et }\\
   & \chi_m^{d_\chi/2}(\gamma)q^m\la_{(d_\chi-1)(g-1)+2}^{-m},\ldots,\chi_m^{d_\chi/2}(\gamma)q^m\la_{g_\chi}^{-m} \text {  si  }  d_\chi \text{ est pair et } i=d_\chi/2.
 \end{align*}
 Rappelons qu'on a
 \begin{align*}
    H^1( C_{\chi,\gamma \al,m}, \Qlb)=\bigoplus_{i=0}^{d-1}   H^1( C, \fc_{\chi^i,\gamma \al}).
\end{align*}
D'après le théorème \ref{thm:MO}, ce dernier comptage est donné par l'évaluation du polynôme $\mathbf{H}_{g_\chi,n_\chi,d_\chi p}$ de $\ZZ[x_1^{\pm 1},\ldots,x_{g_\chi}^{\pm 1},z]$ en $z=q^m$, $x_i=\la_i^m $ pour $1\leq i \leq g$, $x_{j}=\chi_m^i(\gamma) \la_j$ pour $(2i-1)g+3-2i \leq j \leq (2i+1)g  -2i$ pour $ 1\leq i <d_\chi/2$ et $x_j= \chi_m^{d_\chi/2}(\gamma)\la_j^m$ pour $(d_\chi-1)(g-1)+2\leq j \leq g_\chi$ si  $d_\chi $ est pair. C'est encore l'évaluation du polynôme $\psi_m(\chi_m(\gamma)\cdot \mathbf{H}_{g_\chi,n_\chi,d_\chi p})$ en $z=q$ et $x_i=\la_i$ pour tout $1\leq i \leq g_\chi$ (l'action de $\mu_{d_\chi}$, notée $\cdot$, est celle définie au § \ref{S:action mud}). Par définition de $\mathbf{H}_{g,n,p,d_\chi,\tilde e}$, cf. § \ref{S:action mud}, on a:
 \begin{align*}
  \frac1{|\Ga_\chi|}\sum_{\gamma \in J_0(\Fqm)/\Ker(\chi_m)}    \chi_m(\gamma)^{\tilde e}  (\chi_m(\gamma)\cdot \mathbf{H}_{g_\chi,n_\chi,d_\chi p})=  \mathbf{H}_{g,n,p,d_\chi,\tilde e}.
 \end{align*}
 Le résultat s'ensuit.
    \end{preuve}

  \end{paragr}

  \begin{paragr}[Variante.] --- Résumons le principe de la démonstration du théorème \ref{thm:comptage det fix}: on exprime le comptage cherché, à l'aide de la transformation de Fourier du lemme \ref{lem:step 1} et du théorème \ref{thm:iso-fx},  en terme du comptage de fibrés de Hitchin pour certains revêtements afin de pouvoir appliquer le théorème de Mozgovoy-O'Gorman. Dans ce §, on reformule ce qu'on obtient sans user du théorème de Mozgovoy-O'Gorman. Pour cela, on introduit $Q_m^1$ l'ensemble des caractères $\chi$ de $P_0(\Fqm)$ d'ordre divisant $n$ tels que  $\Ker(\chi)\cap P_0^1(\Fqm)\not=\emptyset$. Soit $\chi\in Q_m^1$ et  $\al\in P_0^1(\Fqm)$ tel que $\chi(\al)=1$. On dispose alors d'un revêtement $\rho_{\chi,\al,m}:C_{\chi,\al,m}\to C_0\times_{\Fq}\Fqm $ qui ne dépend pas, à $\Fqm$-isomorphisme près, du choix de $\al$, cf. § \ref{S:variante}.  On le note ci-dessous simplement $\rho_{\chi,m}:C_{\chi,m}\to C_0\times_{\Fq}\Fqm $. Soit $D_{\chi,m}= \rho_{\chi,m} ^* D_{m}$ où  $D_{m}$ est le tiré en arrière de $D_0$ par  $C_0\times_{\Fq}\Fqm \to C_0$. On note aussi $d_\chi|n$ l'ordre de $\chi$. L'énoncé suivant est alors une reformulation du lemme \ref{lem:step 1} et d'une partie de la preuve de la proposition \ref{prop:step 2}.

    \begin{theoreme}
      \label{thm:comparaison comptage}
        Supposons que la caractéristique de $\Fq$ est $>n$ et vérifie l'inégalité \eqref{eq:hyp-p}. Pour tout entier $m\geq 1$ et tout $\be\in P_0(\Fq)$ de degré premier à $n$, on a
        \begin{align*}
          |M_n^\be(C_0,D_0)(\Fqm)|&=  \frac1{|J_0(\Fqm)|} \sum_{\chi\in Q_m^1} \frac{q^{m r_\chi\deg(D)}}{d_\chi}  \chi(\be \varpi_D^{n(n-1)/2})^{-1}   |M_{n_\chi}^e( C_{\chi,m},  D_{\chi,m})(\Fqm) |
        \end{align*}
avec     $n_\chi d_\chi=n$ où  $d_\chi$ est l'ordre de $\chi$  et $r_\chi =n_\chi^2 d_\chi(d_\chi-1)/2$.
    \end{theoreme}
  \end{paragr}

\subsection{Comptage des fibrés de Hitchin de trace nulle et de déterminant fixé}\label{ssec:trace nulle}

\begin{paragr} Soit $e\in \ZZ$ premier à $n$ et  $\be\in P_0^e(\Fq)$. Soit $N_n^\be(C_0,D_0)\subset M_n^\be(C_0,D_0)$ le fermé constitué des fibrés de Hitchin $(\ec,\theta)$ tels que $\trace(\theta)\in H^0(C_0,\oc_{C_0}(D_0))$ est nul. Comme on suppose que la caractéristique de $\Fq$ est $>n$, l'application qui envoie un fibré de Hitchin  $(\ec,\theta)$ et une section $\la\in  H^0(C_0,\oc_{C_0}(D_0))$ sur $(\ec,\theta+\la \Id_{\ec})$ induit un isomorphisme
  \begin{align*}
    N_n^\be(C_0,D_0)\times H^0(C_0,\oc_{C_0}(D_0))\to M_n^\be(C_0,D_0).
  \end{align*}
  Les arguments de \cite{Nitsure} montrent que  $N_n^\be(C_0,D_0)$ est une variété quasi-projective, lisse sur $\Fq$ et purement de dimension
  \begin{align*}
    \dim_{\Fq}(   N_n^\be(C_0,D_0))=(n^2-1)\deg(D).
  \end{align*}
Soit   $N_n^\be(C,D)$ le changement de base à $k$ de  $N_n^\be(C_0,D_0)$.  Il est  connu que, bien que $N_n^\be(C,D)$  ne soit pas propre sur $k$, sa cohomologie à support compact est néanmoins pure: nous ne répéterons pas ici l'argument qui repose sur la propreté du morphisme de Hitchin et l'existence d'une action du groupe multiplicatif qui contracte l'espace  $N_n^\be(C,D)$ sur la fibre centrale.
\end{paragr}

\begin{paragr}[Comptage de $\Fqm$-points.] ---    Le théorème suivant est une conséquence facile du théorème \ref{thm:comptage det fix}. 
  
   \begin{theoreme}\label{thm:comptage tr 0}  Les hypothèses et les notations sont celles du  théorème \ref{thm:comptage det fix}.           Pour tout entier $m\geq 1$ et tout $\be\in P_0(\Fq)$ de degré premier à $n$, on a
      \begin{align*}
    |N_n^\be(C_0,D_0)(\Fqm)|&=   \sum_{\chi\in Q}  \left( q^{r_\chi \deg(D)} \chi_\al^{-1}(\be \varpi_D^{n(n-1)/2} ) \right)^m \psi_m(\tilde{\mathbf{H}}_{g,n,p,d_\chi,\tilde e})[\chi,\al] .
    \end{align*}       
  \end{theoreme}

  \begin{preuve}
    On a
    \begin{align*}
      |N_n^\be(C_0,D_0)(\Fqm)|& = |H^0(C_0\times_{\Fq}\Fqm,\oc_{C_0\times_{\Fq}\Fqm}(D_0))|^{-1}\cdot  |M_n^\be(C_0,D_0)(\Fqm)|\\
      &= \frac{|J_0(\Fqm)|\cdot  |M_n^\be(C_0,D_0)(\Fqm)|}{    q^{m(p+g-1)}\prod_{i=1}^g (1-\la_i^m)(1-q^m\la_i^{-m}) }. 
    \end{align*}
    Le théorème résulte alors du  théorème \ref{thm:comptage det fix} et de l'égalité
    \begin{align*}
    \frac{  \psi_m(\mathbf{H}_{g,n,p,d_\chi,\tilde e})[\chi,\al] }{    q^{m(p+g-1)}\prod_{i=1}^g (1-\la_i^m)(1-q^m\la_i^{-m}) }=\psi_m(\tilde{\mathbf{H}}_{g,n,p,d_\chi,\tilde e})[\chi,\al] .
    \end{align*}
  \end{preuve}
\end{paragr}

\subsection{Polynôme de Poincaré et caractérisque d'Euler}\label{ssec:Poincare}

\begin{paragr}[Polynôme de Poincaré.] --- Pour tout entier $d$, soit $\psi_g(d)$ le nombre d'éléments de $(\ZZ/d\ZZ)^{2g}$ d'ordre $d$. Pour tout multiple $n$ de $d$, on voit que $\psi_g(d)$ est aussi le nombre d'éléments de $(\ZZ/n\ZZ)^{2g}$ d'ordre $d$. On en déduit qu'on a
  \begin{align*}
    \psi_g(n)=\sum_{d|n} \mu(d) (n/d)^{2g}
  \end{align*}
  où $\mu$ est la fonction de Möbius. Le théorème \ref{thm:comptage tr 0} admet le corollaire suivant qui  utilise la pureté de la cohomologie.
  
\begin{corollaire}\label{cor:Poincare} Sous les hypothèses du théorème \ref{thm:comptage tr 0}, pour tout $\be\in P_0(\Fq)$ de degré premier à $n$, on a l'égalité suivante entre polynômes en $u$:

  \begin{align*}
    \sum_{i=0}^{2(n^2-1)\deg(D)} \dim_{\Qlb}(H^i_c(N_n^\be(C,D),\Qlb)) u^i= \sum_{d|n} \psi_g(d) u^{   \deg(D)(d-1)n^2/d} \tilde{\mathbf{H}}^\flat_{g,n,p,d,\tilde e}(-u)
  \end{align*}
  où
  \begin{itemize}
  \item $\tilde e=\deg(\be)+\frac{n(n-1)}{2}\deg(D)$;
  \item  $p=\deg(D)-(2g-2)$;
  \item  $\tilde{\mathbf{H}}^\flat_{g,n,p,d,\tilde e}$ est défini au § \ref{S:specialisation}.
  \end{itemize}
\end{corollaire}

  \begin{preuve}
    Il résulte de la pureté de la cohomologie de $N_n^\be(C,D)$ que les valeurs propres de Frobenius agissant sur $H^i_c(N_n^\be(C,D),\Qlb)$ sont pures de poids $i$ c'est-à-dire que ce sont des entiers algébriques $\la\in \Qlb$ qui vérifient $|\iota(\la)|=q^{i/2}$ pour tout  plongement de $\Qlb$ dans $\CC$. En utilisant la formule des traces de Grothendieck-Lefschetz, on voit que  $|N_n^\be(C_0,D_0)(\Fqm)|$ est la somme sur les entiers $i\geq 0$ de $(-1)^i$ fois la somme des valeurs propres de poids $i$ élevées à la puissance $m$. Mais cette somme est donnée par le membre de droite de l'égalité du théorème \ref{thm:comptage tr 0}. Dans cette expression, pour tout $\chi\in Q$, on voit que $\chi_\al^{-1}(\be \varpi_D^{n(n-1)/2} )$, qui  est une racine de l'unité,  est de poids $0$. Il s'ensuit que $q^{n_\chi^2d_\chi (d_\chi-1)\deg(D)/2}  \chi_\al^{-1}(\be \varpi_D^{n(n-1)/2}) $ est de poids $n_\chi^2d_\chi (d_\chi-1)\deg(D)$. Notons que $d_\chi (d_\chi-1)$ est pair de même que $n_\chi^2d_\chi (d_\chi-1)\deg(D)$.  Par ailleurs les valeurs propres de Frobenius qui interviennent dans l'évaluation du polynôme $\psi_m(\tilde{\mathbf{H}}_{g,n,p,d_\chi,\tilde e})$ sont les valeurs de Frobenius agissant sur $H^1( C_{\chi,\gamma \al,m}, \Qlb)$ et sont donc pures de poids $1$. Il s'ensuit que le polynôme de Poincaré de $H^\bullet_c(N_n^\be(C,D),\Qlb)$ est donné par le polynôme en la variable $u$ suivant:
    \begin{align*}
       \sum_{\chi\in Q}   u^{n_\chi^2d_\chi (d_\chi-1)\deg(D)} \tilde{\mathbf{H}}_{g,n,p,d_\chi,\tilde e}(-u,\ldots,-u, u^2).
    \end{align*}
    Suivant les définitions de \ref{S:specialisation}, on a
    \begin{align*}
      \tilde{\mathbf{H}}_{g,n,p,d_\chi,\tilde e}(-u,\ldots,-u, u^2)= \tilde{\mathbf{H}}_{g,n,p,d_\chi,\tilde e}^\flat(-u)
    \end{align*}
    et ce polynôme ne dépend que de l'ordre $d_\chi$ de $\chi$.
On en déduit le résultat vu que le groupe $Q$ est isomorphe à $(\ZZ/n\ZZ)^{2g}$ par le lemme \ref{lem:ntorsion}.
  \end{preuve}
\end{paragr}

\begin{paragr}[Caractéristique d'Euler-Poincaré.] ---

  \begin{corollaire}\label{cor:Euler}   Soit $\be\in P_0(\Fq)$ de degré premier à $n$.    Si $\deg(D)$ est impair et si  $n\equiv 2 \mod 4 $, la caractéristique d'Euler-Poincaré de $N_n^\be(C,D)$ vaut (sous les hypothèses du théorème \ref{thm:comptage tr 0}):
    \begin{itemize}
    \item $5$ pour $g=1$;
    \item  $-\mu(n) n^{4g-3}$  pour $g\geq 2$.
    \end{itemize}
    Dans les autres cas, elle vaut
    \begin{itemize}
    \item $1$ pour $g=1$;
    \item  $\mu(n) n^{4g-3}$  pour $g\geq 2$.
    \end{itemize}
  \end{corollaire}

  \begin{remarque}
    Si $\deg(D)$ est pair, la caractéristique d'Euler-Poincaré de $N_n^\be(C,D)$ est donc égale à  $\mu(n) n^{4g-3}$ si $g\geq 2$ et $1$ si $g=1$. Lorsque le corps de base est $\CC$, c'est précisément la caractéristique d'Euler-Poincaré de  l'espace des fibrés de Hitchin de rang $n$, relatifs  au faisceau canonique, de déterminant  égal à un fibré en droites fixé de degré premier au rang $n$.   Ce dernier calcul  est dû à Mereb, \cite[second corollaire p.859]{Mereb}. En fait, dans cet article, Mereb calcule la caractéristique d'Euler-Poincaré d'une variété homéomorphe, à savoir une certaine variété caractère associée au groupe $\SL(n)$. L'homéomorphisme en question provient de la théorie de Hodge non-abélienne.
      \end{remarque}
    
          \begin{preuve}    La caractéristique d'Euler-Poincaré de  $N_n^\be(C,D)$ est donnée par l'évaluation en $u=-1$ du polynôme de Poincaré c'est-à-dire, d'après le corollaire \ref{cor:Poincare}, par l'expression
    \begin{align}\label{eq:X EP}
      \sum_{d|n} \psi_g(d)\tilde{\mathbf{H}}_{g,n,p,d,\tilde e}^\flat(1)
    \end{align}
    avec $\tilde e=\deg(\be)+\frac{n(n-1)}{2}\deg(D)$ et $p=\deg(D)-(2g-2)$.
    
    Supposons tout d'abord $g=1$. On a donc $p=\deg(D)$. Supposons $\deg(D)$ impair et $n\equiv 2 \mod 4 $. D'après le lemme \ref{lem:calcul H11 g1} l'expression \eqref{eq:X EP} ci-dessus vaut
    \begin{align*}
      2\psi_1(1)+\psi_1(2)=2+2^{2}-1=5.
    \end{align*}
Dans tous les autres cas, l'expression \eqref{eq:X EP} est égale à $1$ d'après  le lemme \ref{lem:calcul H11 g1}.

Supposons désormais  $g\geq 2$. Posons $\eps=-1$ si $\deg(D)$ est impair et  $n \equiv 2 \mod 4$ et $\eps=1$ sinon.  La combinaison des lemmes \ref{lem:calcul H11} et \ref{lem:last} donne le calcul suivant pour l'expression \eqref{eq:X EP}
    \begin{align*}
      (-1)^{pn }  n^{2g-3} \mu(n)  \eps \sum_{d|n} \psi_g(d) ((-1)^p\eps)^{d}.
    \end{align*}
    Supposons $\deg(D)$ pair. Dans ce cas, on a $(-1)^p=1$ et $\eps=1$. On achève le calcul à l'aide de l'identité:
          \begin{align*}
    \sum_{d|n} \psi_g(d) =n^{2g}.
          \end{align*}
          Supposons $\deg(D)$ impair.  Alors $p$ est impair. Si  $n \equiv 2 \mod 4$, on a $(-1)^p\eps=1$ et on conclut encore par l'identité ci-dessus. Si $n$ est impair, on a $\eps=1$ et $d|n$ est aussi impair. On conclut toujours  par l'identité ci-dessus.  On suppose  $4|n$. On a alors $\mu(n)=0$ ce qui conclut.

\end{preuve}

\begin{lemme} \label{lem:calcul H11 g1} (Cas $g=1$). Soit $d|n$. Soit $\tilde e\in \ZZ$ tel que  $\tilde e-\frac{n(n-1)}{2}\deg(D)$ est premier à $n$. 
  \begin{align*}
    &  \bullet   \tilde{\mathbf{H}}_{1,n,\deg(D),d,\tilde e}^\flat(1)=0 \text{ si }  d\notin\{1,2\};\\
      &   \bullet  \tilde{\mathbf{H}}_{1,n,\deg(D),1,\tilde e}^\flat(1)= \left\lbrace
        \begin{array}{l}
          2 \text{ si } \deg(D) \equiv 1 \mod 2 \text{ et }  n \equiv 2\mod 4 \ ;\\
          1 \text{ sinon ;}
        \end{array}
    \right.\\
   & \bullet    \tilde{\mathbf{H}}_{1,n,\deg(D),2,\tilde e}^\flat(1)= \left\lbrace
        \begin{array}{l}
          1 \text{ si } \deg(D) \equiv 1 \mod 2 \text{ et }  n \equiv 2\mod 4 \ ;\\
          0 \text{ sinon.}
        \end{array}
    \right.
  \end{align*}
    \end{lemme}

    \begin{preuve} On pose $n'=n/d$ et $p=\deg(D)$. En utilisant l'égalité \eqref{eq:identite clef} et le lemme \ref{lem:calcul H sur 1-u} 2.(b) (pour $g=1$), on obtient
      \begin{align*}
       \tilde{\mathbf{H}}_{1,n,p,d,\tilde e}^\flat(1)= \kappa_d \frac1d \sum_{\zeta\in \mu_d}\zeta^{\tilde e}
      \end{align*}
      où $\kappa_d$ vaut $2$ si $p d$ est impair et si $n'$ est congru à $2$ modulo $4$ et $1$ dans tous les autres cas. La somme ci-dessus est nulle sauf si $d|\tilde e$ c'est-à-dire si  $d$ divise le pgcd de $n$ et $\tilde e$.   Or ce pgcd vaut $1$ ou $2$, cf. lemme \ref{lem:last}.

      Pour $d\notin\{1,2\}$ on a donc $\tilde{\mathbf{H}}^\flat_{1,n,p,d,\tilde e}(1)=0$.

      Pour $d=1$, on a $n'=n$ et $\tilde{\mathbf{H}}^\flat_{1,n,p,1,\tilde e}(1)=\kappa_1$; ce dernier  vaut $2$ si $\deg(D)$ est impair et $n$ est congru à $2$ modulo $4$. Il vaut $1$ sinon.

     Considérons enfin $d=2$; on a donc  $\kappa_2=1$. On trouve $\tilde{\mathbf{H}}^\flat_{1,n,p,2,\tilde e}(1)=0$ sauf si  le pgcd de $n$ et $\tilde e$ est $2$ c'est-à-dire, d'après le lemme \ref{lem:last}, si $\deg(D)$ impair et $n$ est congru à $2$ modulo $4$. 
      
    \end{preuve}

    \begin{lemme} \label{lem:calcul H11} On suppose $g\geq 2$. Soit $p\geq 0$ et $\tilde e\in \ZZ$. Soit $d|n$. 
      \begin{enumerate}
         \item Si $d$ et $n/d$ ne sont pas premiers entre eux, on a
        \begin{align*}
          \tilde{\mathbf{H}}^\flat_{g,n,p,d,\tilde e}(1)=0\ ;
        \end{align*}
      \item Si $d$ et $n/d$ sont premiers entre eux, on a
        \begin{align*}
          \tilde{\mathbf{H}}_{g,n,p,d,\tilde e}^\flat(1)= (-1)^{p(n-d) }  n^{2g-3}  \sum_{j|(d,\tilde e)}j \mu(\frac{n}{j})
        \end{align*}
        où $(d,\tilde e)$ désigne le pgcd de $d$ et $\tilde e$.
           \end{enumerate}
    \end{lemme}

    \begin{preuve} Soit $\zeta\in \mu_d$ et $r=n/d$. D'après le lemme \ref{lem:calcul H sur 1-u} 2 .(a) et l'identité \ref{eq:identite clef},  on a 
    \begin{align}\label{eq:somme mu}
 \tilde{\mathbf{H}}_{g,n,p,d,\tilde e}^\flat(1)=\frac1d \sum_{\zeta\in \mu_d} \zeta^{\tilde e}      (-1)^{(n-d)p}\mu(r) r^{2g-3}     \prod_{i=1}^{d-1}(1-\zeta^{ir})^{g-1} (1-\zeta^{-ir})^{g-1} .
      \end{align}

      \begin{lemme} \label{lem:petit calcul}Soit $\zeta\in \mu_d$. On a
            \begin{align*}
           \prod_{i=1}^{d-1}(1-\zeta^{ir})=d
            \end{align*}
            si  $r$ et $d$ sont premiers entre eux et $\zeta$ est d'ordre $d$. Dans tous les autres cas, on a
              \begin{align*}
           \prod_{i=1}^{d-1}(1-\zeta^{ir})=0.
            \end{align*}
      \end{lemme}

      \begin{preuve} Si $\zeta^r$ n'est pas d'ordre $d$, on obtient clairement $0$. Si  $\zeta^r$ est d'ordre $d$, alors il en est de même de $\zeta$ et $r$ et $d$ sont premiers entre eux. Par conséquent, $\zeta^{ir}$ parcourt $\mu_d\setminus\{1\}$ pour $1\leq i\leq d-1$. On conclut à l'aide de l'identité:
\begin{align*}
  \prod_{\zeta\in \mu_d,\zeta\not=1} (1-\zeta)=d.
\end{align*}
\end{preuve}

Si $d$ et $r$ ne sont pas premiers entre eux, l'égalité \eqref{eq:somme mu} et le lemme \ref{lem:petit calcul} donne la nullité dans l'assertion 1.
 On suppose désormais que $d$ et $r=n/d$ sont premiers entre eux. Soit $\mu_{d}^*\subset \mu_d$ le sous-ensemble des racines primitives $d$-ièmes de l'unité. En utilisant l'égalité \eqref{eq:somme mu} et le lemme \ref{lem:petit calcul}, puis le lemme \ref{lem:calcul2} ci-dessous, on obtient 
\begin{align*}
  \tilde{\mathbf{H}}_{g,n,p,d_\chi,\tilde e}^\flat(1)&=(-1)^{p(n-d) } \mu(\frac{n}{d}) n^{2g-3}  \left(\sum_{\zeta \in \mu_d^*} \zeta^{\tilde e} \right)\\
  &=   (-1)^{p(n-d) }  n^{2g-3}  \sum_{j|(d,\tilde e)}j \mu(\frac{n}{j}).
\end{align*}

\begin{lemme}\label{lem:calcul2}
  Soit $d\geq 1$ un entier et $i\in \ZZ$. Soit $(d,i)$ le pgcd de $d$ et $i$. On a
    \begin{align*}
            \sum_{\zeta \in \mu_{d}^*} \zeta^{i} =\sum_{j|(d,i)}j \mu(\frac{d}{j}).
    \end{align*}
\end{lemme}

\begin{preuve} Soit $\delta_i$ la fonction sur $\NN^*$ caractéristique des diviseurs de $i$. Pour tout  $d\geq 1$ entier, soit
  \begin{align*}
    \rho_i(d)=  \sum_{\zeta \in \mu_{d}^*} \zeta^{i}.
  \end{align*}
  On a
  \begin{align*}
    d \delta_i(d)&=\sum_{\zeta\in \mu_d} \zeta^i\\
    &=\sum_{j|d}  \rho_i(j).
  \end{align*}
Par inversion de Möbius, on a  $\rho_i(d)=\sum_{j|d}    j\delta_i(j) \mu(\frac{d}{j})=\sum_{j|(d,i)} j \mu(\frac{d}{j})$.
  
\end{preuve}

\end{preuve}

\begin{lemme}\label{lem:last}
Soit $n\geq 1$ un entier, $e\in \ZZ$ premier à $n$ et $\tilde e= e+n(n-1)\deg(D)/2$.
  \begin{enumerate}
  \item Si $\deg(D)$ est impair et $n \equiv 2 \mod 4$,  alors le pgcd de $\tilde e$ et $n$ vaut $2$ et pour tout diviseur $d|n$ on a
    \begin{align*}
          \sum_{j|(d,\tilde e)}j \mu(\frac{n}{j})=(-1)^{d-1}\mu(n).
    \end{align*}
  \item Dans tous les autres cas, le pgcd de $\tilde e$ et $n$  vaut $1$ et pour tout $d|n$ on a 
    \begin{align*}
          \sum_{j|(d,\tilde e)}j \mu(\frac{n}{j})=\mu(n).
    \end{align*}
  \end{enumerate}
\end{lemme}

\begin{preuve}
  Soit $f$ un diviseur commun de $\tilde e$ et $n$.

  Si $n$ est impair, on a $f|n(n-1)\deg(D)/2$ donc $f|e$ donc $f=1$. Donc $(n,\tilde e)=1$.

  Supposons $n$ pair et écrivons $n=2n'$. Si $f|n'$ alors $f| n(n-1)\deg(D)/2$ et $f|e$ et de nouveau $f=1$. Il s'ensuit que le pgcd  $(n,\tilde e)$ vaut $1$ ou $2$.

  Si  $(n,\tilde e)=1$, l'égalité de l'assertion 2 est évidente.

  On suppose désormais $(n,\tilde e)=2$. On a donc $n$ et $\tilde e$ pairs, $n'$ et $e$ impairs. On a donc aussi $\deg(D)$ impair. On observe que $n$ est congru à $2$ modulo 4 et $\mu(n)=-\mu(n')$. Réciproquement si $\deg(D)$ impair et que  $n$ est congru à $2$ modulo $4$, on voit que $\tilde e$ est pair et donc   $(n,\tilde e)=2$.

  Soit $d$ un diviseur de $n$. Si $d$ est impair, on a $(d,\tilde e)=1$ et 
  \begin{align*}
          \sum_{j|(d,\tilde e)}j \mu(\frac{n}{j})=\mu(n).
    \end{align*}
    Si   $d$ est pair, on a $(d,\tilde e)=2$ et
      \begin{align*}
          \sum_{j|(d,\tilde e)}j \mu(\frac{n}{j})=    \mu(n) +2\mu(n')=-\mu(n).
    \end{align*}
\end{preuve}
\end{paragr}
\bibliography{biblio}
\bibliographystyle{alpha}

\begin{flushleft}
Pierre-Henri Chaudouard \\
Université Paris Cité\\
IMJ-PRG \\
Bâtiment Sophie Germain\\
8 place Aurélie Nemours\\
F-75013 PARIS 
France\\
\medskip
Institut universitaire de France (IUF)\\
\medskip
email:\\
Pierre-Henri.Chaudouard@imj-prg.fr \\
\end{flushleft}

\end{document}